\documentclass[final,3p,times]{elsarticle}

\biboptions{sort&compress}

\usepackage{txfonts}
\usepackage{fleqn}

\usepackage{amsthm}
\usepackage{bm}
\usepackage{multirow}
\usepackage{array}
\usepackage{xcolor}
\colorlet{myrem}{black}

\usepackage[%
pdftex,%
colorlinks=true,%
bookmarks=true,%
citecolor=blue,%
urlcolor=blue]{hyperref} 
\usepackage[caption=false]{subfig}

\usepackage{amsmath}
\usepackage{mathrsfs}  

\newcommand{\field}[1]{\mathbb{#1}}
\newcommand{\OP}[1]{\mathscr{#1}}

\newcommand{\tp}{\intercal}
\DeclareMathOperator*{\supp}{supp}
\let\Re\relax
\DeclareMathOperator{\Re}{Re}
\let\Im\relax
\DeclareMathOperator{\Im}{Im}

\DeclareMathOperator{\Span}{Span}

\graphicspath{{./figs_cmpt/}}
\DeclareGraphicsExtensions{.pdf,.jpeg,.png}
\newcommand{\wtilde}[1]{\widetilde{#1}}

\newcommand{\fs}[1]{\mathsf{#1}}

\DeclareMathOperator{\diag}{diag}

\newcommand{\ovl}[1]{\overline{#1}}

\newcommand{\bigO}[1]{\mathop{\mathscr{O}}\left(#1\right)}
\newcommand{\melem}[1]{\mathfrak{#1}}

\let\Re\relax
\DeclareMathOperator{\Re}{Re}
\let\Im\relax
\DeclareMathOperator{\Im}{Im}
\newcommand{\vv}[1]{\boldsymbol{#1}}
\newcommand{\vs}[1]{\boldsymbol{#1}}

\newcommand{\what}[1]{\widehat{#1}}

\newtheorem{theorem}{Theorem}

\newtheorem{defn}{Definition}
\newtheorem{rem}{Remark}
\usepackage{enumitem}
\journal{Journal of Computational Physics}
\begin{document}

\begin{frontmatter}



\title{Transparent boundary condition and its effectively local approximation 
for the Schr\"{o}dinger equation on a rectangular computational domain}


\author[phy]{Samardhi Yadav}
\ead{Samardhi@physics.iitd.ac.in}
\author[phy,opc]{Vishal Vaibhav}
\ead{vishal.vaibhav@gmail.com}

\address[phy]{Department of Physics, Indian Institute of Technology Delhi, Hauz
Khas, New Delhi–110016, India}%
\address[opc]{Optics and Photonics Center, Indian Institute of Technology Delhi,
Hauz Khas, New Delhi–110016, India}%
\begin{abstract}
The transparent boundary condition for the free Schr\"{o}dinger equation 
on a rectangular computational domain requires implementation of an operator of the
form $\sqrt{\partial_t-i\triangle_{\Gamma}}$ where $\triangle_{\Gamma}$ is the 
Laplace-Beltrami operator. It is known that this operator is nonlocal in 
time as well as space which poses a significant challenge in developing an
efficient numerical method of solution. The computational complexity of the
existing methods scale with the number of time-steps which can be attributed to
the nonlocal nature of the boundary operator. In this work, we report an 
effectively local approximation for the boundary operator such that the
resulting complexity remains independent of number of time-steps. At the heart
of this algorithm is a Pad\'e approximant based rational approximation of
certain fractional operators that handles corners of the domain adequately. 
For the spatial discretization, we use a Legendre-Galerkin spectral method with
a new boundary adapted basis which ensures that the resulting linear system 
is banded. A compatible boundary-lifting procedure is also presented which
accommodates the segments as well as the corners on the boundary. The proposed
novel scheme can be implemented within the framework of any one-step
time marching schemes. In particular, we demonstrate these ideas for two
one-step methods, namely, the backward-differentiation formula of order 1 (BDF1) and the 
trapezoidal rule (TR). For the sake of comparison, we also present a convolution 
quadrature based scheme conforming to the one-step methods which is 
computationally expensive but serves as a golden standard. 
Finally, several numerical tests are presented to 
demonstrate the effectiveness of our novel method as well as to verify the 
order of convergence empirically. 
\end{abstract}


\end{frontmatter}

\tableofcontents
\section*{Notations}
\label{sec:notations}
The set of non-zero positive real numbers ($\field{R}$) is denoted by
$\field{R}_+$. For any complex number $\zeta$, $\Re(\zeta)$ and $\Im(\zeta)$ 
refer to the real and the imaginary parts of $\zeta$, respectively. The open 
interval $(-1,1)$ is denoted by $\field{I}$.
 
\section{Introduction}
This paper is devoted to the numerical solution of the two-dimensional free 
Schr\"{o}dinger equation. The free Schr\"odinger equation and its various 
generalizations are of significant physical interest because they appear in 
many areas of application such as quantum mechanics, propagation of 
electromagnetic waves under paraxial approximation 
(or slow varying envelope approximation for pulse
propagation in optical fibers)~\cite{AK2003} and parabolic approximation to one-way
propagation equation in underwater acoustics~\cite{LM1998}.

The initial value problem (IVP) corresponding to the free 
Schr\"{o}dinger equation formulated on $\field{R}^2$ reads as
\begin{equation}\label{eq:2D-SE}
i\partial_tu+\triangle u=0,\quad(\vv{x},t)\in\field{R}^2\times\field{R}_+,
\end{equation}
where the initial condition $u(\vv{x},0)=u_0(\vv{x})$ is assumed to be compactly 
supported. In order to ensure the uniqueness of the solution, a Sommerfeld-like 
radiation condition at infinity is imposed~\cite{CiCP2008}. For the 
numerical solution of this problem, we first need to choose 
a bounded computational domain, denoted by $\Omega_i$, such that 
$\supp u_0(\vv{x})\subset\Omega_i$. The next step is to determine appropriate 
boundary condition at the fictitious boundary ($\partial\Omega_i$)
such that the solution of the resulting initial-boundary value 
problem (IBVP) on $\Omega_i$ remains identical to that of original problem 
which is formulated on an unbounded domain. Such boundary
conditions are known as the \emph{transparent boundary conditions} (TBCs) which 
usually appear as \emph{Dirichlet-to-Neumann} (DtN) or
\emph{Neumann-to-Dirichlet} (NtD) maps at the fictitious 
boundary~\cite{CiCP2008}. In cases where it is possible to obtain these maps in
a closed form, they turn out to be nonlocal with respect to space 
as well as time. The specific form of
the transparent or nonreflecting boundary maps depend on the nature of the
computational domain chosen. For practical reasons, the rectangular domain happens 
to be a natural choice for the computational domain. In this work, we address the 
challenges involved in choosing $\Omega_i$ to be a rectangular domain for which an 
exact nonreflecting boundary condition is known at a continuous
level~\cite{FP2011,V2019}. In particular, we address the problem of designing a 
numerical scheme which provides a viable trade-off between computational 
complexity and accuracy when it comes to dealing with the nonlocal boundary 
maps at a discrete level. At the outset, let us state that
the new algorithm presented in this work 
circumvents the need to store the entire history of the field on the
boundary $\partial\Omega_i$ unlike the existing methods~\cite{FP2011} such that
the overall complexity turns out to be dominated by that of the linear solver
for the interior problem. This effective localization of the DtN 
maps rely on a \emph{novel} Pad\'e approximant based approach that handles the
corners adequately.


Several authors have contributed towards the 
formulation of exact TBCs for the two-dimensional Schr\"{o}dinger 
equation~\cite{Menza1996,Menza1997,S2002,HH2004,FP2011,JPA2018,V2019}.
Exact TBCs for the free Schr\"odinger equation on convex domains with smooth 
boundary was provided by Sch\"adle~\cite{S2002} in terms of a single and a 
double layer potential. Han and Huang~\cite{HH2004} provided the exact TBCs 
in terms of Hankel functions for the circular domain. {\color{myrem}Antoine~\textit{et~al.} 
provided asymptotic approximations of a microdifferential TBC for the case of linear Schr\"odinger equation 
in $\field{R}^2$~\cite{AB2001}. Szeftel developed absorbing boundary conditions and showed 
well posedness of the overall IBVP for the 
linear Schr\"odinger equation on convex domains with smooth boundary in $\field{R}^d$~\cite{S2005}.} 
For more general cases including the nonlinear Schr\"odinger equation on convex domains,
Antoine~\textit{et~al.} proposed some higher-order absorbing boundary conditions
with the use of pseudo-differential operators~\cite{ABK2012,ABK2013}. The transparent 
boundary operator for the Schr\"odinger equation on a 
rectangular domain has the form $\sqrt{\partial_t-i\triangle_{\Gamma}}$ where 
$\triangle_{\Gamma}$ is the Laplace-Beltrami operator. Note that these operators are 
nonlocal, both in time as well as space. To the best of our knowledge, an integral 
representation of this operator first appeared in~\cite{Menza1996} and more 
recently derived in~\cite{V2014,V2019}. The earliest numerical implementation 
for the operator $\sqrt{\partial_t-i\triangle_{\Gamma}}$ was given by 
Menza~\cite{Menza1997} who introduced a Pad\'e approximant based 
rational approximation for the complex square root function to obtain effectively 
local boundary conditions. However, this method did not address the problems
arising at the corners of the rectangular domain properly\footnote{We discuss these issues in
Sec.~\ref{sec:CT-CP} of this paper}. This issue was first resolved by Feschenko and 
Popov~\cite{FP2011} who derived the exact TBCs for the two-dimensional Schr\"odinger 
equation on a rectangular domain by introducing an auxiliary function which satisfies a 
one-dimensional Schr\"odinger equation along the boundary. As elaborated 
in the work of the second author~\cite{V2019}, the boundary conditions obtained 
by Feschenko and Popov can be directly obtained from the operator 
$\sqrt{\partial_t-i\triangle_{\Gamma}}$ by reducing it to the time-fractional operator 
$\sqrt{\partial_t}$ acting on an auxiliary function with the aforementioned
properties.
The corners are handled by imposing the TBCs for the one dimensional Schr\"odinger
equation at the end points of the segment under consideration. 
Observe that the nonlocality in space is reduced to solving one-dimensional 
Schr\"odinger equations along the boundary segments of the rectangular computational 
domain but the nonlocality with respect to the temporal variable still poses a serious 
challenge. The numerical approach presented by Feschenko and 
Popov becomes computationally expensive on account of the
increasing number of linear systems to be solved on the boundary segments with 
increasing time-steps\footnote{The authors in~\cite{FP2011} suggested truncating the history
of the field at the boundary in the discrete convolution in order to reduce the
complexity, however, no estimates for the effective depth of history are available 
in their work. We have deferred these considerations in this manuscript because
such a proposal needs to be first evaluated for the 1D problem before attempting
the 2D problem.}. The storage requirements also
scale linearly with the number of time-steps which makes it a formidable task to
carryout the numerical solution with smaller temporal grid parameters (assuming a
fixed duration of time). 

To remedy this problem, we introduce an effectively local form of these TBCs by 
inserting a Pad\'e approximant based rational approximation for 
$\sqrt{\partial_t}$ operator acting on the auxiliary function. The resulting
numerical scheme is referred to as the \emph{novel-Pad\'e} (NP) method 
to contrast it with Pad\'e approximant based approach given by Menza~\cite{Menza1997}.
Let us further emphasize that Menza's method involved writing an
effectively local form for the full transparent boundary operator 
($\sqrt{\partial_t-i\triangle_{\Gamma}}$) via rational approximation as opposed
to the temporal fractional operator $\sqrt{\partial_t}$. The NP method ensures
that the computational cost and the memory requirement does not grow with
increasing time-steps. At the discrete level, the TBCs are formulated as a 
Robin-type boundary condition at each time-step. For the spatial discretization,
we use a Legendre-Galerkin spectral method where we develop a new basis,
referred to as the boundary adapted basis~\cite{S1994}, in terms
of the Legendre polynomials which ensures that the resulting linear system within 
the framework of Galerkin method is sparse. This approach, however, requires 
a boundary-lifting procedure in order to homogenize the Robin-type boundary 
conditions. To demonstrate the effectiveness of our method, we first provide 
a convolution quadrature (CQ)~\cite{Lubich1986} based implementation of the 
time-fractional operator $\sqrt{\partial_t}$ followed by the NP method. It must
be noted that the former scheme is nonlocal and suffers from the same problems as
that of Feschenko and Popov~\cite{FP2011} (where the quadrature weights are
derived according to~\cite{BP1991}) while the latter is effectively local.
Further, the CQ scheme considered in this paper ensures that the 
time-discretization of the interior problem remains compatible with that of 
the boundary conditions unlike that of~\cite{BP1991} which was found to suffer
from stability issues~\cite{Mayfield1989}.

Before we conclude this section, we would like to contrast our method with that 
of perfectly matched layers (PMLs) which introduces an additional layer of 
finite-thickness outside $\Omega_i$. 
The method of PMLs introduces inhomogeneous terms which correspond to the 
spatial derivatives in the PDEs~\cite{Zheng2007,AGT2020}.
Therefore, within the spectral-Galerkin method, the resulting linear system
necessarily becomes dense. In contrast, our method leads to a banded linear 
system with TBCs imposed on $\partial\Omega_i$. A thorough comparison with the method of
PMLs would require consideration of more general PDEs (such as presence of 
variable potentials etc.) where spectral methods do not offer any 
competitive advantage (in the realm of non-iterative solvers). Therefore, 
this consideration is beyond the scope of this paper. However, let us remark
that our numerical scheme will remain competitive with that of PMLs if the TBCs
considered in this paper continue to hold for the general system. This 
conclusion is based on the 
observation that the one-dimensional Schr\"odinger equations on the boundary
segments have constant coefficients leading to banded linear systems whose cost
of computation is significantly smaller than that of the interior problem. 


The paper is organized in the 
following manner: Sec.~\ref{sec:tbcs} presents the
derivation of TBCs followed by a discussion of numerically tractable forms of
the TBCs. This section also addresses the development of corner conditions whose
discrete form is useful in developing boundary adapted basis. In 
Sec.~\ref{sec:numerical-implementation}, we discuss the time-discretization of 
TBCs and the numerical solution of the IBVP using Legendre-Galerkin method along
a step-by-step analysis of the computational complexity. We confirm the 
efficiency and stability of our numerical schemes with several 
numerical tests presented in Sec.~\ref{sec:numerical-experiments}. Finally, we 
conclude this paper in Sec.~\ref{sec:conclusion}.
\section{Transparent Boundary Conditions}\label{sec:tbcs}
\begin{figure}[!htbp]
\begin{center}
\def\myscale{1}
\includegraphics[scale=\myscale]{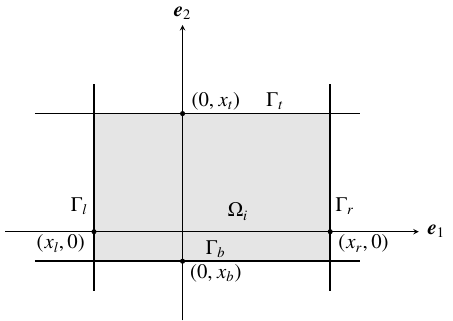}
\end{center}
\caption{\label{fig:rect-domain} The figure shows a rectangular domain
with boundary segments parallel to one of the axes.}
\end{figure}

Consider a rectangular computational domain ($\Omega_i$) with boundary 
segments parallel to one of the axes (see Fig.~\ref{fig:rect-domain}) with
$\Omega_i =(x_l,x_r)\times(x_b,x_t)$ referred to as the~\emph{interior} 
domain. Consider the decomposition of the field $u(\vv{x},t)$ such that
$u(\vv{x},t)\in\fs{L}^2(\field{R}^2)=\fs{L}^2(\Omega_i)\oplus \fs{L}^2(\Omega_e)$ 
where $\Omega_{e}=\field{R}^2\setminus\overline{\Omega}_i$ is referred to as 
the~\emph{exterior} domain. An equivalent formulation of the IVP 
in~\eqref{eq:2D-SE} can be stated as  
\begin{equation}\label{eq:split-IVP}
\begin{split}
\text{interior problem}:\quad
&\left\{\begin{aligned}
&i\partial_tv+\triangle v=0,\quad(\vv{x},t)\in\Omega_i\times\field{R}_+,\\
&v(\vv{x},0)= u_0(\vv{x}), \quad(\vv{x},t)\in\Omega_i,\\
\end{aligned}\right.\\
\text{exterior problem}:\quad
&\left\{\begin{aligned}
&i\partial_tw+\triangle w=0,\quad(\vv{x},t)\in\Omega_e\times\field{R}_+,\\
&w(\vv{x},0)=0,\quad\vv{x}\in\Omega_e,\\
&\lim_{|\vv{x}|\rightarrow\infty}\sqrt{|\vv{x}|}
\left(\nabla
w\cdot\frac{\vv{x}}{|\vv{x}|}+e^{-i\frac{\pi}{4}}\partial_t^{1/2}w\right)=0;
\end{aligned}\right.\\
\text{continuity relations}:\quad 
&v(\vv{x},t)|_{\Gamma}=w(\vv{x},t)|_{\Gamma},\quad 
\partial_{n}v(\vv{x},t)|_{\Gamma}=\partial_{n}w(\vv{x},t)|_{\Gamma},
\end{split}
\end{equation}
where $\Gamma$ denotes the common boundary of the interior and the exterior
domains. {\color{myrem}Note that the initial data $v(\vv{x},0)= u_0(\vv{x})$ is assumed 
to be compactly supported~\footnote{This criteria can not be relaxed for our novel 
Pad\'e approach.} such that $\supp u_0(\vv{x})\subset\Omega_i$.}
For the sake of completeness, we would like to present 
a detailed discussion of the derivation of the TBCs first presented 
in~\cite{V2019}. To this end, let $(\zeta_1,\zeta_2)$ be the covariables 
corresponding to $(x_1,x_2)$ to be used in the two dimensional Fourier 
transform. Let us denote the one dimensional Fourier transform of 
$w(x_1,x_2,t)$ with respect to $x_2$ by $\wtilde{w}(x_1,\zeta_2,t)$:
\begin{equation*}
\begin{split}
&\wtilde{w}(x_1,\zeta_2,t)=\OP{F}_{x_2}[w(x_1,x_2,t)]
\equiv\int_{\field{R}}w(x_1,x_2,t)e^{-i\zeta_2x_2}dx_2.
\end{split}
\end{equation*}
Let us denote the Laplace transform of $\wtilde{w}(x_1,\zeta_2,t)$ with 
respect to $t$ by $\wtilde{W}(x_1,\zeta_2,z)$:
\begin{equation*}
\begin{split}
&\wtilde{W}(x_1,\zeta_2,t)=\OP{L}_{t}[\wtilde{w}(x_1,\zeta_2,t)]
\equiv\int_{\field{R}_+}\wtilde{w}(x_1,\zeta_2,t)e^{-zt}dt.
\end{split}
\end{equation*}
First we focus on deriving the TBCs on the segment $\Gamma_r$ by employing 
the Fourier transform with respect to $x_2$ and then taking the Laplace 
transform of the right exterior problem:
\begin{equation}
(\partial_{x_1}^2+\alpha^2)\wtilde{W}(x_1,\zeta_2,z)=0,\quad x_1\in(x_r,\infty),
\end{equation}
where $\alpha=\sqrt{iz-\zeta_2^2}$ such that $\sqrt{\cdot}$ denotes the branch
with $\Im(\alpha)>0$. The solution takes the form
\begin{equation}
\wtilde{W}({x}_1,\zeta_2,z)
=c\exp(+i\alpha x_1) 
+d\exp(-i\alpha x_1),
\end{equation}
where $c$ and $d$ are independent of $x_1$. In order to have a bounded 
solution on the right exterior domain, $d$ must be zero. Further, $c$ can be
eliminated using the derivative of the field to obtain
\begin{equation}
\partial_{x_1}\wtilde{W}(x_1,\zeta_2,z)=i\alpha\wtilde{W}(x_1,\zeta_2,z).
\end{equation}
To facilitate the inverse Laplace transform of the equation above, we proceed as
follows:
\begin{equation}
\partial_{x_1}\wtilde{w}(x_1,\zeta_2,t)=
\OP{L}^{-1}[\partial_{x_1}\wtilde{W}(x_1,\zeta_2,z)]
=\OP{L}^{-1}[i\alpha^{-1}]\star\OP{L}^{-1}[\alpha^2\wtilde{W}(x_1,\zeta_2,z)],
\end{equation}
where `$\star$' represents the convolution operation with respect $t$. Let us 
observe that $\OP{L}^{-1}\left[{1}/{\sqrt{z}}\right](t)=H(t)/{\sqrt{\pi t}}$ where 
$H(t)$ denotes the Heaviside step-function. In the following, we use 
the fractional operators which are denoted by 
$\partial_t^{\alpha},\;\alpha\in\field{R}$ (see Miller and Ross~\cite{MR1993}).
For $\alpha<0$, we obtain the Riemann-Liouville fractional integral of order 
$|\alpha|$ while $\alpha>0$ corresponds to fractional derivatives which are 
introduced as integral order derivatives of fractional integrals (discussed 
in~\ref{app:frac-op}). Using the shifting property of Laplace transforms and 
introducing fractional operators of order $1/2$, 
we can express the result as
\begin{equation}\label{eq:dtn-map}
\partial_{x_1}\wtilde{w}(x_1,\zeta_2,t)
=e^{i\pi/4}e^{-i \zeta_2^2t}\partial_t^{-1/2}
 e^{i\zeta_2^2t}\left[i\partial_t\wtilde{w}(x_1,\zeta_2,t)
-\zeta_2^2\wtilde{w}(x_1,\zeta_2,t)\right].
\end{equation}
In order to obtain the inverse Fourier transform of the above relations, let us 
introduce the kernel 
\begin{equation}
\mathcal{G}(x_2,t)=\frac{e^{-i\pi/4}}{\sqrt{4\pi t}}
\exp\left[i\frac{x_2^2}{4t}\right],
\end{equation}
defined in a distributional sense with its Fourier transform given by 
$\wtilde{G}(\zeta_2,t)=e^{-i\zeta_2^2 t}$.
Introducing the integral definition for $1/2$-order fractional integral 
and taking inverse Fourier transform of~\eqref{eq:dtn-map} yields
\begin{equation}\label{eq:IPS-TBC}
\partial_{x_1}{w}(\vv{x},t)
=\frac{e^{i\pi/4}}{\sqrt{\pi}}\int_0^t\int_{\field{R}}
\left[i\partial_{\tau}{w}(x_1,x'_2,\tau)+\partial^2_{x'_2}w(x_1,x'_2,\tau)
\right]
\frac{\mathcal{G}(x_2-x'_2,t-\tau)}{\sqrt{t-\tau}}dx_2'd\tau.
\end{equation}
The relationship above is referred to as the~\emph{Dirichlet-to-Neumann} (DtN)
map relating the Dirichlet data ${w}(x_1,x_2,t)$ to the Neumann data 
$\partial_{x_1}w(\vv{x},t)$ on the right segment. Using the following result 
from calculus of fractional operators
\begin{equation}\label{eq:ISP-Op3}
\partial_t\left[e^{-i\zeta_2^2t}\partial_t^{-1/2}e^{i\zeta_2^2t}
\wtilde{w}(x_1,\zeta_2,t)\right]
=e^{-i\zeta_2^2t}\partial_t^{-1/2}e^{i\zeta_2^2t}\partial_t
\wtilde{w}(x_1,\zeta_2,t),
\end{equation}
we can recast~\eqref{eq:IPS-TBC} by taking out the operators from square 
bracketed term as
\begin{equation}
\partial_{x_1}{w}(\vv{x},t)
=-(\partial_{t}-i\partial^2_{x_2})\frac{e^{-i\pi/4}}{\sqrt{\pi}}
\int_0^t\int_{\field{R}}w(x_1,x'_2,\tau)
\frac{\mathcal{G}(x_2-x'_2,t-\tau)}{\sqrt{t-\tau}}dx_2'd\tau.
\end{equation}
In order to express the DtN map compactly, we introduce the notation 
\begin{equation}\label{eq:Op-ISP-TBC}
(\partial_t-i\partial_{x_2}^2)^{-1/2}f({x}_2,t)\\
=\frac{1}{\sqrt{\pi}}\int_0^t\int_{\field{R}}f(x'_2,\tau)
\frac{\mathcal{G}(x_2-x'_2,t-\tau)}{\sqrt{t-\tau}}dx_2'd\tau.
\end{equation}
so that
\begin{equation*}
(\partial_t-i\partial_{x_2}^2)[(\partial_t-i\partial_{x_2}^2)^{-1/2}f]
=(\partial_t-i\partial_{x_2}^2)^{1/2}f.
\end{equation*}
Next, invoking the continuity relations defined in~\eqref{eq:split-IVP}, the DtN 
map on the segment $\Gamma_r$ can now be expressed in the compact form as
\begin{equation}
\partial_{x_1}{v}(\vv{x},t)+e^{-i\pi/4}(\partial_{t}
-i\partial^2_{x_2})^{1/2}v(\vv{x},t)=0.
\end{equation}
One can obtain the transparent boundary conditions for the rest of the
boundary segments of the domain $\Omega_i$ in a similar fashion which yields 
the following equivalent formulation of the IVP corresponding to~\eqref{eq:2D-SE} 
on the computational domain $\Omega_i$:
\begin{equation}\label{eq:2D-SE-CT}
\left\{\begin{aligned}
&i\partial_tu+\triangle u=0,\quad (\vv{x},t)\in\Omega_i\times\field{R}_+,\\
&u(\vv{x},0)=u_0(\vv{x})\in \fs{L}^2(\Omega_i),\quad\supp\,u_0\subset\Omega_i,\\
&\partial_{n}{u}+
e^{-i\pi/4}(\partial_{t}-i\partial^2_{x_2})^{1/2}u=0,
\quad\vv{x}\in\Gamma_l\cup\Gamma_r,\,t>0,\\
&\partial_{n}{u}+
e^{-i\pi/4}(\partial_{t}-i\partial^2_{x_1})^{1/2}u=0,
\quad\vv{x}\in\Gamma_b\cup\Gamma_t,\,t>0.
\end{aligned}\right.
\end{equation}
Note that the operators present in the DtN maps are non-local both in space and 
time. For the numerical solution of this IBVP, a suitable numerical
implementation needs to be developed. Before we address the numerical aspects,
we develop a suitable representation for these operators at the continuous level
that is compatible with the rectangular domain. The earliest approach was due 
to Menza~\cite{Menza1997} 
which is discussed in Sec.~\ref{sec:CT-CP}. This approach is based on 
Pad\'e approximant for the square root function, however, the treatment at 
the corners is inexact. A more recent approach was presented by Feshchenko and 
Popov~\cite{FP2011} who developed a representation in terms of fractional operators. 
We discuss this approach in Sec.~\ref{sec:CT-CQ} where we adopt a more direct
way to arrive at the results as presented by the second author~\cite{V2019}.
Finally, we present one of the main contributions of this paper in
Sec.~\ref{sec:CT-NP} where a novel Pad\'{e} based approach is described that is
capable of handling the corners of the domain in a natural manner.
\subsection{Fractional operator based approach}\label{sec:CT-CQ}
\begin{figure}[!hbt]
\begin{center}
\def\myscale{0.75}
\includegraphics[scale=\myscale]{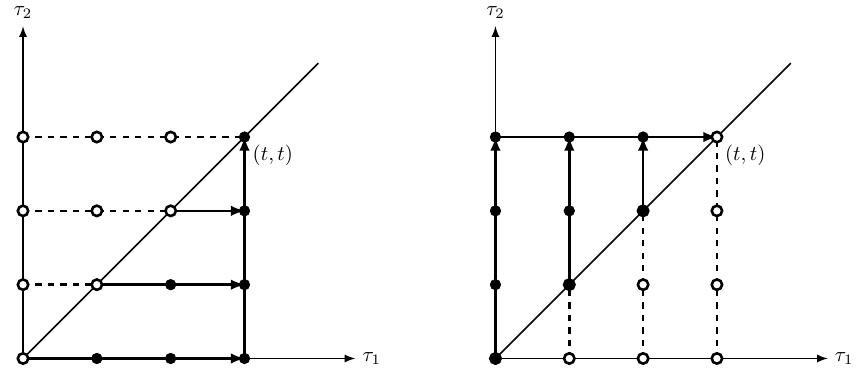}
\end{center}
\caption{\label{fig:IVP-auxi}A schematic depiction of the evolution of the 
auxiliary field $\varphi(x_1,x_2,\tau_1,\tau_2)$ in the 
$(\tau_1,\tau_2)$-plane is provided in this figure where the plot on the right 
corresponds $\vv{x}\in\Gamma_r\cup\Gamma_l$ and the plot on the left corresponds 
$\vv{x}\in\Gamma_t\cup\Gamma_b$. The filled circles depict the evolution of the 
auxiliary field $\varphi(x_1,x_2,\tau_1,\tau_2)$ either above or below the diagonal
in the $(\tau_1,\tau_2)$-plane starting from the diagonal which also serves as 
initial conditions for solving IVPs corresponding to the auxiliary function. 
The TBCs for the auxiliary field require the history of the auxiliary field 
at the corner points which makes empty circles relevant. Note that these
values at the corners can be taken from the adjacent segment of the boundary 
where it is already being computed and this is depicted by broken lines. 
Note that the vertical/horizontal lines where the arrows end corresponds to the
history of the auxiliary field needed for the TBCs on $\partial\Omega_i$ 
in the current time ($t$).}
\end{figure}

In this section, we focus on expressing the DtN maps in terms of (time-)fractional
operators. To this end, let us recall the form of the DtN map given 
by~\eqref{eq:dtn-map} in terms of the Fourier variables so that
\begin{equation}
\begin{split}
\partial_{x_1}\wtilde{w}(x_1,\zeta_2,t)
&=e^{i\pi/4}e^{-i \zeta_2^2t}\partial_t^{-1/2}e^{i\zeta_2^2t}
\left[i\partial_t\wtilde{w}(x_1,\zeta_2,t)-\zeta_2^2
\wtilde{w}(x_1,\zeta_2,t)\right]
=ie^{i\pi/4}e^{-i\zeta_2^2t}\partial_t^{1/2}e^{i\zeta_2^2t}
\wtilde{w}(x_1,\zeta_2,t)\\
&=-e^{-i\pi/4}\partial_{t'}^{1/2}
e^{-i\zeta_2^2(t-t')}\left.\wtilde{w}(x_1,\zeta_2,t')\right|_{t'=t}.
\end{split}
\end{equation}
where the last step uses the composition rule for fractional operators. Introducing 
the auxiliary function $\varphi(x_1,x_2,\tau_1,\tau_2)$ such that
\begin{equation}\label{eq:aux-fn}
\OP{F}_{x_2}[\varphi](x_1,x_2,\tau_1,\tau_2)
=e^{-i\zeta_2^2(\tau_2-\tau_1)}\wtilde{w}(x_1,\zeta_2,\tau_1),
\end{equation}
we can express the DtN map on the segment $\Gamma_r$ as: 
\begin{equation}
\partial_{x_1}w(x_1,x_2,t)=-\left.e^{-i\pi/4}\partial_{\tau_1}^{1/2}
\varphi(x_1,x_2,\tau_1,\tau_2)\right|_{\tau_1=\tau_2=t}.
\end{equation}
Numerical implementation of this fractional derivative requires the history of the
auxiliary function from the start of the computations. {\color{myrem}In order to compute 
the history, we introduce~\footnote{The idea of 
introducing an ODE with a solution satisfysing the condition instead of writing 
the condition itself first appeared in the work of Feschenko and 
Popov~\cite{FP2011} and later in~\cite{M2013}.} an IVP for the auxiliary 
function $\varphi(x_1,x_2,\tau_1,\tau_2)$ by observing the condition in~\eqref{eq:aux-fn}:}
\begin{equation}\label{eq:ivp1-aux-fn}
[i\partial_{\tau_2}+\partial^2_{x_2}]\varphi(x_1,x_2,\tau_1,\tau_2)=0,
\quad\tau_2\in(\tau_1,t],\quad x_1 = x_r,
\end{equation}
where the diagonal values $\varphi(x_1,x_2,\tau_1,\tau_1)=w(x_1,x_2,\tau_1)$ 
serve as the initial condition. 
Note that the history is needed 
for all $\tau_1\in[0,t]$ which means that we must solve the IVP for all 
$\tau_1\in[0,t]$ over $\Gamma_r$. It is interesting to note that the IVP
satisfied by the auxiliary function on the segment $\Gamma_r$ turns out to be 
a one-dimensional Schr\"{o}dinger equation. The exact boundary conditions are well-known 
for this equation when the initial data is compactly supported. 
The formulation of the TBCs for the IVP~\eqref{eq:ivp1-aux-fn} requires us to first understand the support 
(with respect to $x_2$) of the auxiliary function on the segment $\Gamma_r$ at $\tau_2=0$. From the 
knowledge of the compact support of the initial field
$u_0(\vs{x})$, we can establish the support of the auxiliary function at $\tau_2=0$. Let 
$\wtilde{u}_0(\vs{\zeta})=\OP{F}_{(x_1,x_2)}u_0(\vv{x})$, then taking the Fourier 
transform of the IVP defined by~\eqref{eq:2D-SE}, we obtain
\begin{equation}
\partial_t\wtilde{u}(\vs{\zeta})+i(\zeta_1^2+\zeta_2^2)\wtilde{u}(\vs{\zeta})=0.
\end{equation}
Taking into account the initial condition, the solution to this IVP becomes 
$\wtilde{u}(\vs{\zeta},t)= e^{-i(\zeta_1^2+\zeta_2^2)t}\wtilde{u}_0(\vs{\zeta})$
which yields
\begin{equation}\label{eq:aux-func-defn}
\begin{split}
&\OP{F}_{x_1,x_2}[\varphi](x_1,x_2,\tau_1,\tau_2)=e^{-i\zeta_2^2(\tau_2-\tau_1)}
\wtilde{u}(\zeta_1,\zeta_2,\tau_1)=e^{-i\zeta_2^2(\tau_2-\tau_1)}
e^{-i(\zeta_1^2+\zeta_2^2)\tau_1}\wtilde{u}_0(\vs{\zeta}),\\
&\varphi(x_1,x_2,\tau_1,\tau_2)=\frac{1}{(2\pi)^2}\int_{\field{R}^2}
e^{i\vs{\zeta}\cdot\vv{x}-i\zeta_1^2\tau_1-i\zeta_2^2\tau_2}
\wtilde{u}_0(\vs{\zeta})d^2\vs{\zeta}.
\end{split}
\end{equation}
From the definition of the auxiliary function written above as inverse Fourier 
transform consisting of $\wtilde{u}_0(\vs{\zeta})$, we can establish that 
the $\supp_{x_2}\varphi(x_r,x_2,\tau_1,0)\subset[x_b,x_t]$
by considering the support of the initial field at $\tau_2=0$.
It is straightforward to see that the auxiliary function also satisfies the IVP
\begin{equation}\label{eq:ivp2-aux-fn}
[i\partial_{\tau_1}+\partial^2_{x_1}]\varphi(x_1,x_2,\tau_1,\tau_2)=0,
\quad\tau_1\in(\tau_2,t],
\end{equation}
where diagonal values $\varphi(x_1,x_2,\tau_2,\tau_2)=w(x_1,x_2,\tau_2)$ serves
as the initial condition. Note that we solve the IVP for all $\tau_2\in[0,t]$.
The transparent boundary condition at the end points of the segment $\Gamma_r$ becomes:
\begin{equation}
\partial_{x_2}\varphi(x_1,x_2,\tau_1,\tau_2)
\pm e^{-i\pi/4}\partial^{1/2}_{\tau_2}\varphi(x_1,x_2,\tau_1,\tau_2)=0,
\end{equation}
where the sign is determined by $x_2\in\{x_t, x_b\}$, respectively. 
Note that the fractional derivative ($\partial_{\tau_2}^{1/2}$) requires the
knowledge of $\varphi(x_1,x_2,\tau_1,\tau_2)$ at the end points of $\Gamma_r$ 
from the start of the computation. The DtN map in terms of a time-fractional derivative
operator for the IVP~\eqref{eq:2D-SE-CT} can be expressed as
\begin{equation}\label{eq:maps-cq}
\begin{split}
&\partial_{x_1}{u}(\vv{x},t)=\left.-{e^{-i\pi/4}}\partial^{1/2}_{\tau_1}
\varphi(x_1,x_2,\tau_1,\tau_2)\right|_{\tau_1,\tau_2=t},\quad (x_1,x_2)\in\Gamma_r,\\
&\partial_{x_1}{u}(\vv{x},t)=\left.+{e^{-i\pi/4}}\partial^{1/2}_{\tau_1}
\varphi(x_1,x_2,\tau_1,\tau_2)\right|_{\tau_1,\tau_2=t},\quad (x_1,x_2)\in\Gamma_l,\\
&\partial_{x_2}{u}(\vv{x},t)=\left.-{e^{-i\pi/4}}\partial^{1/2}_{\tau_2}
\varphi(x_1,x_2,\tau_1,\tau_2)\right|_{\tau_1,\tau_2=t},\quad (x_1,x_2)\in\Gamma_t,\\
&\partial_{x_2}{u}(\vv{x},t)=\left.+{e^{-i\pi/4}}\partial^{1/2}_{\tau_2}
\varphi(x_1,x_2,\tau_1,\tau_2)\right|_{\tau_1,\tau_2=t},\quad (x_1,x_2)\in\Gamma_b.
\end{split}
\end{equation}
The IVPs for the auxiliary function are summarised below
\begin{equation}\label{eq:ivps-cq}
\begin{split}
&[i\partial_{\tau_1}+\partial_{x_1}^2]\varphi(x_1,x_2,\tau_1,\tau_2)=0,
\quad (x_1,x_2)\in\Gamma_t\cup\Gamma_b,\;\tau_1\in(\tau_2,t],\\
&[i\partial_{\tau_2}+\partial_{x_2}^2]\varphi(x_1,x_2,\tau_1,\tau_2)=0,
\quad (x_1,x_2)\in\Gamma_l\cup\Gamma_r,\;\tau_2\in(\tau_1,t].
\end{split}
\end{equation}
Let us introduce the set of corner points of the computational domain
as $\{\Gamma_{rt},\Gamma_{rb},\Gamma_{lt},\Gamma_{lb}\}$ where 
$\Gamma_{ij}=\Gamma_i\cap\Gamma_j$ with $i,j\in\{r,t,l,b\}$. The transparent boundary conditions 
at corner points for the IVPs listed in \eqref{eq:ivps-cq} are as follows
\begin{equation}\label{eq:maps-cq-auxi}
\begin{split}
&\partial_{x_1}\varphi(x_1,x_2,\tau_1,\tau_2)
+e^{-i\pi/4}\partial^{1/2}_{\tau_1}\varphi(x_1,x_2,\tau_1,\tau_2)=0,
\quad (x_1,x_2)\in\{\Gamma_{rt},\Gamma_{rb}\},\\
&\partial_{x_1}\varphi(x_1,x_2,\tau_1,\tau_2)
-e^{-i\pi/4}\partial^{1/2}_{\tau_1}\varphi(x_1,x_2,\tau_1,\tau_2)=0,
\quad (x_1,x_2)\in\{\Gamma_{lt},\Gamma_{lb}\},\\
&\partial_{x_2}\varphi(x_1,x_2,\tau_1,\tau_2)
+e^{-i\pi/4}\partial^{1/2}_{\tau_2}\varphi(x_1,x_2,\tau_1,\tau_2)=0,
\quad (x_1,x_2)\in\{\Gamma_{rt},\Gamma_{lt}\},\\
&\partial_{x_2}\varphi(x_1,x_2,\tau_1,\tau_2)
-e^{-i\pi/4}\partial^{1/2}_{\tau_2}\varphi(x_1,x_2,\tau_1,\tau_2)=0.
\quad (x_1,x_2)\in\{\Gamma_{rb},\Gamma_{lb}\},\\
\end{split}
\end{equation}
The history required for the fractional derivatives present in the DtN maps
can be understood  with the help of the schematic shown in Fig.~\ref{fig:IVP-auxi}.
The two IVPs for the auxiliary field $\varphi(x_1,x_2,\tau_1,\tau_2)$ advances the 
field either above or below the diagonal in the $(\tau_1,\tau_2)$-plane starting 
from the diagonal which also serves as initial conditions for solving IVPs. This 
is depicted by filled circles where arrows denote the direction of evolution. 
The TBCs present in~\eqref{eq:maps-cq-auxi} to solve the IVPs for the auxiliary 
functions on the boundary segments require the history of the auxiliary function 
at corners from the start of the computations which makes the empty circles 
relevant. Note that these values at the corners can be taken from the adjacent 
segment of the boundary where it is already being computed. This is depicted by 
broken lines in Fig.~\ref{fig:IVP-auxi}. 
\subsubsection{Corner conditions}\label{sec:CT-CQ-CC}
Having discussed the boundary conditions satisfied by the field on the four 
segments of the computational domain ($\Omega_i$) as stated
in~\eqref{eq:maps-cq}, we turn to the constraints on the field at the corners 
($\{\Gamma_{rt},\Gamma_{rb},\Gamma_{lt},\Gamma_{lb}\}$) which will be used later
for solving the IBVP in~\eqref{eq:2D-SE-CT} with a Galerkin method. Let us understand this by
considering the corner $\Gamma_{rt}$ and
recall the DtN map  over segment $\Gamma_r$ from~\eqref{eq:maps-cq}
\begin{equation*}
\partial_{x_1}{u}(x_1,x_2,t)=\left.-{e^{-i\pi/4}}\partial^{1/2}_{\tau_1}
\varphi(x_1,x_2,\tau_1,\tau_2)\right|_{\tau_1=\tau_2=t}.
\end{equation*}
Applying $\partial_{x_2}$ to obtain the additional constraints on the field at 
corners, we have
\begin{equation*}
\partial_{x_2}\partial_{x_1}{u}(x_1,x_2,t)=\left.-{e^{-i\pi/4}}\partial^{1/2}_{\tau_1}
\partial_{x_2}\varphi(x_1,x_2,\tau_1,\tau_2)\right|_{\tau_1=\tau_2=t}.
\end{equation*}
Now, plugging-in the value of $\partial_{x_2}\varphi(x_1,x_2,\tau_1,\tau_2)$ 
from~\eqref{eq:maps-cq-auxi} at the corner point $\Gamma_{rt}$, we obtain the following
corner condition
\begin{equation}\label{eq:cc-cq}
\partial_{x_1}\partial_{x_2}{u}(x_1,x_2,t)=\left.e^{-i\pi/4}{e^{-i\pi/4}}\partial^{1/2}_{\tau_1}
\partial^{1/2}_{\tau_2}\varphi(x_1,x_2,\tau_1,\tau_2)\right|_{\tau_1=\tau_2=t}.
\end{equation}
\subsection{Novel-Pad\'e based approach}\label{sec:CT-NP}
In this section, we would like to explore the Pad\'e approximant based representation 
for the $1/2$-order temporal derivative in~\eqref{eq:maps-cq} and its
other ramifications. It is well-known that a rational approximation to the
square-root function would allow one to obtain an effectively local approximation 
for the fractional derivative operator. We contrast this with the situation where 
the history grows with each time-step making it expensive computationally as well 
as from a storage point of view. 

The diagonal Pad\'{e} approximants for $\sqrt{1-z}$, with negative real axis as
the branch-cut, can be used to obtain rational approximation $R_M(z)$ for the 
function $\sqrt{z}$ as
\begin{equation}\label{eq:sqrt-pade}
R_M(z)
=b_0-\sum_{k=1}^M\frac{b_k}{z+\eta^2_k},\quad\text{where}\quad 
\left\{\begin{aligned}
& b_0=2M+1,\;b_k = \frac{2\eta^2_k(1+\eta^2_k)}{2M+1},\\
&\eta_k = \tan\theta_k,\;\theta_k =\frac{k\pi}{2M+1},\;k=1,2,\ldots,M.
\end{aligned}\right.
\end{equation}
Here, the expression for the diagonal Pad\'e approximant is taken from~\cite{V2011}. 
Let us start by introducing some 
compact notations for the auxiliary functions as:
\begin{equation}
\begin{split}
&\varphi_{a_1}(x_2,\tau_1,\tau_2)=\varphi(x_{a_1},x_2,\tau_1,\tau_2),\quad a_1\in\{l,r\},\\
&\varphi_{a_2}(x_1,\tau_1,\tau_2)=\varphi(x_1,x_{a_2},\tau_1,\tau_2),\quad a_2\in\{b,t\}.
\end{split}
\end{equation}
Let $z\in\field{C}$, then the Laplace transform of the auxiliary function can be denoted as  
\begin{equation}
\OP{L}_{\tau_1}[\varphi_{a_1}(x_2,\tau_1,\tau_2)]=\varPhi_{a_1}(x_2,z,\tau_2),\quad 
\OP{L}_{\tau_2}[\varphi_{a_2}(x_1,\tau_1,\tau_2)]=\varPhi_{a_2}(x_1,\tau_1,z). 
\end{equation}
The DtN maps described in~\eqref{eq:maps-cq} on the segments $\Gamma_{a_1}$ and $\Gamma_{a_2}$ 
can be restated as 
\begin{equation}
\begin{split}
&\partial_{x_1}{u}(x_1,x_2,t)\pm e^{-i\pi/4}\left.\OP{L}^{-1}_{z}\left[
\sqrt{z}\varPhi_{a_1}(x_2,z,\tau_2)\right]\right|_{\tau_1=\tau_2=t}=0,\\
&\partial_{x_2}{u}(x_1,x_2,t)\pm e^{-i\pi/4}\left.\OP{L}^{-1}_{z}\left[
\sqrt{z}\varPhi_{a_2}(x_2,\tau_1,z)\right]\right|_{\tau_1=\tau_2=t}=0.
\end{split}
\end{equation}
Inserting the rational approximation of order $M$ for the function $\sqrt{z}$ as
\begin{equation}\label{eq:maps-npade1}
\begin{split}
&\partial_{x_1}{u}(x_1,x_2,t)\pm e^{-i\pi/4}\left.\OP{L}^{-1}_{z}\left[
R_M(z)\varPhi_{a_1}(x_2,z,\tau_2)\right]\right|_{\tau_1=\tau_2=t}\approx 0,\\
&\partial_{x_2}{u}(x_1,x_2,t)\pm e^{-i\pi/4}\left.\OP{L}^{-1}_{z}\left[
R_M(z)\varPhi_{a_2}(x_2,\tau_1,z)\right]\right|_{\tau_1=\tau_2=t}\approx 0.
\end{split}
\end{equation}
The map expressed above is still nonlocal, however, the rational nature of 
the expression allows us to obtain an effectively local approximation for the 
DtN operation. {\color{myrem}Following Lindmann's approach~\cite{L1985}, we 
introduce the auxiliary fields $\varphi_{k,a_1}(x_2,\tau_1,\tau_2)$ and
$\varphi_{k,a_2}(x_1,\tau_1,\tau_2)$ at the boundary segments $\Gamma_{a_1}$
and $\Gamma_{a_2}$, respectively, such that}
\begin{equation}
\begin{split}
&(z+\eta^2_k)^{-1}\varPhi_{a_1}(x_2,z,\tau_2)=\varPhi_{k,a_1}(x_2,z,\tau_2),\\
&(z+\eta^2_k)^{-1}\varPhi_{a_2}(x_1,\tau_1,z)=\varPhi_{k,a_2}(x_1,\tau_1,z),
\end{split}
\end{equation}
where $k=1,2,\ldots,M$. In the physical space, every $\varphi_{k,a_1}$ and 
$\varphi_{k,a_2}$ satisfy the following ODEs: 
\begin{equation}\label{eq:ode-auxi-npade1}
\begin{split}
&(\partial_{\tau_1}+\eta^2_k)\varphi_{k,a_1}(x_2,\tau_1,\tau_2)=
\varphi_{a_1}(x_2,\tau_1,\tau_2),\\
&(\partial_{\tau_2}+\eta^2_k)\varphi_{k,a_2}(x_1,\tau_1,\tau_2)=
\varphi_{a_2}(x_1,\tau_1,\tau_2),
\end{split}
\end{equation}
with the initial conditions assumed to be $\varphi_{k,a_1}(x_2,0,\tau_2)=0$ and
$\varphi_{k,a_2}(x_1,\tau_1,0)=0$, respectively. 
The DtN maps for the interior problem present in~\eqref{eq:maps-npade1} now reads as 
\begin{equation}\label{eq:maps-pade}
\begin{split}
& \partial_{x_1}{u}(\vv{x},t)\pm e^{-i\pi/4}\left[ b_0 u(\vv{x},t) -
  \sum_{k=1}^M b_k\varphi_{k,a_1}(x_2,t,t)\right]\approx 0,\\
& \partial_{x_2}{u}(\vv{x},t)\pm e^{-i\pi/4}\left[ b_0 u(\vv{x},t) -
  \sum_{k=1}^M b_k\varphi_{k,a_2}(x_1,t,t)\right]\approx 0.
\end{split}
\end{equation}
The solution to the ODEs~\eqref{eq:ode-auxi-npade1} reads as
\begin{equation}\label{eq:phi_k}
\begin{split}
& \varphi_{k,a_1}(x_2,\tau_1,\tau_2)
  =\int_0^{\tau_1}e^{-\eta^2_k(\tau_1-s_1)}\varphi_{a_1}(x_2,s_1,\tau_2)ds_1,
   \quad a_1\in\{l,r\},\\
& \varphi_{k,a_2}(x_1,\tau_1,\tau_2)
  =\int_0^{\tau_2}e^{-\eta^2_k(\tau_2-s_2)}\varphi_{a_2}(x_1,\tau_1,s_2)ds_2,
   \quad a_2\in\{t,b\}.
\end{split}
\end{equation}
It is straightforward to conclude that auxiliary fields satisfy the same IVPs which are 
satisfied by auxiliary functions
\begin{equation}\label{eq:ivp-pade-auxi}
\begin{split}
&[i\partial_{\tau_2}+\partial_{x_2}^2]\varphi_{k,a_1}(x_2,\tau_1,\tau_2)=0,
\quad (x_1,x_2)\in\Gamma_l\cup\Gamma_r,\;\tau_2\in(\tau_1,t],\\
&[i\partial_{\tau_1}+\partial_{x_1}^2]\varphi_{k,a_2}(x_1,\tau_1,\tau_2)=0,
\quad (x_1,x_2)\in\Gamma_t\cup\Gamma_b,\;\tau_1\in(\tau_2,t].
\end{split}
\end{equation}
The transparent boundary conditions for these IVPs can be written as:
\begin{equation}\label{eq:maps-pade-auxi}
\begin{split}
&\partial_{x_2}\varphi_{k,a_1}(x_b,\tau_1,\tau_2)
-e^{-i\pi/4}\partial^{1/2}_{\tau_2}\varphi_{k,a_1}(x_b,\tau_1,\tau_2)=0,\\
&\partial_{x_2}\varphi_{k,a_1}(x_t,\tau_1,\tau_2)
+e^{-i\pi/4}\partial^{1/2}_{\tau_2}\varphi_{k,a_1}(x_t,\tau_1,\tau_2)=0,\\
&\partial_{x_1}\varphi_{k,a_2}(x_l,\tau_1,\tau_2)
-e^{-i\pi/4}\partial^{1/2}_{\tau_1}\varphi_{k,a_2}(x_l,\tau_1,\tau_2)=0,\\
&\partial_{x_1}\varphi_{k,a_2}(x_r,\tau_1,\tau_2)
+e^{-i\pi/4}\partial^{1/2}_{\tau_1}\varphi_{k,a_2}(x_r,\tau_1,\tau_2)=0.
\end{split}
\end{equation}
Once again, we use the Pad\'e approximants based representation for the 
$1/2$-order temporal derivative operator present in the TBCs for 
the auxiliary fields described in~\eqref{eq:maps-pade-auxi}.
Introducing the auxiliary fields $\psi_{k,k',a_1,a_2}(\tau_1,\tau_2)$ and 
$\psi_{k,k',a_2,a_1}(\tau_1,\tau_2)$ at the end points of the boundary segments 
$\Gamma_{a_1}$ and $\Gamma_{a_2}$, respectively, such that 
\begin{equation}
\begin{split}
& (z+\eta^2_{k'})^{-1}\varPhi_{k,a_1}(x_{a_2},\tau_1,z)=\Psi_{k,k',a_1,a_2}(\tau_1,z),
  \quad \;a_1\in\{l,r\},\\
& (z+\eta^2_{k'})^{-1}{\varPhi}_{k,a_2}(x_{a_1},\tau_2,z)={\Psi}_{k,k', a_2,a_1}(z,\tau_2),
  \quad \;a_2\in\{b,t\}.
\end{split}
\end{equation}
In the physical space, every $\psi_{k,k',a_1,a_2}$ and $\psi_{k,k',a_2,a_1}$ satisfy
the following ODEs:
\begin{equation}\label{eq:ode-auxi-npade2}
\begin{split}
& (\partial_{\tau_2}+\eta^2_{k'})\psi_{k,k',a_1,a_2}(\tau_1,\tau_2)=
    \varphi_{k,a_1}(x_{a_2},\tau_1,\tau_2),\\
& (\partial_{\tau_1}+\eta^2_{k'})\psi_{k,k',a_2,a_1}(\tau_1,\tau_2)= 
    \varphi_{k,a_2}(x_{a_1},\tau_1,\tau_2),
\end{split}
\end{equation}
with the initial conditions assumed to be $\psi_{k,k',a_1,a_2}(\tau_1,0)=0$ and
$\psi_{k,k',a_2,a_1}(0,\tau_2)=0$, respectively. 
The DtN maps for the auxiliary fields~\eqref{eq:maps-pade-auxi} now reads as 
\begin{equation}\label{eq:dtn-using-psi}
\begin{split}
& \partial_{x_2}\varphi_{k,a_1}(x_{a_2},\tau_1,\tau_2)\pm e^{-i\pi/4}\left[
  b_0\varphi_{k,a_1}(x_{a_2},\tau_1,\tau_2) -
  \sum_{k'=1}^M b_{k'}\psi_{k,k',a_1,a_2}(\tau_1,\tau_2)\right]=0,\\
& \partial_{x_1}\varphi_{k,a_2}(x_{a_1},\tau_1,\tau_2)\pm e^{-i\pi/4}\left[ 
  b_0\varphi_{k,a_2}(x_{a_1},\tau_1,\tau_2) -
  \sum_{k'=1}^M b_{k'}\psi_{k,k',a_2,a_1}(\tau_1,\tau_2)\right]=0.
\end{split}
\end{equation}
The solution to the ODEs~\eqref{eq:ode-auxi-npade2} reads as
\begin{equation}
\begin{split}
\psi_{k,k',a_1,a_2}(\tau_1,\tau_2)
&=\int_0^{\tau_2}e^{-\eta^2_{k'}(\tau_2-s_2)}\varphi_{k,a_1}(x_{a_2},\tau_1,s_2)ds_2,\\
\psi_{k,k',a_2,a_1}(\tau_1,\tau_2)
&=\int_0^{\tau_1}e^{-\eta^2_{k'}(\tau_1-s_1)}\varphi_{k,a_2}(x_{a_1},s_1,\tau_2)ds_1.\\
\end{split}
\end{equation}
\begin{rem}\label{rem:aux-func-psi}
It is interesting to note that fields $\psi_{k,k',a_1,a_2}(\tau_1,\tau_2)$ and 
$\psi_{k,k',a_2,a_1}(\tau_1,\tau_2)$ are transpose of each other. We can verify it by 
inserting the value of auxiliary fields $\varphi_{k,a_1}$ and $\varphi_{k,a_2}$ 
from~\eqref{eq:phi_k} as:
\begin{equation}
\begin{split}
\psi_{k,k',a_1,a_2}(\tau_1,\tau_2)
&=\int_0^{\tau_2}e^{-\eta^2_{k'}(\tau_2-s_2)}\varphi_{k,a_1}(x_{a_2},\tau_1,s_2)ds_2\\
&=\int_0^{\tau_2}e^{-\eta^2_{k'}(\tau_2-s_2)}\int_0^{\tau_1}e^{-\eta^2_{k}(\tau_1-s_1)}
\varphi(x_{a_1},x_{a_2},s_1,s_2)ds_1 ds_2\\
&=\int_0^{\tau_1}e^{-\eta^2_{k}(\tau_1-s_1)}\int_0^{\tau_2}e^{-\eta^2_{k'}(\tau_2-s_2)}
\varphi(x_{a_1},x_{a_2},s_1,s_2)ds_2 ds_1\\
&=\int_0^{\tau_1}e^{-\eta^2_{k}(\tau_1-s_1)}
\varphi_{k',a_2}(x_{a_1},s_1,\tau_2)ds_1
=\psi_{k',k,a_2,a_1}(\tau_1,\tau_2).
\end{split}
\end{equation}
This observation allows us to develop a rather efficient numerical scheme from a 
storage point of view.
\end{rem}

\begin{figure}[!htbp]
\begin{center}
\includegraphics[scale=0.75]{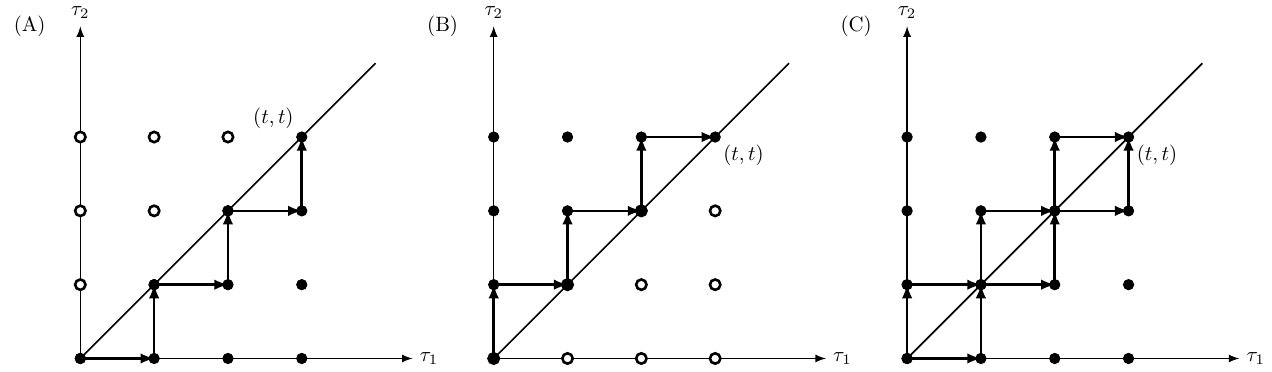}
\end{center}
\caption{\label{fig:IVP-auxi-pade}A schematic depiction of the evolution of the 
auxiliary fields $\varphi_{k,a_1}(x_2,\tau_1,\tau_2),\;
\varphi_{k,a_2}(x_1,\tau_1,\tau_2)$ and $\psi_{k,k',a_1,a_2}(\tau_1,\tau_2)$ in
the $(\tau_1,\tau_2)$-plane is provided in this figure. The plots (A)
and (B) depict the evolution of the fields $\varphi_{k,a_2}(x_1,\tau_1,\tau_2)$ 
and $\varphi_{k,a_1}(x_2,\tau_1,\tau_2)$ on the boundary segments $\Gamma_{a_2}$
and $\Gamma_{a_1}$, respectively. The plot (C) depicts the evolution of the field 
$\psi_{k,k',a_1,a_2}(\tau_1,\tau_2)$ which can be achieved by moving either below 
or above the diagonal.} 
\end{figure}

In order to advance the auxiliary fields ($\psi$) from diagonal to diagonal point in 
$(\tau_1,\tau_2)$-plane, we need to establish the ODEs to advance the auxiliary fields 
in the directions orthogonal to those mentioned in~\eqref{eq:ode-auxi-npade2}. One can
achieve that by observing
\begin{equation}
\begin{split}
[\partial_{\tau_1}+\eta^2_{k}]\psi_{k,k',a_1,a_2}(\tau_1,\tau_2)
&=\int_0^{\tau_2}e^{-\eta^2_{k'}(\tau_2-s_2)}[\partial_{\tau_1}+\eta^2_{k}]
\varphi_{k,a_1}(x_{a_2},\tau_1,s_2)ds_2\\
&=\int_0^{\tau_2}e^{-\eta^2_{k'}(\tau_2-s_2)}\varphi(x_{a_1},x_{a_2},\tau_1,s_2)ds_2
=\varphi_{k',a_2}(x_{a_1},\tau_1,\tau_2),\\
[\partial_{\tau_2}+\eta^2_{k}]\psi_{k,k',a_2,a_1}(\tau_1,\tau_2)
&=\int_0^{\tau_1}e^{-\eta^2_{k'}(\tau_1-s_1)}[\partial_{\tau_2}+\eta^2_{k}]
\varphi_{k,a_2}(x_{a_1},s_1,\tau_2)ds_1\\
&=\int_0^{\tau_1}e^{-\eta^2_{k'}(\tau_1-s_1)}\varphi(x_{a_1},x_{a_2},s_1,\tau_2)ds_1
=\varphi_{k',a_1}(x_{a_1},\tau_1,\tau_2).
\end{split}
\end{equation}
These ODEs for the auxiliary fields can now be summarized as
\begin{equation}\label{eq:ode-auxi-npade3}
\begin{split}
&(\partial_{\tau_1}+\eta^2_{k})\psi_{k,k',a_1,a_2}(\tau_1,\tau_2)
=\varphi_{k',a_2}(x_{a_1},\tau_1,\tau_2),\\
&(\partial_{\tau_2}+\eta^2_{k})\psi_{k,k',a_2,a_1}(\tau_1,\tau_2)
=\varphi_{k',a_1}(x_{a_2},\tau_1,\tau_2).
\end{split}
\end{equation}
The evolution of the auxiliary fields $\varphi_{k,a_1}(x_2,\tau_1,\tau_2),\;
\varphi_{k,a_2}(x_1,\tau_1,\tau_2)$ and $\psi_{k,k',a_1,a_2}(\tau_1,\tau_2)$
in the $(\tau_1,\tau_2)$-plane can be understood with the help of the schematic 
shown in Fig.~\ref{fig:IVP-auxi-pade}. In this schematic, Fig.~\ref{fig:IVP-auxi-pade}(A)
and Fig.~\ref{fig:IVP-auxi-pade}(B) depict the evolution of the fields 
$\varphi_{k,a_2}(x_1,\tau_1,\tau_2)$ and $\varphi_{k,a_1}(x_2,\tau_1,\tau_2)$ 
on the boundary segments $\Gamma_{a_2}$
and $\Gamma_{a_1}$, respectively. The diagonal to diagonal computation 
of the fields $\varphi_{k,a_1}(x_2,\tau_1,\tau_2)$ and
$\varphi_{k,a_2}(x_1,\tau_1,\tau_2)$ consists of first advancing the fields using the IVPs
established in~\eqref{eq:ivp-pade-auxi} and then using the ODEs in~\eqref{eq:ode-auxi-npade1} 
for the second movement. Our novel Pad\'e approach makes the numerical scheme efficient 
from a storage point of view which is obvious from the schematic presented. Similarly,
Fig.~\ref{fig:IVP-auxi-pade}(C) depicts the evolution of the field 
$\psi_{k,k',a_1,a_2}(\tau_1,\tau_2)$ which can be achieved by moving either below or 
above the diagonal. 
\subsubsection{Corner conditions}\label{sec:CT-NP-CC}
With in the novel Pad\'e approach, the constraints satisfied by the field at the 
corner can be obtained following a similar approach as
that of the Sec.~\ref{sec:CT-CQ-CC}. Consider the corner $\Gamma_{rt}$ and
recall the DtN map  over segment $\Gamma_r$ from~\eqref{eq:maps-pade}
\begin{equation*}
\partial_{x_1}{u}(\vv{x},t)+ e^{-i\pi/4}\left[ 
b_0 u(\vv{x},t) - \sum_{k=1}^M b_k\varphi_{k,a_1}(x_2,t,t)\right]=0.
\end{equation*}
Applying $\partial_{x_2}$ to obtain the additional constraints on the field at 
corners, we have
\begin{equation*}
\partial_{x_2}\partial_{x_1}u(\vv{x},t)
+e^{-i\pi/4}\left[ b_0\partial_{x_2}u(\vv{x},t)-
\sum_{k=1}^M b_k\partial_{x_2}\varphi_{k,a_1}(x_2,t,t)\right]=0.
\end{equation*}
Plugging-in the value of $\partial_{x_2}\varphi_{k,a_1}(x_2,t,t)$ 
from~\eqref{eq:maps-pade-auxi} at corner point $\Gamma_{rt}$ we get
\begin{equation*}
\partial_{x_2}\partial_{x_1} u(\vv{x},t)
+e^{-i\pi/4}b_0\partial_{x_2}u(\vv{x},t)+b_ke^{-i\pi/2}\sum_{k=1}^M 
\left[b_0\varphi_{k,a_1}(x_2,t,t) 
-\sum_{k'=1}^M b_{k'}\psi_{k,k',a_1,a_2}(t,t)\right]=0.
\end{equation*}
The corner condition at $\Gamma_{rt}$ becomes:
\begin{equation}\label{eq:cc-npade}
\partial_{x_2}\partial_{x_1}u(\vv{x},t)
+e^{-i\pi/4}b_0\left[\partial_{x_2}u(\vv{x},t)
+\partial_{x_1}u(\vv{x},t)+e^{-i\pi/4}b_0u(\vv{x},t)\right]\\
-e^{-i\pi/2}\sum_{k=1}^M\sum_{k'=1}^M b_kb_{k'}\psi_{k,k',a_1,a_2}(t,t)=0.
\end{equation}
\subsection{Difficulty with conventional Pad\'{e} based approach}\label{sec:CT-CP}
Following Menza's approach~\cite{Menza1997}, a Pad\'e approximant based 
representation for the operators 
$(\partial_{t}-i\partial^2_{x_j})^{1/2},\,j=1,2$ can be developed by considering 
the rational function in~\eqref{eq:sqrt-pade}. Using the Fourier representation
in space and Laplace in time, we can write the DtN map on the segments
$\Gamma_{a_1}$ as
\begin{equation}
\partial_{x_1}u(x_1,x_2,t)\pm e^{-i\pi/4}\OP{L}_z^{-1}\OP{F}^{-1}_{\zeta_2}
\left[R_M(z+i\zeta_2^2)\wtilde{U}(x_1,\zeta_2,z)\right]=0,\quad x_1\in\{x_l,x_r\}.
\end{equation}
Introducing the auxiliary fields 
$\{\varphi_{k,a_1}(x_2,t)
\leftrightarrow\wtilde{\varPhi}_{k,a_1}(\zeta_2,z),\,z\in\field{C},\,\zeta_2\in\field{R}\}$ 
at the boundary segments denoted by $\Gamma_{a_1}$ such that 
\begin{equation}
(z+i\zeta_2^2+\eta^2_k)^{-1}\wtilde{U}(x_1,\zeta_2,z)=\wtilde{\varPhi}_{k,a_1}(\zeta_2,z),
\quad k=1,2,\ldots,M.
\end{equation}
so that the DtN map becomes
\begin{equation}
\partial_{x_1}{u}(\vv{x},t)\pm e^{-i\pi/4}\left[ b_0w(\vv{x},t) -
\sum_{k=1}^M b_k\varphi_{k,a_1}(x_2,t)\right]=0,
\end{equation}
and every $\varphi_{k,a_1}$ satisfies the IVP given by
\begin{equation}\label{eq:auxi-cp}
(i\partial_t + \partial^2_{x_2} +i\eta^2_k)\varphi_{k,a_1}(x_2,t)= iu(x_1,x_2,t),
\quad x_1\in\{x_l,x_r\},
\end{equation}
with the initial conditions set to $\varphi_{k,a_1}(x_2,0)=0$. The equation
above is a driven one-dimensional Schr\"odinger equation which can be handled 
in the standard way if $u(x_1,x_2,t)$ does not hit the end points of
$\Gamma_{a_1}$. Under this assumption, the transparent boundary conditions at
the end points are given by
\begin{equation}\label{eq:auxi-cp-corner}
\begin{split}
&\partial_{x_2}\varphi_{k,a_1}(x_t,t)
+e^{-i\pi/4}e^{-\eta_k^2 t}\partial_t^{1/2}
\left[e^{+\eta_k^2 t}\varphi_{k,a_1}(x_t,t)\right]=0,\\
&\partial_{x_2}\varphi_{k,a_1}(x_b,t)
-e^{-i\pi/4}e^{-\eta_k^2 t}\partial_t^{1/2}
\left[e^{+\eta_k^2 t}\varphi_{k,a_1}(x_b,t)\right]=0.
\end{split}
\end{equation}
These conditions become inaccurate as soon as the support of  $u(x_1,x_2,t)$
extends to the corners. While in principle, it is possible to remedy the problem
here (irrespective of the computational cost), we do not pursue this line of 
investigation on account
of the fact that the value of the field $u(x_1,x_2,t)$ is required 
outside of $\Omega_i$ in order to write the exact nonreflecting boundary 
conditions for $\varphi_{k,a_1}$.

\section{Numerical Implementation}\label{sec:numerical-implementation}
In this section, we address the complete numerical solution 
of the initial boundary-value problem (IBVP) stated in~\eqref{eq:2D-SE-CT}.
In order to formulate the discrete linear system corresponding to the IBVP, we 
need to introduce temporal as well as spatial discretization of the problem. 
We first consider the temporal derivative which is discretized using one-step 
methods, namely, backward differentiation formula of order 1 (BDF1) and the 
trapezoidal rule (TR). Subsequently, a compatible temporal discretization 
scheme is developed for the boundary conditions where we strive to formulate 
the novel boundary maps at the time-discrete level as a Robin-type 
boundary condition for each time-step. Finally, the resulting 
semi-discrete partial differential equation is solved using a Legendre-Galerkin 
method where the boundary condition is enforced exactly by designing a 
boundary-adapted basis in terms of the Legendre polynomials followed by a 
boundary lifting process.

For the computational domain $\Omega_i$, we introduce a reference domain 
$\Omega_i^{\text{ref.}}=\field{I}\times\field{I}$ keeping in view the typical 
domain of definition of the orthogonal polynomials being used.
In order to describe the associated linear maps between the reference domain and 
the actual computational domain, we introduce the variables 
$y_1,y_2\in\Omega_i^{\text{ref.}}$ such that
\begin{equation}
\left\{\begin{aligned}
&x_1 = J_1 y_1+\bar{x}_1,\quad J_1 = \frac{1}{2}(x_r-x_l),\quad
\bar{x}_1=\frac{1}{2}(x_l+x_r), & \beta_1=J_1^{-2},\\
&x_2 = J_2 y_2+\bar{x}_2,\quad J_2 = \frac{1}{2}(x_t-x_b),\quad
\bar{x}_2=\frac{1}{2}(x_b+x_t), & \beta_2=J_2^{-2}.
\end{aligned}\right.
\end{equation}
Let $\Delta t$ denote the time-step. 
In the rest of the section, with a slight abuse of notation, we switch to the 
variables $y_1\in\field{I}$ and $y_2\in\field{I}$ for all the discrete approximations 
of the dependent variables. For instance, $u^j(\vv{y})$ is taken to approximate 
$u(\vv{x},j\Delta t)$ for $j=0,1,2,\ldots,N_t-1$. The temporal discretization of the 
interior problem using one-step methods is enumerated below:
\begin{itemize}
\item The BDF1 based discretization is given by
\begin{equation}\label{eq:2d-se-bdf1}
i\frac{u^{j+1}-u^{j}}{\Delta t}+\triangle u^{j+1}=0
\implies\left(\beta_1\partial^2_{y_1}
+\beta_2\partial^2_{y_2}\right)u^{j+1}+i\rho u^{j+1}
=i\rho u^j,\quad \rho = 1/\Delta t.
\end{equation}
\item The TR based discretization is given by
\begin{equation}\label{eq:2d-disc-tr}
i\frac{u^{j+1}-u^{j}}{\Delta t}+\triangle u^{j+1/2}=0
\implies\left(\beta_1\partial^2_{y_1}+\beta_2\partial^2_{y_2}\right)v^{j+1}+i\rho v^{j+1}
=i\rho u^j,\quad \rho = 2/\Delta t,
\end{equation}
where we have used the staggered samples of the field defined as
\begin{equation}
v^{j+1} = u^{j+1/2}=\frac{u^{j+1}+u^j}{2},\quad v^0=0.
\end{equation}
\end{itemize}
The rest of this section is organized as follows: Sec.~\ref{sec:dtn-maps}
addresses the temporal discretization of the boundary maps which is followed 
by a discussion of the resulting spatial problem in Sec.~\ref{sec:IBVP}.
\subsection{Discretizing the boundary conditions}\label{sec:dtn-maps}
As mentioned above, we seek a compatible temporal discretization of the novel 
boundary conditions discussed in earlier sections in order to solve the IBVP. To 
this end, we first discuss the convolution quadrature (CQ) based discretization 
of the $1/2$-order time-fractional derivative present in the TBCs developed in 
Sec.~\ref{sec:CT-CQ} and its other ramifications. Here, we would like to
emphasize that an unsuitable temporal discretization of the convolution operators
may destroy the overall stability of the numerical scheme~\cite{Mayfield1989}. 
In order to avoid this, we try to match the time-stepping methods in the CQ based 
discretization of the boundary conditions to that of the interior problem 
(namely, BDF1 and TR). Next, we discuss the temporal discretization 
of the boundary maps developed in Sec.~\ref{sec:CT-NP} using Pad\'{e} approximants 
based representation for the fractional operator present in the TBCs. This
method is labelled as `NP' which expands to novel Pad\'e.

In order to facilitate the temporal discretization of the TBCs, we introduce the 
equispaced samples of the auxiliary function, $\varphi^{j,k}(y_1,y_2)$, to 
approximate $\varphi(x_1,x_2,j\Delta t,k\Delta t)$ for $j,k=0,1,\ldots N_t-1$. It 
follows from the definition~\eqref{eq:aux-func-defn} that the diagonal values of 
the auxiliary function at the discrete level are determined by the interior field, i.e.,
$\varphi^{j,j}(y_1,y_2) = u^j(y_1,y_2)$.

\subsubsection{CQ--BDF1}
The discrete scheme for the boundary maps is said to be `CQ--BDF1' if the 
underlying one-step method in the CQ scheme is BDF1. Let 
$(\omega_k)_{k\in\field{N}_0}$ denote the
corresponding quadrature weights, then
\begin{equation}
\omega_j=\left[{\left(j-{3}/{2}\right)}/{j}\right]\omega_{j-1},\quad j\geq 1,
\quad\text{with}\quad\omega_0=1.
\end{equation}
Here, the quadrature weights are obtained using semi-discrete formulation as
discussed in~\cite{SV2023}. The resulting discretization scheme for 
the DtN map described in~\eqref{eq:maps-cq} in terms of $\varphi^{j,k}(y_1,y_2)$ 
stated as a Robin-type boundary condition reads as
\begin{equation}
\begin{split}
&\partial_{y_1}u^{j+1}+\sqrt{\rho/\beta_1}e^{-i\pi/4}u^{j+1}
=-\sqrt{\rho/\beta_1}e^{-i\pi/4}\OP{B}\left[\{\varphi^{k,j+1}(y_r,y_2)\}_{k=0}^j\right],\\
&\partial_{y_1}u^{j+1}-\sqrt{\rho/\beta_1}e^{-i\pi/4}u^{j+1}
=+\sqrt{\rho/\beta_1}e^{-i\pi/4}\OP{B}\left[\{\varphi^{k,j+1}(y_l,y_2)\}_{k=0}^j\right],\\
&\partial_{y_2}u^{j+1}+\sqrt{\rho/\beta_2}e^{-i\pi/4}u^{j+1}
=-\sqrt{\rho/\beta_2}e^{-i\pi/4}\OP{B}\left[\{\varphi^{j+1,k}(y_1,y_t)\}_{k=0}^j\right],\\
&\partial_{y_2}u^{j+1}-\sqrt{\rho/\beta_2}e^{-i\pi/4}u^{j+1}
=+\sqrt{\rho/\beta_2}e^{-i\pi/4}\OP{B}\left[\{\varphi^{j+1,k}(y_1,y_b)\}_{k=0}^j\right],
\end{split}
\end{equation}
where we have used the fact that $\varphi^{j+1,j+1}(y_1,y_2)=u^{j+1}(y_1,y_2)$ and 
\begin{equation}
\begin{split}
&\mathcal{B}^{j+1}_{a_1}(y_2)
=\OP{B}\left[\{\varphi^{k,j+1}(y_{a_1},y_2)\}_{k=0}^j\right]
=\sum_{k=1}^{j+1}\omega_{k}\varphi^{j+1-k,j+1}(y_{a_1},y_2),\quad
a_1\in\{r,l\},\\
&\mathcal{B}^{j+1}_{a_2}(y_1)
=\OP{B}\left[\{\varphi^{j+1,k}(y_1,y_{a_2})\}_{k=0}^j\right]
=\sum_{k=1}^{j+1}\omega_{k}\varphi^{j+1,j+1-k}(y_1,y_{a_2}),\quad
a_2\in\{b,t\}.
\end{split}
\end{equation}
We may refer to these functions as the \emph{history functions} on account of
the fact that, for the latest time-step, they can be assumed to be determined 
by the quantities that are already determined in the previous time-steps.
Set $\alpha_k=\sqrt{\rho/\beta_k}e^{-i\pi/4}$ for $k=1,2$, then the discretized
maps take the form
\begin{equation}\label{eq:maps-cq-bdf1}
\begin{split}
&\partial_{y_1}u^{j+1}+\alpha_1 u^{j+1}=-\alpha_1\mathcal{B}^{j+1}_r(y_2),\quad
 \partial_{y_2}u^{j+1}+\alpha_2 u^{j+1}=-\alpha_2\mathcal{B}^{j+1}_t(y_1),\\
&\partial_{y_1}u^{j+1}-\alpha_1 u^{j+1}=+\alpha_1\mathcal{B}^{j+1}_l(y_2),\quad
 \partial_{y_2}u^{j+1}-\alpha_2 u^{j+1}=+\alpha_2\mathcal{B}^{j+1}_b(y_1).
\end{split}
\end{equation}
To obtain the history of the field $\varphi(y_1,y_2,\tau_1,\tau_2)$ needed to
compute the history functions in~\eqref{eq:maps-cq-bdf1}, we employ a 
Legendre-Galerkin method to the solve the IVPs listed in~\eqref{eq:ivps-cq} 
by formulating the boundary conditions given in~\eqref{eq:maps-cq-auxi} as a 
Robin-type map. Let us start by introducing a compact notation obtained by 
restricting the auxiliary function in the following manner:
\begin{equation}
\varphi^{p,q}_{a_1}(y_2)=\varphi^{p,q}(y_{a_1},y_2),\quad
\varphi^{p,q}_{a_2}(y_1)=\varphi^{p,q}(y_1,y_{a_2}),\quad 
\varphi^{p,q}_{a_1a_2}=\varphi^{p,q}(y_{a_1},y_{a_2}),\quad 
a_1\in\{l,r\},\;a_2\in\{b,t\}.
\end{equation}
The BDF1-based discretization of the IVPs~\eqref{eq:ivps-cq} for the auxiliary 
functions read as follows: 
\begin{equation}
\begin{split}
&(y_1,y_2)\in{\Gamma}^{\text{ref.}}_b\cup{\Gamma}^{\text{ref.}}_t:\quad
\left\{\begin{aligned}
&i\frac{\varphi^{p+1,q}_{a_2}-\varphi^{p,q}_{a_2}}{\Delta t}
+\beta_1\partial_{y_1}^2\varphi^{p+1,q}_{a_2}=0
\implies-\alpha^{-2}_1\partial_{y_1}^2\varphi^{p+1,q}_{a_2}+\varphi^{p+1,q}_{a_2}
=\varphi^{p,q}_{a_2},\\
&p=q,q+1,\ldots,j,\quad q = 0,1,\ldots,j;
\end{aligned}\right.\\
&(y_1,y_2)\in\Gamma^{\text{ref.}}_l\cup\Gamma^{\text{ref.}}_r:\quad
\left\{\begin{aligned}
&i\frac{\varphi^{m,n+1}_{a_1}-\varphi^{m,n}_{a_1}}{\Delta t}
+\beta_2\partial_{y_2}^2\varphi^{m,n+1}_{a_1}=0\implies
-\alpha^{-2}_2\partial_{y_2}^2\varphi^{m,n+1}_{a_1}+\varphi^{m,n+1}_{a_1}
=\varphi^{m,n}_{a_1},\\
&n=m,m+1,\ldots,j,\quad m = 0,1,\ldots,j.
\end{aligned}\right.
\end{split}
\end{equation}
The BDF1-based discretization of the DtN map described in~\eqref{eq:maps-cq-auxi} stated
as Robin-type boundary condition reads as:
\begin{equation}
\begin{split}
&\left.(\partial_{y_1}+\alpha_1)\varphi^{p+1,q}_{a_2}(y_1)\right|_{y_1=y_r} =-\alpha_1\OP{B}\left[\{\varphi^{k,q}_{a_2}(y_r)\}_{k=0}^p\right],\\
&\left.(\partial_{y_1}-\alpha_1)\varphi^{p+1,q}_{a_2}(y_1)\right|_{y_1=y_l} =+\alpha_1\OP{B}\left[\{\varphi^{k,q}_{a_2}(y_l)\}_{k=0}^p\right],\\
&\left.(\partial_{y_2}+\alpha_2)\varphi^{m,n+1}_{a_1}(y_2)\right|_{y_1=y_t} =-\alpha_2\OP{B}\left[\{\varphi^{m,k}_{a_1}(y_t)\}_{k=0}^n\right],\\
&\left.(\partial_{y_2}-\alpha_2)\varphi^{m,n+1}_{a_1}(y_2)\right|_{y_1=y_b} =+\alpha_2\OP{B}\left[\{\varphi^{m,k}_{a_1}(y_b)\}_{k=0}^m\right],
\end{split}
\end{equation}
where 
\begin{equation}
\begin{split}
&\mathcal{B}^{p+1,q}_{a_2,a_1}
 =\OP{B}\left[\{\varphi^{k,q}_{a_2}(y_{a_1})\}_{k=0}^p\right]
 =\sum_{k=1}^{p+1}\omega_{k}\varphi^{p+1-k,q}(y_{a_1},y_{a_2})
 =\sum_{k=1}^{p+1}\omega_{k}\varphi^{p+1-k,q}_{a_1a_2},\\
&\mathcal{B}^{m,n+1}_{a_1,a_2}
 =\OP{B}\left[\{\varphi^{m,k}_{a_1}(y_{a_2})\}_{k=0}^n\right]
 =\sum_{k=1}^{n+1}\omega_{k}\varphi^{m,n+1-k}(y_{a_1},y_{a_2})
 =\sum_{k=1}^{n+1}\omega_{k}\varphi^{m,n+1-k}_{a_1a_2}.
\end{split}
\end{equation}
The maps can now be expressed as
\begin{equation}
\begin{split}
&
\left.(\partial_{y_1}+\alpha_1)\varphi^{p+1,q}_{a_2}(y_1)\right|_{y_1=y_r}=-\alpha_1\mathcal{B}^{p+1,q}_{a_2,r},\quad
\left.(\partial_{y_2}+\alpha_2)\varphi^{m,n+1}_{a_1}(y_2)\right|_{y_1=y_t}=-\alpha_2\mathcal{B}^{m,n+1}_{a_1,t},\\
&
\left.(\partial_{y_1}-\alpha_1)\varphi^{p+1,q}_{a_2}(y_1)\right|_{y_1=y_l}=+\alpha_1\mathcal{B}^{p+1,q}_{a_2,l},\quad 
\left.(\partial_{y_2}-\alpha_2)\varphi^{m,n+1}_{a_1}(y_2)\right|_{y_1=y_b}=+\alpha_2\mathcal{B}^{m,n+1}_{a_1,b}.
\end{split}
\end{equation}

\emph{Corner conditions}--To solve the IBVP, we need to look at the discrete 
form of the corner conditions developed in Sec.~\ref{sec:CT-CQ-CC}. On the 
reference domain considered for the numerical implementation, we express the 
corner condition present in~\eqref{eq:cc-cq} as
\begin{equation}
\sqrt{\beta_1}\sqrt{\beta_2}\partial_{y_2}\partial_{y_1} u(x_1,x_2,t) = 
\left.e^{-i\pi/4}{e^{-i\pi/4}}\partial^{1/2}_{\tau_1}
\partial^{1/2}_{\tau_2}\varphi(x_1,x_2,\tau_1,\tau_2)\right|_{\tau_1=\tau_2=t}.
\end{equation}
At the discrete level, the corner condition at $\Gamma_{rt}$ takes the form
\begin{equation*}
\partial_{y_2}\partial_{y_1} u^{j+1} 
=\sqrt{\frac{\rho}{\beta_1}}e^{-i\pi/4}\sqrt{\frac{\rho}{\beta_2}}e^{-i\pi/4} 
\sum_{k=0}^{j+1}\omega_{k} \left[ \sum_{k'=0}^{j+1}\omega_{k'}\varphi^{j+1-k,j+1-k'}(y_r,y_t)\right]
=\alpha_1\alpha_2
\sum_{k=0}^{j+1}\omega_{k} \left[
\sum_{k'=0}^{j+1}\omega_{k'}\varphi^{j+1-k,j+1-k'}(y_r,y_t)\right].
\end{equation*}
In order to separate the current time-step value, we proceed in the following
manner:
\begin{equation}
\begin{split}
\partial_{y_2}\partial_{y_1} u^{j+1} 
&=\alpha_1\alpha_2
\sum_{k=0}^{j+1}\omega_{k} \left[ \varphi^{j+1-k,j+1}(y_r,y_t)
+\sum_{k'=1}^{j+1}\omega_{k'}\varphi^{j+1-k,j+1-k'}(y_r,y_t)\right]\\
&=\alpha_1\alpha_2\left[
\sum_{k=0}^{j+1}\omega_{k}\varphi^{j+1-k,j+1}(y_r,y_t)
+\sum_{k=0}^{j+1}\sum_{k'=1}^{j+1}\omega_{k}\omega_{k'}\varphi^{j+1-k,j+1-k'}(y_r,y_t)\right]\\
&=\alpha_1\alpha_2\left[ \varphi^{j+1,j+1}(y_r,y_t) + \mathcal{B}_r^{j+1}(y_t)
+\mathcal{B}_t^{j+1}(y_r) 
+\sum_{k=1}^{j+1}\sum_{k'=1}^{j+1}
\omega_{k}\omega_{k'}\varphi^{j+1-k,j+1-k'}(y_r,y_t)\right].
\end{split}
\end{equation}
Plugging-in the value for the history functions from~\eqref{eq:maps-cq-bdf1} in the above 
equation, we obtain
\begin{equation}\label{eq:cc-cq-bdf1}
\begin{split}
&
\left.(\partial_{y_2}+\alpha_2)(\partial_{y_1}+\alpha_1)u^{j+1}\right|_{(y_r,y_t)}=+\alpha_1\alpha_2\mathcal{C}^{j+1}_{rt},\quad
\left.(\partial_{y_2}-\alpha_2)(\partial_{y_1}+\alpha_1)u^{j+1}\right|_{(y_r,y_b)}=-\alpha_1\alpha_2\mathcal{C}^{j+1}_{br},\\
&
\left.(\partial_{y_2}+\alpha_2)(\partial_{y_1}-\alpha_1)u^{j+1}\right|_{(y_l,y_t)}=-\alpha_1\alpha_2\mathcal{C}^{j+1}_{tl},\quad
\left.(\partial_{y_2}-\alpha_2)(\partial_{y_1}-\alpha_1)u^{j+1}\right|_{(y_l,y_b)}=+\alpha_1\alpha_2\mathcal{C}^{j+1}_{lb},
\end{split}
\end{equation}
where
\begin{equation*}
\mathcal{C}^{j+1}_{a_1a_2}=\sum_{k=1}^{j+1}\sum_{k'=1}^{j+1}\omega_{k}\omega_{k'}\varphi^{j+1-k,j+1-k'}(y_{a_1},y_{a_2})
,\quad a_1\in\{r,l\},\quad a_2\in\{t,b\}.
\end{equation*}

\subsubsection{CQ--TR}
The discrete scheme is said to be `CQ--TR' if the underlying one-step method
in the CQ scheme is TR. Let $(\omega_k)_{k\in\field{N}_0}$ denote the
corresponding quadrature weights, then 
\begin{equation}
\omega_j=
\begin{cases}
C_{j/2},      & j\;\text{even},\\
-C_{(j-1)/2}, & j\;\text{odd},\\
\end{cases}
\quad\text{where}\quad C_n=\frac{1\cdot 3\cdots(2n-1)}{n!2^n}.
\end{equation}
The quadrature weights can also be
generated using the following recurrence relation~\cite{SV2023}
\begin{equation}
(j+1)\omega_{j+1}=(j-1)\omega_{j-1}-\omega_j,\quad j\geq1,
\end{equation}
with $\omega_0=1,\omega_1=-1$. The resulting discretization scheme for 
the DtN map described in~\eqref{eq:maps-cq} in terms of
$\varphi^{j,k}(y_1,y_2)$ stated as a Robin-type boundary condition (consistent
with the staggered samples of the interior field) reads as
\begin{equation}
\begin{split}
&\partial_{y_1}v^{j+1}+\alpha_1 v^{j+1}=-\frac{\alpha_1}{2}\left(\OP{B}
\left[\{\varphi^{k,j+1}(y_r,y_2)\}_{k=0}^j\right]
+ \OP{B}\left[\{\varphi^{k,j}(y_r,y_2)\}_{k=0}^{j-1}\right]\right),\\
&\partial_{y_1}v^{j+1}-\alpha_1 v^{j+1}=+\frac{\alpha_1}{2}\left(\OP{B}
\left[\{\varphi^{k,j+1}(y_l,y_2)\}_{k=0}^j\right]
+ \OP{B}\left[\{\varphi^{k,j}(y_l,y_2)\}_{k=0}^{j-1}\right]\right),\\
&\partial_{y_2}v^{j+1}+\alpha_2 v^{j+1}=-\frac{\alpha_2}{2}\left(\OP{B}
\left[\{\varphi^{k,j+1}(y_1,y_t)\}_{k=0}^j\right]
+ \OP{B}\left[\{\varphi^{k,j}(y_1,y_t)\}_{k=0}^{j-1}\right]\right),\\
&\partial_{y_2}v^{j+1}-\alpha_2 v^{j+1}=+\frac{\alpha_2}{2}\left(\OP{B}
\left[\{\varphi^{k,j+1}(y_1,y_b)\}_{k=0}^j\right]
+ \OP{B}\left[\{\varphi^{k,j}(y_1,y_b)\}_{k=0}^{j-1}\right]\right),\\
\end{split}
\end{equation}
where the history operator $\OP{B}$ is defined in the similar manner as in the case of BDF1.
Let $ \mathcal{B}^{j+1/2}=(\mathcal{B}^{j+1}+\mathcal{B}^{j})/2$, then the maps 
can be compactly written as
\begin{equation}\label{eq:maps-cq-tr}
\begin{split}
&
\partial_{y_1}v^{j+1}+\alpha_1v^{j+1}=-\alpha_1\mathcal{B}^{j+1/2}_r(y_2),\quad
\partial_{y_2}v^{j+1}+\alpha_2 v^{j+1}=-\alpha_2\mathcal{B}^{j+1/2}_t(y_1),\\
&
\partial_{y_1}v^{j+1}-\alpha_1 v^{j+1}=+\alpha_1\mathcal{B}^{j+1/2}_l(y_2),\quad
\partial_{y_2}v^{j+1}-\alpha_2 v^{j+1}=+\alpha_2\mathcal{B}^{j+1/2}_b(y_1).
\end{split}
\end{equation}
The TR-based discretization of the IVPs~\eqref{eq:ivps-cq} for the auxiliary 
function reads as follows: 
\begin{equation}
\begin{split}
&(y_1,y_2)\in{\Gamma}^{\text{ref.}}_b\cup{\Gamma}^{\text{ref.}}_t:\quad
\left\{
\begin{aligned}
&i\frac{\varphi^{p+1,q}_{a_2}-\varphi^{p,q}_{a_2}}{\Delta t}
+\beta_1\partial_{y_1}^2\frac{\varphi^{p+1,q}_{a_2}+\varphi^{p,q}_{a_2} }{2}=0
\implies
-\alpha^{-2}_1\partial_{y_1}^2\varphi^{p+1/2,q}_{a_2}+\varphi^{p+1/2,q}_{a_2}
=\varphi^{p,q}_{a_2},\\
&p=q,q+1,\ldots,j,\quad q = 0,1,\ldots,j;
\end{aligned}\right.\\
&(y_1,y_2)\in\Gamma^{\text{ref.}}_l\cup\Gamma^{\text{ref.}}_r:\quad
\left\{\begin{aligned}
&i\frac{\varphi^{m,n+1}_{a_1}-\varphi^{m,n}_{a_1}}{\Delta t}
+\beta_2\partial_{y_2}^2\frac{\varphi^{m,n+1}_{a_1}+\varphi^{m,n}_{a_1}}{2}=0
\implies
-\alpha^{-2}_2\partial_{y_2}^2\varphi^{m,n+1/2}_{a_1}+\varphi^{m,n+1/2}_{a_1}
=\varphi^{m,n}_{a_1},\\
&n=m,m+1,\ldots,j,\quad m = 0,1,\ldots,j.
\end{aligned}\right.
\end{split}
\end{equation}
The TR-based discretization of the DtN map described in~\eqref{eq:maps-cq-auxi} stated
as Robin-type boundary condition reads as:
\begin{equation}
\begin{split}
& \left.(\partial_{y_1}+\alpha_1)\varphi^{p+1/2,q}_{a_2}(y_1)\right|_{y_1=y_r} 
=-\frac{\alpha_1}{2}\left(\OP{B}\left[\{\varphi^{k,q}_{a_2}(y_r)\}_{k=0}^p\right]
+ \OP{B}\left[\{\varphi^{k,q}_{a_2}(y_r)\}_{k=0}^{p-1}\right]\right),\\
& \left.(\partial_{y_1}-\alpha_1)\varphi^{p+1/2,q}_{a_2}(y_1)\right|_{y_1=y_l} 
=+\frac{\alpha_1}{2}\left(\OP{B}\left[\{\varphi^{k,q}_{a_2}(y_l)\}_{k=0}^p\right]
+ \OP{B}\left[\{\varphi^{k,q}_{a_2}(y_l)\}_{k=0}^{p-1}\right]\right),\\
& \left.(\partial_{y_2}+\alpha_2)\varphi^{m,n+1/2}_{a_1}(y_2)\right|_{y_1=y_t} 
=-\frac{\alpha_2}{2}\left(\OP{B}\left[\{\varphi^{m,k}_{a_1}(y_t)\}_{k=0}^n\right]
+ \OP{B}\left[\{\varphi^{m,k}_{a_1}(y_t)\}_{k=0}^{p-1}\right]\right),\\
& \left.(\partial_{y_2}-\alpha_2)\varphi^{m,n+1/2}_{a_1}(y_2)\right|_{y_1=y_b} 
=+\frac{\alpha_2}{2}\left(\OP{B}\left[\{\varphi^{m,k}_{a_1}(y_b)\}_{k=0}^m\right]
+ \OP{B}\left[\{\varphi^{m,k}_{a_1}(y_b)\}_{k=0}^{p-1}\right]\right).
\end{split}
\end{equation}
where  operator $\OP{B}$ is defined as
\begin{equation}
\mathcal{B}^{p+1,q}_{a_2,a_1}
 =\OP{B}\left[\{\varphi^{k,q}_{a_2}(y_{a_1})\}_{k=0}^p\right]
 =\sum_{k=1}^{p+1}\omega_{k}\varphi^{p+1-k,q}_{a_1a_2},\quad
\mathcal{B}^{m,n+1}_{a_1,a_2}
 =\OP{B}\left[\{\varphi^{m,k}_{a_1}(y_{a_2})\}_{k=0}^n\right]
 =\sum_{k=1}^{n+1}\omega_{k}\varphi^{m,n+1-k}_{a_1a_2}.
\end{equation}
so that the maps can now be expressed as
\begin{equation}
\begin{split}
& 
\left.(\partial_{y_1}+\alpha_1)\varphi^{p+1/2,q}_{a_2}(y_1)\right|_{y_1=y_r}
=-\alpha_1\mathcal{B}^{p+1/2,q}_{a_2,r},\quad
\left.(\partial_{y_2}+\alpha_2)\varphi^{m,n+1/2}_{a_1}(y_2)\right|_{y_1=y_t}
=-\alpha_2\mathcal{B}^{m,n+1/2}_{a_1,t},\\
& 
\left.(\partial_{y_1}-\alpha_1)\varphi^{p+1/2,q}_{a_2}(y_1)\right|_{y_1=y_l}
=+\alpha_1\mathcal{B}^{p+1/2,q}_{a_2,l},\quad
\left.(\partial_{y_2}-\alpha_2)\varphi^{m,n+1/2}_{a_1}(y_2)\right|_{y_1=y_b}
=+\alpha_2\mathcal{B}^{m,n+1/2}_{a_1,b}.
\end{split}
\end{equation}
\emph{Corner conditions}--Following along the similar lines as that of BDF1, we can summarize 
the corner conditions consistent with staggered samples of the field as:
\begin{equation}
\begin{split}
&
\left.(\partial_{y_2}+\alpha_2)(\partial_{y_1}+\alpha_1)v^{j+1}\right|_{(y_r,y_t)}=+\alpha_1\alpha_2\mathcal{C}^{j+1/2}_{rt},\quad
\left.(\partial_{y_2}-\alpha_2)(\partial_{y_1}+\alpha_1)v^{j+1}\right|_{(y_r,y_b)}=-\alpha_1\alpha_2\mathcal{C}^{j+1/2}_{br},\\
&
\left.(\partial_{y_2}+\alpha_2)(\partial_{y_1}-\alpha_1)v^{j+1}\right|_{(y_l,y_t)}=-\alpha_1\alpha_2\mathcal{C}^{j+1/2}_{tl},\quad
\left.(\partial_{y_2}-\alpha_2)(\partial_{y_1}-\alpha_1)v^{j+1}\right|_{(y_l,y_b)}=+\alpha_1\alpha_2\mathcal{C}^{j+1/2}_{lb},\\
\end{split}
\end{equation}
where 
\begin{equation*}
\mathcal{C}^{j+1/2}_{a_1a_2}=\frac{1}{2}\left[
\sum_{k=1}^{j+1}\sum_{k'=1}^{j+1}\omega_{k}\omega_{k'}\varphi^{j+1-k,j+1-k'}(y_{a_1},y_{a_2})
+\sum_{k=1}^{j}\sum_{k'=1}^{j}\omega_{k}\omega_{k'}\varphi^{j-k,j-k'}(y_{a_1},y_{a_2})
\right]. 
\end{equation*}

\subsubsection{NP--BDF1}
The discrete scheme for the boundary maps is said to be `NP--BDF1' if the 
temporal derivatives in the realization of the boundary maps, dubbed as `NP', 
developed in Sec.~\ref{sec:CT-NP} are discretized using the one-step method 
BDF1. To achieve this, we start by writing the time-discrete form of the DtN 
maps present in~\eqref{eq:maps-pade} in the reference domain as
\begin{equation}
\begin{split}
&
\sqrt{\beta_1}\partial_{y_1}u^{j+1} \pm e^{-i\pi/4}\left[ b_0u^{j+1} -
\sum_{k=1}^M b_k\varphi^{j+1,j+1}_{k,a_1}\right]=0,\\
& 
\sqrt{\beta_2}\partial_{y_2}u^{j+1} \pm e^{-i\pi/4}\left[ b_0u^{j+1} -
\sum_{k=1}^M b_k\varphi^{j+1,j+1}_{k,a_2}\right]=0,
\end{split}
\end{equation}
where we recall $a_1\in\{r,l\}$ and $a_2\in\{t,b\}$. In order to turn these
equations into a Robin-type boundary condition, we need to compute 
the discrete samples $\varphi^{j+1,j+1}_{k,a_1}$ and $\varphi^{j+1,j+1}_{k,a_2}$ 
by considering the ODEs established for these auxiliary function earlier. For 
the temporal discretization of these ODEs, we try to match the time discretization 
scheme to that of the interior scheme to avoid any kind of instability as a 
result of incompatible discretization of the boundary maps. In the following, 
we discuss the BDF1-based discretization of~\eqref{eq:ode-auxi-npade1}: 
Setting $\rho=1/\Delta t$ and for $k=1,2,\ldots,M$, we have
\begin{equation}
\left\{\begin{aligned}
&\frac{\varphi^{j+1,j+1}_{k,a_1}-\varphi^{j,j+1}_{k,a_1}}{\Delta t}
+\eta^2_k\varphi^{j+1,j+1}_{k,a_1}
=\varphi^{j+1,j+1}\implies
\varphi^{j+1,j+1}_{k,a_1} 
=\frac{\rho}{(\rho+\eta^2_k)}\varphi^{j,j+1}_{k,a_1}
+\frac{1}{(\rho+\eta^2_k)}u^{j+1},\\
&\frac{\varphi^{j+1,j+1}_{k,a_2}-\varphi^{j+1,j}_{k,a_2}}{\Delta t} 
+\eta^2_k\varphi^{j+1,j+1}_{k,a_2}
=\varphi^{j+1,j+1}\implies
\varphi^{j+1,j+1}_{k,a_2} 
=\frac{\rho}{(\rho+\eta^2_k)}\varphi^{j+1,j}_{k,a_2}
+\frac{1}{(\rho+\eta^2_k)}u^{j+1}.
\end{aligned}\right.
\end{equation}
Introduce the scaling $\ovl{\eta}_k=\eta_k/\sqrt{\rho}$ so that
\begin{equation}
\varphi^{j+1,j+1}_{k,a_1} 
=\frac{1}{(1+\ovl{\eta}^2_k)}\varphi^{j,j+1}_{k,a_1}
+\frac{1/\rho}{(1+\ovl{\eta}^2_k)}u^{j+1},\quad
\varphi^{j+1,j+1}_{k,a_2} 
=\frac{1}{(1+\ovl{\eta}^2_k)}\varphi^{j+1,j}_{k,a_2}
+\frac{1/\rho}{(1+\ovl{\eta}^2_k)}u^{j+1}.
\end{equation}
Next, we plug-in the values of $\varphi^{j+1,j+1}_{k,a_1}$ and $\varphi^{j+1,j+1}_{k,a_2}$ 
in the boundary maps to obtain
\begin{equation}
\begin{split}
& \partial_{y_1} u^{j+1}\pm\alpha_1\left[
\ovl{b}_0 +\frac{1}{\rho}\sum_{k=1}^M\Gamma_{k}\right]u^{j+1}
=\mp\alpha_1\sum_{k=1}^M \Gamma_{k}\varphi_{k,a_1}^{j,j+1},\\
& \partial_{y_2} u^{j+1}\pm\alpha_2\left[
\ovl{b}_0+\frac{1}{\rho}\sum_{k=1}^M\Gamma_{k}\right]u^{j+1}
=\mp\alpha_2\sum_{k=1}^M \Gamma_{k}\varphi_{k,a_2}^{j+1,j},
\end{split}
\end{equation}
where 
\begin{equation}
\Gamma_{k} =-{\ovl{b}_{k}}/{(1+\ovl{\eta}^2_{k})},
\quad \ovl{b}_0=b_0/\sqrt{\rho},\quad \ovl{b}_k =b_k/\sqrt{\rho},
\quad k=1,2,\ldots,M.
\end{equation}
The maps can now be compactly written as a Robin-type boundary condition as
\begin{equation}
\partial_{y_1} u^{j+1}\pm\alpha_1\varpi u^{j+1}=\mp\alpha_1\mathcal{B}^{j+1}_{a_1}(y_2),\quad
\partial_{y_2} u^{j+1}\pm\alpha_2\varpi u^{j+1}=\mp\alpha_2\mathcal{B}^{j+1}_{a_2}(y_1),
\end{equation}
where
\begin{equation}\label{eq:varpi-defn}
\varpi=\ovl{b}_0 + \frac{1}{\rho}\sum_{k=1}^M\Gamma_{k}
=\ovl{b}_0-\frac{1}{\rho}\sum_{k=1}^M\frac{\ovl{b}_k}{1+\ovl{\eta}^2_k}=\sqrt{\rho}R_M(\rho),
\end{equation}
and
\begin{equation}
\mathcal{B}^{j+1}_{a_1}(y_2)=\sum_{k=1}^M \Gamma_{k}\varphi_{k,a_1}^{j,j+1},\quad 
\mathcal{B}^{j+1}_{a_2}(y_1)=\sum_{k=1}^M \Gamma_{k}\varphi_{k,a_2}^{j+1,j},
\end{equation}
are the so-called history functions in this setting. To obtain the discrete samples 
$\varphi_{k,a_1}^{j,j+1}$ and $\varphi_{k,a_2}^{j+1,j}$
needed to compute the history functions, we employ a Legendre-Galerkin method to solve
the IVPs listed in~\eqref{eq:ivp-pade-auxi} by formulating the boundary conditions
present in~\eqref{eq:maps-pade-auxi} as a Robin-type boundary map. The BDF1 based 
discretization of the IVPs take the form
\begin{equation}
\begin{split}
&(y_1,y_2)\in{\Gamma}^{\text{ref.}}_b\cup{\Gamma}^{\text{ref.}}_t:\quad
i\frac{\varphi^{j+1,j}_{k,a_2}-\varphi^{j,j}_{k,a_2}}{\Delta t}
+\beta_1\partial_{y_1}^2\varphi^{j,j}_{k,a_2}=0\implies 
-\alpha^{-2}_1\partial_{y_1}^2\varphi^{j+1,j}_{k,a_2}+\varphi^{j+1,j}_{k,a_2}
=\varphi^{j,j}_{k,a_2},\\
& (y_1,y_2)\in\Gamma^{\text{ref.}}_l\cup\Gamma^{\text{ref.}}_r:\quad
i\frac{\varphi^{j,j+1}_{k,a_1}-\varphi^{j,j}_{k,a_1}}{\Delta t}
+\beta_1\partial_{y_2}^2\varphi^{j,j}_{k,a_1}=0\implies 
-\alpha^{-2}_1\partial_{y_2}^2\varphi^{j,j+1}_{k,a_1}+\varphi^{j,j+1}_{k,a_1}
=\varphi^{j,j}_{k,a_1}.
\end{split}
\end{equation}
The discrete form of the DtN maps~\eqref{eq:dtn-using-psi} for these IVPs can be written as:
\begin{equation}\label{eq:maps-pade-auxi-bdf1}
\begin{split}
& \sqrt{\beta_2}\partial_{y_2}\varphi_{k,a_1}^{j,j+1} \pm e^{-i\pi/4}\left[ 
b_0\varphi_{k,a_1}^{j,j+1} - \sum_{k'=1}^M b_{k'}\psi^{j,j+1}_{k,k',a_1,a_2}\right]=0,\\
& \sqrt{\beta_1}\partial_{y_1}\varphi_{k,a_2}^{j+1,j} \pm e^{-i\pi/4}\left[ 
b_0\varphi_{k,a_2}^{j+1,j} - \sum_{k'=1}^M b_{k'}\psi^{j+1,j}_{k,k',a_2,a_1}\right]=0.
\end{split}
\end{equation}
In order to obtain the history functions, we invoke BDF1 based discretization
of~\eqref{eq:ode-auxi-npade2} and~\eqref{eq:ode-auxi-npade3} to write
\begin{equation}\label{eq:diag-psi-bdf1}
\begin{split}
& \frac{\psi^{j+1,j+1}_{k,k',a_1a_2}-\psi^{j,j+1}_{k,k',a_1,a_2}}{\Delta t} 
 +\eta^2_{k}\psi^{j+1,j+1}_{k,k',a_1,a_2}
 =\varphi_{k',a_2}^{j+1,j+1}\implies 
 \psi^{j+1,j+1}_{k,k',a_1a_2} 
 =\frac{1}{(1+\ovl{\eta}^2_{k})}\psi^{j,j+1}_{k,k',a_1,a_2}
 +\frac{1/\rho}{(1+\ovl{\eta}^2_{k})}\varphi_{k',a_2}^{j+1,j+1},\\
& \frac{\psi^{j,j+1}_{k,k',a_1,a_2}-\psi^{j,j}_{k,k',a_1,a_2}}{\Delta t} 
 +\eta^2_{k'}\psi^{j,j+1}_{k,k',a_1,a_2}
 =\varphi_{k,a_1}^{j,j+1}\implies 
 \psi^{j,j+1}_{k,k',a_1a_2} 
  =\frac{1}{(1+\ovl{\eta}^2_{k'})}\psi^{j,j}_{k,k',a_1,a_2}
  +\frac{1/\rho}{(1+\ovl{\eta}^2_{k'})}\varphi_{k,a_1}^{j,j+1}.
\end{split}
\end{equation}
Plugging-in the values of $\psi^{j,j+1}_{k,k',a_1,a_2}$ and $\psi^{j+1,j}_{k,k',a_2,a_1}$
in~\eqref{eq:maps-pade-auxi-bdf1}, we obtain
\begin{equation}
\begin{split}
&\partial_{y_2}\varphi_{k,a_1}^{j,j+1}\pm\alpha_2\left[
\ovl{b}_0+\frac{1}{\rho}\sum_{k'=1}^M\Gamma_{k'}\right]\varphi_{k,a_1}^{j,j+1}
=\mp\alpha_2\sum_{k'=1}^M \Gamma_{k'}\psi_{k,k',a_1,a_2}^{j,j},\\
&\partial_{y_1}\varphi_{k,a_2}^{j+1,j}\pm \alpha_1\left[
\ovl{b}_0+\frac{1}{\rho}\sum_{k'=1}^M\Gamma_{k'}\right]\varphi_{k,a_2}^{j+1,j}
=\mp\alpha_1\sum_{k'=1}^M \Gamma_{k'}\psi_{k,k',a_2,a_1}^{j,j}.
\end{split}
\end{equation}
Once again, we can express the maps compactly as
\begin{equation}
\partial_{y_2}\varphi_{k,a_1}^{j,j+1}\pm\alpha_2\varpi\varphi_{k,a_1}^{j,j+1}
=\mp\alpha_2\mathcal{B}^{j+1}_{k,a_1,a_2},\quad
\partial_{y_1}\varphi_{k,a_2}^{j+1,j}\pm\alpha_1\varpi\varphi_{k,a_2}^{j+1,j}
=\mp\alpha_1\mathcal{B}^{j+1}_{k,a_2,a_1},
\end{equation}
where the history functions in this setting are given by
\begin{equation}
\mathcal{B}^{j+1}_{k,a_1,a_2}=\sum_{k'=1}^M \Gamma_{k'}\psi_{k,k',a_1,a_2}^{j,j},
\quad \mathcal{B}^{j+1}_{k,a_2,a_1}=\sum_{k'=1}^M \Gamma_{k'}\psi_{k,k',a_2,a_1}^{j,j}
=\sum_{k'=1}^M \Gamma_{k'}\psi_{k',k,a_1,a_2}^{j,j}.
\end{equation}
\emph{Corner conditions}--To solve the IBVP, we need to develop the discrete form 
of the corner conditions developed in Sec.~\ref{sec:CT-NP}. On the reference 
domain considered for the numerical 
implementation, we express the corner conditions present in \eqref{eq:cc-npade} as
\begin{equation}
\begin{split}
&\left[\partial_{y_2}\partial_{y_1}+\alpha_2\ovl{b}_0\partial_{y_1}
+\alpha_1\ovl{b}_0\partial_{y_2}+\alpha_1\alpha_2\ovl{b}_0^2\right]u(x_{a_1},x_{a_2},t)
-\alpha_1\alpha_2\sum_{k=1}^M\sum_{k'=1}^M\ovl{b}_k\ovl{b}_{k'}\psi_{k,k',a_1,a_2}(t,t)=0.
\end{split}
\end{equation}
At the discrete level, the corner conditions take the form
\begin{equation}\label{eq:cc-bdf1-npade}
\begin{split}
&\left[\partial_{y_2}\partial_{y_1}+\alpha_2\ovl{b}_0\partial_{y_1}
+\alpha_1\ovl{b}_0\partial_{y_2}+\alpha_1\alpha_2\ovl{b}_0^2\right]u^{j+1}
-\alpha_1\alpha_2\sum_{k=1}^M\sum_{k'=1}^M\ovl{b}_k\ovl{b}_{k'}\psi^{j+1,j+1}_{k,k',a_1,a_2}=0.
\end{split}
\end{equation}
The diagonal-to-diagonal update of auxiliary fields $\psi$ is already computed 
in~\eqref{eq:diag-psi-bdf1} which can be simplified as
\begin{equation}
\begin{split}
\psi^{j+1,j+1}_{k,k',a_1a_2} 
&=\frac{1}{(1+\ovl{\eta}^2_{k})(1+\ovl{\eta}^2_{k'})}\psi^{j,j}_{k,k',a_1a_2}
+\frac{(1/\rho)^2}{(1+\ovl{\eta}^2_{k})(1+\ovl{\eta}^2_{k'})}u^{j+1}\\
&\quad+\frac{1/\rho}{(1+\ovl{\eta}^2_{k})(1+\ovl{\eta}^2_{k'})}\varphi_{k,a_1}^{j,j+1}
+\frac{1/\rho}{(1+\ovl{\eta}^2_{k})(1+\ovl{\eta}^2_{k'})}\varphi_{k',a_2}^{j+1,j}.
\end{split}
\end{equation}
The last term in~\eqref{eq:cc-bdf1-npade} can be expressed in terms of boundary sums as 
\begin{equation}
\begin{split}
\sum_{k}\sum_{k'}\ovl{b}_k\ovl{b}_{k'}\psi^{j+1,j+1}_{k,k',a_1a_2} 
&=\sum_{k}\sum_{k'}
\left[\frac{\ovl{b}_k\ovl{b}_{k'}}{(1+\ovl{\eta}^2_{k})(1+\ovl{\eta}^2_{k'})}\psi^{j,j}_{k,k',a_1a_2}\right]\\
&\quad+\frac{1}{\rho}\left[\sum_{k}\Gamma_{k}\right]\mathcal{B}^{j+1}_{a_1}(y_2)
+\frac{1}{\rho}\left[\sum_{k}\Gamma_{k}\right]\mathcal{B}^{j+1}_{a_2}(y_1)
+\frac{1}{\rho^2}\left[\sum_{k}\Gamma_{k}\right]^2 u^{j+1}.
\end{split}
\end{equation}
Plugging-in the value of $\psi^{j+1,j+1}_{k,k',a_1,a_2}$ in~\eqref{eq:cc-bdf1-npade} and 
collecting all the terms at $(j+1)$-th time-step on the left hand side, we
arrive at the expression 
\begin{equation}
\begin{split}
&\left[\partial_{y_2}\partial_{y_1}+\alpha_2\varpi\partial_{y_1}
+\alpha_1\varpi\partial_{y_2}+\alpha_1\alpha_2\ovl{b}_0^2\right]u^{j+1}
+\alpha_1\alpha_2[2\varpi (\varpi-\ovl{b}_0)-(\varpi-\ovl{b}_0)^2]u^{j+1}\\
&=\left[\partial_{y_2}\partial_{y_1}+\alpha_2\varpi\partial_{y_1}
+\alpha_1\varpi\partial_{y_2}+\alpha_1\alpha_2\ovl{b}_0^2\right]u^{j+1}
+\alpha_1\alpha_2[\varpi^2-\ovl{b}_0^2]u^{j+1}\\
&=(\partial_{y_2}+\alpha_2\varpi)(\partial_{y_1}+\alpha_1\varpi)u^{j+1}.
\end{split}
\end{equation}
With this observation, the corner conditions turns out to be 
\begin{equation}
\begin{split}
&\left.(\partial_{y_2}+\alpha_2\varpi)(\partial_{y_1}+\alpha_1\varpi)u^{j+1}\right|_{(y_r,y_t)}=+\alpha_1\alpha_2\mathcal{C}^{j+1}_{rt},\quad
\left.(\partial_{y_2}-\alpha_2\varpi)(\partial_{y_1}+\alpha_1\varpi)u^{j+1}\right|_{(y_r,y_b)}=-\alpha_1\alpha_2\mathcal{C}^{j+1}_{br},\\
&\left.(\partial_{y_2}+\alpha_2\varpi)(\partial_{y_1}-\alpha_1\varpi)u^{j+1}\right|_{(y_l,y_t)}=-\alpha_1\alpha_2\mathcal{C}^{j+1}_{tl},\quad
\left.(\partial_{y_2}-\alpha_2\varpi)(\partial_{y_1}-\alpha_1\varpi)u^{j+1}\right|_{(y_l,y_b)}=+\alpha_1\alpha_2\mathcal{C}^{j+1}_{lb},\\
\end{split}
\end{equation}
where
\begin{equation*}
\mathcal{C}^{j+1}_{a_1a_2}=\sum_{k}\sum_{k'}
\left[
\frac{\ovl{b}_k\ovl{b}_{k'}}{(1+\ovl{\eta}^2_{k})(1+\ovl{\eta}^2_{k'})}\psi^{j,j}_{k,k',a_1a_2}
\right],\quad a_1\in\{r,l\},\quad a_2\in\{t,b\}.
\end{equation*}
\subsubsection{NP--TR}
The discrete scheme for the boundary maps is said to be `NP--TR' if the 
temporal derivatives in the realization of the boundary maps, dubbed as `NP', 
developed in Sec.~\ref{sec:CT-NP} are discretized using the one-step method TR. 
The time-discrete form of DtN maps consistent with the staggered samples of field
present in~\eqref{eq:maps-pade} on the reference domain reads as 
\begin{equation}
\begin{split}
& \sqrt{\beta_1}\partial_{y_1}v^{j+1} \pm e^{-i\pi/4}\left[ b_0 v^{j+1} -
    \sum_{k=1}^M b_k\varphi^{j+1/2,j+1/2}_{k,a_1}\right]=0,\\
& \sqrt{\beta_2}\partial_{y_2}v^{j+1} \pm e^{-i\pi/4}\left[ b_0 v^{j+1} -
    \sum_{k=1}^M b_k\varphi^{j+1/2,j+1/2}_{k,a_2}\right]=0.
\end{split}
\end{equation}
In order to obtain a Robin-type boundary condition, we need to compute the 
discrete values of $\varphi^{j+1/2,j+1/2}_{k,a_1}$ and 
$\varphi^{j+1/2,j+1/2}_{k,a_2}$ as in the last section. 
The TR-based discretization of~\eqref{eq:ode-auxi-npade1} reads as
\begin{equation}
\begin{split}
& 
\frac{\varphi^{j+1,j+1}_{k,a_1}-\varphi^{j,j+1}_{k,a_1}}{\Delta t} 
+\eta^2_k\frac{\varphi^{j+1,j+1}_{k,a_1}+\varphi^{j,j+1}_{k,a_1}}{2}
=\frac{\varphi_{a_1}^{j+1,j+1}+\varphi_{a_1}^{j,j+1}}{2},\\
& 
\frac{\varphi^{j+1,j+1}_{k,a_2}-\varphi^{j+1,j}_{k,a_2}}{\Delta t} 
+\eta^2_k\frac{\varphi^{j+1,j+1}_{k,a_2}+\varphi^{j+1,j}_{k,a_2}}{2}
=\frac{\varphi_{a_2}^{j+1,j+1}+\varphi_{a_2}^{j+1,j}}{2},
\end{split}
\end{equation}
which simplifies to 
\begin{equation}
\begin{split}
& 
\varphi^{j+1,j+1}_{k,a_1} 
=\left(\frac{1-\ovl{\eta}^2_k}{1+\ovl{\eta}^2_k}\right)\varphi^{j,j+1}_{k,a_1}
+\frac{2/\rho}{(1+\ovl{\eta}^2_k)}\left[v_{a_1}^{j+1}
+\frac{\varphi_{a_1}^{j,j+1}-\varphi_{a_1}^{j,j}}{2}\right],\\
& 
\varphi^{j+1,j+1}_{k,a_2} 
=\left(\frac{1-\ovl{\eta}^2_k}{1+\ovl{\eta}^2_k}\right)\varphi^{j+1,j}_{k,a_2}
+\frac{2/\rho}{(1+\ovl{\eta}^2_k)}\left[v_{a_2}^{j+1}
+\frac{\varphi_{a_2}^{j+1,j}-\varphi_{a_2}^{j,j}}{2}\right].
\end{split}
\end{equation}
The staggered samples then work out to be
\begin{equation}
\begin{split}
& \varphi^{j+1/2,j+1/2}_{k,a_1} =\frac{1}{2}\left[\left(
\frac{1-\ovl{\eta}^2_k}{1+\ovl{\eta}^2_k}\right)\varphi^{j,j+1}_{k,a_1}
+\varphi^{j,j}_{k,a_1}\right] +\frac{1/\rho}{(1+\ovl{\eta}^2_k)}\left[v_{a_1}^{j+1}
+\frac{\varphi_{a_1}^{j,j+1}-\varphi_{a_1}^{j,j}}{2}\right],\\
& \varphi^{j+1/2,j+1/2}_{k,a_2} =\frac{1}{2}\left[\left(
\frac{1-\ovl{\eta}^2_k}{1+\ovl{\eta}^2_k}\right)\varphi^{j+1,j}_{k,a_2}+\varphi^{j,j}_{k,a_2}\right]
+\frac{1/\rho}{(1+\ovl{\eta}^2_k)}\left[v_{a_2}^{j+1}
+\frac{\varphi_{a_2}^{j+1,j}-\varphi_{a_2}^{j,j}}{2}\right].
\end{split}
\end{equation}
Plugging-in the values 
of $\varphi^{j+1/2,j+1/2}_{k,a_1}$ and $\varphi^{j+1/2,j+1/2}_{k,a_2}$, the maps can be 
compactly written as
\begin{equation}
\partial_{y_1} v^{j+1}\pm\alpha_1\varpi
v^{j+1}=\mp\alpha_1\mathcal{B}^{j+1/2}_{a_1}(y_2),\quad
\partial_{y_2} v^{j+1}\pm\alpha_2\varpi v^{j+1}
=\mp\alpha_2\mathcal{B}^{j+1/2}_{a_2}(y_1),
\end{equation}
where the history functions in this setting are given by
\begin{equation}\label{eq:hist-func-np-tr-segment}
\begin{split}
& \mathcal{B}^{j+1/2}_{a_1}(y_2)=\sum_{k=1}^M \left(\frac{-\ovl{b}_k}{2}
\left[\left(\frac{1-\ovl{\eta}^2_k}{1+\ovl{\eta}^2_k}\right)
{\varphi}^{j,j+1}_{k,a_1}+{\varphi}^{j,j}_{k,a_1}\right]+\frac{\Gamma_k}{\rho}
\left[\frac{\varphi_{a_1}^{j,j+1}-u_{a_1}^{j,j}}{2}\right]\right),\\
& \mathcal{B}^{j+1/2}_{a_2}(y_1)=\sum_{k=1}^M \left(\frac{-\ovl{b}_k}{2}
\left[\left(\frac{1-\ovl{\eta}^2_k}{1+\ovl{\eta}^2_k}\right)
{\varphi}^{j+1,j}_{k,a_2}+{\varphi}^{j,j}_{k,a_2}\right]+\frac{\Gamma_k}{\rho}
\left[\frac{\varphi_{a_2}^{j+1,j}-u_{a_1}^{j,j}}{2}\right]\right).
\end{split}
\end{equation}
To obtain the discrete samples $\varphi_{k,a_1}^{j,j+1}$ and $\varphi_{k,a_2}^{j+1,j}$
needed to compute the history functions, we employ a Legendre-Galerkin method to solve
the IVPs~\eqref{eq:ivp-pade-auxi} using a TR-based discretization as
\begin{equation}
\begin{split}
&(y_1,y_2)\in{\Gamma}^{\text{ref.}}_b\cup{\Gamma}^{\text{ref.}}_t:\quad
i\frac{\varphi^{j+1,j}_{k,a_2}-\varphi^{j,j}_{k,a_2}}{\Delta t}
+\beta_1\partial_{y_1}^2\varphi^{j+1/2,j}_{k,a_2}=0\implies
-\alpha^{-2}_1\partial_{y_1}^2\varphi^{j+1/2,j}_{k,a_2}+\varphi^{j+1/2,j}_{k,a_2}
=\varphi^{j,j}_{k,a_2},\\
&(y_1,y_2)\in\Gamma^{\text{ref.}}_l\cup\Gamma^{\text{ref.}}_r:\quad
i\frac{\varphi^{j,j+1}_{k,a_1}-\varphi^{j,j}_{k,a_1}}{\Delta t}
+\beta_1\partial_{y_2}^2\varphi^{j,j+1/2}_{k,a_1}=0\implies 
-\alpha^{-2}_1\partial_{y_2}^2\varphi^{j,j+1/2}_{k,a_1}+\varphi^{j,j+1/2}_{k,a_1}
=\varphi^{j,j}_{k,a_1}.
\end{split}
\end{equation}
The discrete form of the DtN maps~\eqref{eq:dtn-using-psi} for these IVPs can be written as:
\begin{equation}\label{eq:maps-pade-auxi-tr}
\begin{split}
& \sqrt{\beta_2}\partial_{y_2}\varphi_{k,a_1}^{j,j+1/2} \pm e^{-i\pi/4}\left[ 
b_0\varphi_{k,a_1}^{j,j+1/2} - \sum_{k'=1}^M b_{k'}\psi^{j,j+1/2}_{k,k',a_1,a_2}\right]=0,\\
& \sqrt{\beta_1}\partial_{y_1}\varphi_{k,a_2}^{j+1/2,j} \pm e^{-i\pi/4}\left[ 
b_0\varphi_{k,a_2}^{j+1/2,j} - \sum_{k'=1}^M b_{k'}\psi^{j+1/2,j}_{k,k',a_2,a_1}\right]=0.
\end{split}
\end{equation}
In order to arrive at the Robin-type boundary conditions, we invoke TR-based 
discretization of~\eqref{eq:ode-auxi-npade2} to write
\begin{equation}\label{eq:aux-psi-update}
\begin{split}
& \frac{\psi^{j,j+1}_{k,k',a_1a_2}-\psi^{j,j}_{k,k',a_1a_2}}{\Delta t} 
+\eta^2_{k'}\psi^{j,j+1/2}_{k,k',a_1a_2}
=\varphi_{k,a_1}^{j,j+1/2}\implies
\psi^{j,j+1}_{k,k',a_1a_2} 
=\frac{(1-\ovl{\eta}^2_{k'})}{(1+\ovl{\eta}^2_{k'})}\psi^{j,j}_{k,k',a_1a_2}
+\frac{2/\rho}{(1+\ovl{\eta}^2_{k'})}\varphi_{k,a_1}^{j,j+1/2},\\
&\frac{\psi^{j+1,j}_{k,k',a_2a_1}-\psi^{j,j}_{k,k',a_2a_1}}{\Delta t} 
+\eta^2_{k'}\psi^{j+1/2,j}_{k,k',a_2a_1}
=\varphi_{k,a_1}^{j+1/2,j}\implies
\psi^{j+1,j}_{k,k',a_2a_1} 
=\frac{(1-\ovl{\eta}^2_{k'})}{(1+\ovl{\eta}^2_{k'})}\psi^{j,j}_{k,k',a_2a_1}
+\frac{2/\rho}{(1+\ovl{\eta}^2_{k'})}\varphi_{k,a_2}^{j+1/2,j},
\end{split}
\end{equation}
so that the staggered samples work out to be
\begin{equation}\label{eq:diag-psi-tr}
\psi^{j,j+1/2}_{k,k',a_1a_2} 
=\frac{1}{(1+\ovl{\eta}^2_{k'})}\psi^{j,j}_{k,k',a_1a_2}
+\frac{1/\rho}{(1+\ovl{\eta}^2_{k'})}\varphi_{k,a_1}^{j,j+1/2},\quad
\psi^{j+1/2,j}_{k,k',a_2a_1} 
=\frac{1}{(1+\ovl{\eta}^2_{k'})}\psi^{j,j}_{k,k',a_2a_1}
+\frac{1/\rho}{(1+\ovl{\eta}^2_{k'})}\varphi_{k,a_2}^{j+1/2,j}.
\end{equation}
Plugging-in the values of $\psi^{j,j+1/2}_{k,k',a_1,a_2}$ and $\psi^{j+1/2,j}_{k,k',a_2,a_1}$
in~\eqref{eq:maps-pade-auxi-tr}, we get
\begin{equation}
\begin{split}
&\partial_{y_2}\varphi_{k,a_1}^{j,j+1/2}\pm\alpha_2\left[
\ovl{b}_0+\frac{1}{\rho}\sum_{k'=1}^M\Gamma_{k'}\right]\varphi_{k,a_1}^{j,j+1/2}
=\mp\alpha_2\sum_{k'=1}^M \Gamma_{k'}\psi_{k,k',a_1,a_2}^{j,j},\\
& \partial_{y_1}\varphi_{k,a_2}^{j+1/2,j}\pm \alpha_1\left[
\ovl{b}_0+\frac{1}{\rho}\sum_{k'=1}^M\Gamma_{k'}\right]\varphi_{k,a_2}^{j+1/2,j}
=\mp\alpha_1\sum_{k'=1}^M \Gamma_{k'}\psi_{k,k',a_2,a_1}^{j,j}.
\end{split}
\end{equation}
Finally, we can express these maps compactly as
\begin{equation}
\partial_{y_2}\varphi_{k,a_1}^{j,j+1/2}\pm\alpha_2\varpi\varphi_{k,a_1}^{j,j+1/2}
=\mp\alpha_2\mathcal{B}^{j+1}_{k,a_1,a_2},\quad
\partial_{y_1}\varphi_{k,a_2}^{j+1/2,j}\pm\alpha_1\varpi\varphi_{k,a_2}^{j+1/2,j}
=\mp\alpha_1\mathcal{B}^{j+1}_{k,a_2,a_1}.
\end{equation}
where the history functions are given by
\begin{equation}
\mathcal{B}^{j+1}_{k,a_1,a_2}=\sum_{k'=1}^M \Gamma_{k'}\psi_{k,k',a_1,a_2}^{j,j},
\quad \mathcal{B}^{j+1}_{k,a_2,a_1}=\sum_{k'=1}^M \Gamma_{k'}\psi_{k,k',a_2,a_1}^{j,j}
=\sum_{k'=1}^M \Gamma_{k'}\psi_{k',k,a_1,a_2}^{j,j},
\end{equation}
where the second result uses the fact that
$\psi_{k,k',a_2,a_1}^{j,j}=\psi_{k',k,a_1,a_2}^{j,j}$ (see Remark~\ref{rem:aux-func-psi}). 

Addressing the diagonal-to-diagonal update for the $\psi$-function in the order 
$\psi^{j,j}_{k,k',a_1a_2}\to\psi^{j,j+1}_{k,k',a_1a_2}\to\psi^{j+1,j+1}_{k,k',a_1a_2}$, 
we use~\eqref{eq:aux-psi-update} and the discrete version
of~\eqref{eq:ode-auxi-npade3} which reads as
\begin{equation}\label{eq:diag-to-diag-psi}
\begin{split}
& \frac{\psi^{j+1,j+1}_{k,k',a_1a_2}-\psi^{j,j+1}_{k,k',a_1a_2}}{\Delta t} 
+\eta^2_{k}\psi^{j+1/2,j+1}_{k,k',a_1a_2}
=\varphi_{k',a_2}^{j+1/2,j+1},\\
& \psi^{j+1,j+1}_{k,k',a_1a_2} 
=\frac{(1-\ovl{\eta}^2_{k})}{(1+\ovl{\eta}^2_{k})}\psi^{j,j+1}_{k,k',a_1a_2}
+\frac{2/\rho}{(1+\ovl{\eta}^2_{k})}\left[
\frac{\varphi_{k',a_2}^{j+1,j+1}+\varphi_{k',a_2}^{j,j+1}}{2}\right].
\end{split}
\end{equation}
To obtain the discrete samples $\varphi_{a_1}^{j,j+1}$ and $\varphi_{a_2}^{j+1,j}$
needed to compute the history functions in~\eqref{eq:hist-func-np-tr-segment}, 
we employ a Legendre-Galerkin method to solve
the IVPs listed in~\eqref{eq:ivps-cq} using a TR-based discretization as
\begin{equation}
\begin{split}
&(y_1,y_2)\in{\Gamma}^{\text{ref.}}_b\cup{\Gamma}^{\text{ref.}}_t:\quad 
i\frac{\varphi^{j+1,j}_{a_2}-\varphi^{j,j}_{a_2}}{\Delta t}
+\beta_1\partial_{y_1}^2\varphi^{j+1/2,j}_{a_2}=0\implies
-\alpha^{-2}_1\partial_{y_1}^2\varphi^{j+1/2,j}_{a_2}+\varphi^{j+1/2,j}_{a_2}
=\varphi^{j,j}_{a_2},\\
&(y_1,y_2)\in\Gamma^{\text{ref.}}_l\cup\Gamma^{\text{ref.}}_r:\quad
i\frac{\varphi^{j,j+1}_{a_1}-\varphi^{j,j}_{a_1}}{\Delta t}
+\beta_1\partial_{y_2}^2\varphi^{j,j+1/2}_{a_1}=0\implies
-\alpha^{-2}_1\partial_{y_2}^2\varphi^{j,j+1/2}_{a_1}+\varphi^{j,j+1/2}_{a_1}
=\varphi^{j,j}_{a_1}.
\end{split}
\end{equation}
The DtN maps for the IVPs above are given in~\eqref{eq:maps-cq-auxi} where 
we approximate the fractional derivative according to the novel Pad\'e scheme 
as discussed in Sec.~\ref{sec:CT-NP}. To this end, let us introduce the 
auxiliary functions $\psi_{k,a_2a_1}(\tau_1,\tau_2)$ and 
$\psi_{k,a_1a_2}(\tau_1,\tau_2)$ such that
\begin{equation}
(\partial_{\tau_1}+\eta^2_{k})\psi_{k,a_2a_1}(\tau_1,\tau_2)
=\varphi_{a_2}(x_1,\tau_1,\tau_2),\quad
(\partial_{\tau_2}+\eta^2_{k})\psi_{k,a_1a_2}(\tau_1,\tau_2)
=\varphi_{a_1}(x_2,\tau_1,\tau_2).
\end{equation}
As it turns out, the auxiliary functions, 
$\psi_{k,a_2a_1}(\tau_1,\tau_2)$ and $\psi_{k,a_1a_2}(\tau_1,\tau_2)$ become 
identical to $\varphi_{k,a_1}(x_2,\tau_1,\tau_2)$ 
and $\varphi_{k,a_2}(x_1,\tau_1,\tau_2)$, respectively 
(see Remark~\ref{rem:additional-aux-psi}). Following 
along the similar lines as above, we can write the discrete Robin-type 
boundary conditions for these IVPs as
\begin{equation}
\partial_{y_2}\varphi_{a_1}^{j,j+1/2}\pm\alpha_2\varpi\varphi_{a_1}^{j,j+1/2}
=\mp\alpha_2\mathcal{B}^{j+1}_{a_1,a_2},\quad
\partial_{y_1}\varphi_{a_2}^{j+1/2,j}\pm\alpha_1\varpi\varphi_{a_2}^{j+1/2,j}
=\mp\alpha_1\mathcal{B}^{j+1}_{a_2,a_1},
\end{equation}
where the history functions are given by
\begin{equation}
\mathcal{B}^{j+1}_{a_1,a_2} 
=\sum_{k=1}^M \Gamma_{k}\phi_{k,a_2}^{j,j}(y_{a_1}),\quad
\mathcal{B}^{j+1}_{a_2,a_1} 
=\sum_{k=1}^M\Gamma_{k}\phi_{k,a_1}^{j,j}(y_{a_2}).
\end{equation}
\begin{rem}\label{rem:additional-aux-psi}
From the definition of $\psi_{k,a_2a_1}(\tau_1,\tau_2)$ and 
$\psi_{k,a_1a_2}(\tau_1,\tau_2)$ it follows that
\begin{equation}
\begin{split}
\psi_{k,a_2a_1}(\tau_1,\tau_2)
&
=\int_0^{\tau_1}e^{-\eta^2_k(\tau_1-s)}\varphi(x_1,x_2,\tau_2,z)ds
=\varphi_{k,a_1}(x_{a_2},\tau_1,\tau_2),\quad a_2 \in \{b,t\},\\
\psi_{k,a_1a_2}(\tau_1,\tau_2)
& 
=\int_0^{\tau_2}e^{-\eta^2_k(\tau_2-s)}\varphi(x_1,x_2,s,\tau_1)ds
=\varphi_{k,a_2}(x_{a_1},\tau_1,\tau_2),\quad a_1 \in \{l,r\}.
\end{split}
\end{equation}
Note that the discrete samples of the auxiliary functions 
$\varphi_{k,a_1}(x_2,\tau_1,\tau_2)$ and $\varphi_{k,a_2}(x_1,\tau_1,\tau_2)$
are already being computed in connection with the Robin-type boundary condition
for the interior field.
\end{rem}
\emph{Corner conditions}--Following along the similar lines as discussed in the 
case of NP--BDF1, we can formulate the time-discrete
corner conditions for NP--TR as well. At the discrete level, we write the corner 
conditions consistent with staggered samples of the field present 
in~\eqref{eq:cc-npade} as
\begin{equation}\label{eq:cc-tr-npade}
\begin{split}
&\left[\partial_{y_2}\partial_{y_1}+\alpha_2\ovl{b}_0\partial_{y_1}
+\alpha_1\ovl{b}_0\partial_{y_2}+\alpha_1\alpha_2\ovl{b}_0^2\right]v^{j+1}
-\alpha_1\alpha_2\sum_{k=1}^M\sum_{k'=1}^M\ovl{b}_k\ovl{b}_{k'}\psi^{j+1/2,j+1/2}_{k,k',a_1,a_2}=0.
\end{split}
\end{equation}
The diagonal-to-diagonal update of auxiliary fields $\psi$ is already computed 
in~\eqref{eq:diag-psi-tr} which can be simplified as
\begin{equation}
\begin{split}
\psi^{j+1,j+1}_{k,k',a_1a_2} 
&=\frac{(1-\ovl{\eta}^2_{k})}{(1+\ovl{\eta}^2_{k})}\frac{(1-\ovl{\eta}^2_{k'})}{(1+\ovl{\eta}^2_{k'})}
\psi^{j,j}_{k,k',a_1a_2}\\
&\quad+\frac{(1-\ovl{\eta}^2_{k})}{(1+\ovl{\eta}^2_{k})}\frac{1/\rho}{(1+\ovl{\eta}^2_{k'})}
\left[\varphi_{k,a_1}^{j,j+1}+\varphi_{k,a_1}^{j,j}\right]
+\frac{1/\rho}{(1+\ovl{\eta}^2_{k})}\left(\frac{1-\ovl{\eta}^2_{k'}}{1+\ovl{\eta}^2_{k'}}\right)
\left[\varphi^{j+1,j}_{k',a_2}+\varphi^{j,j}_{k',a_2}\right],\\
&\quad+\frac{1/\rho}{(1+\ovl{\eta}^2_{k})}\frac{1/\rho}{(1+\ovl{\eta}^2_{k'})}
\left[u^{j+1}+u^j+\varphi_{a_2}^{j+1,j}+\varphi_{a_2}^{j,j+1}\right].
\end{split}
\end{equation}
The last term in~\eqref{eq:cc-tr-npade} can be expressed in terms of the history
functions~\eqref{eq:hist-func-np-tr-segment} to obtain
\begin{equation}
\begin{split}
\sum_{k}\sum_{k'}\ovl{b}_k\ovl{b}_{k'}\psi^{j+1/2,j+1/2}_{k,k',a_1a_2} 
&=\sum_{k}\sum_{k'}\frac{1}{2}\ovl{b}_k\ovl{b}_{k'}
\left[\left(\frac{1-\ovl{\eta}^2_{k}}{1+\ovl{\eta}^2_{k}}\right)\left(
\frac{1-\ovl{\eta}^2_{k'}}{1+\ovl{\eta}^2_{k'}}\right)
\psi^{j,j}_{k,k',a_1a_2}+\psi^{j,j}_{k,k',a_1a_2}\right]\\
&-\frac{1}{\rho}\left[\sum_{k}\Gamma_{k}\right]
\sum_{k}\ovl{b}_k\left(\frac{-\ovl{\eta}^2_{k}}{1+\ovl{\eta}^2_{k}}\right)
\left[\varphi_{k,a_1}^{j,j}+\varphi^{j,j}_{k,a_2}\right]
+\frac{1}{\rho^2}\left[\sum_{k}\Gamma_{k}\right]^2u^j\\
&+\frac{1}{\rho}\left[\sum_{k}\Gamma_{k}\right]\mathcal{B}^{j+1/2}_{a_1}(y_2)
+\frac{1}{\rho}\left[\sum_{k}\Gamma_{k}\right]\mathcal{B}^{j+1/2}_{a_2}(y_1)
+\frac{1}{\rho^2}\left[\sum_{k}\Gamma_{k}\right]^2 v^{j+1}.
\end{split}
\end{equation}
Plugging-in the value of $\psi^{j+1/2,j+1/2}_{k,k',a_1,a_2}$ in~\eqref{eq:cc-tr-npade} and 
collecting all the terms at $(j+1)$-th time-step on the left hand side, we arrive at the expression 
\begin{equation}
\begin{split}
&\left[\partial_{y_2}\partial_{y_1}+\alpha_2\varpi\partial_{y_1}
+\alpha_1\varpi\partial_{y_2}+\alpha_1\alpha_2\ovl{b}_0^2\right]v^{j+1}
+\alpha_1\alpha_2[2\varpi (\varpi-\ovl{b}_0)-(\varpi-\ovl{b}_0)^2]v^{j+1}\\
&=\left[\partial_{y_2}\partial_{y_1}+\alpha_2\varpi\partial_{y_1}
+\alpha_1\varpi\partial_{y_2}+\alpha_1\alpha_2\ovl{b}_0^2\right]v^{j+1}
+\alpha_1\alpha_2[\varpi^2-\ovl{b}_0^2]v^{j+1}\\
&=(\partial_{y_2}+\alpha_2\varpi)(\partial_{y_1}+\alpha_1\varpi){v}^{j+1}.
\end{split}
\end{equation}
With this observation, the corner conditions turns out to be 
\begin{equation}
\begin{split}
&\left.(\partial_{y_2}+\alpha_2\varpi)(\partial_{y_1}+\alpha_1\varpi)v^{j+1}\right|_{(y_r,y_t)}=+\alpha_1\alpha_2\mathcal{C}^{j+1/2}_{rt},\quad
\left.(\partial_{y_2}-\alpha_2\varpi)(\partial_{y_1}+\alpha_1\varpi)v^{j+1}\right|_{(y_r,y_b)}=-\alpha_1\alpha_2\mathcal{C}^{j+1/2}_{br},\\
&\left.(\partial_{y_2}+\alpha_2\varpi)(\partial_{y_1}-\alpha_1\varpi)v^{j+1}\right|_{(y_l,y_t)}=-\alpha_1\alpha_2\mathcal{C}^{j+1/2}_{tl},\quad
\left.(\partial_{y_2}-\alpha_2\varpi)(\partial_{y_1}-\alpha_1\varpi)v^{j+1}\right|_{(y_l,y_b)}=+\alpha_1\alpha_2\mathcal{C}^{j+1/2}_{lb},
\end{split}
\end{equation}
where
\begin{equation}
\begin{split}
\mathcal{C}^{j+1/2}_{a_1a_2}
&=\alpha_1\alpha_2\sum_{k}\sum_{k'}\frac{1}{2}\ovl{b}_k\ovl{b}_{k'}
\left[\left(\frac{1-\ovl{\eta}^2_{k}}{1+\ovl{\eta}^2_{k}}\right)\left(
\frac{1-\ovl{\eta}^2_{k'}}{1+\ovl{\eta}^2_{k'}}\right)
\psi^{j,j}_{k,k',a_1a_2}+\psi^{j,j}_{k,k',a_1a_2}\right]\\
&\quad-\alpha_1\alpha_2\frac{1}{\rho}\left[\sum_{k}\Gamma_{k}\right]
\sum_{k}\frac{1}{2}\ovl{b}_k\left(\frac{1-\ovl{\eta}^2_{k}}{1+\ovl{\eta}^2_{k}}-1\right)
\left[\varphi_{k,a_1}^{j,j}+\varphi^{j,j}_{k,a_2}\right]
+\alpha_1\alpha_2\frac{1}{\rho^2}\left[\sum_{k}\Gamma_{k}\right]^2u^j,
\end{split}
\end{equation}
for $a_1\in\{r,l\}$ and $a_2\in\{t,b\}$.
\subsection{Numerical solution of the IBVP}\label{sec:IBVP}
So far, we achieved the temporal discretization of the novel boundary conditions 
developed in this work together with the temporal discretization of the interior 
problem using one step time-stepping methods namely, BDF1 and TR. For the spatial 
discretization, we use a Legendre-Galerkin spectral method where we develop a new basis,
referred to as the boundary adapted basis~\cite{S1994}, in terms
of the Legendre polynomials to arrive at a banded linear system.
Recall that the Robin-type formulation of boundary maps discussed in this work 
turns out to be non-homogeneous in nature. The formulation of the linear system
using a Legendre-Galerkin method requires us to convert the TBCs from 
non-homogeneous to homogeneous form which is achieved by using a boundary 
lifting procedure. Further, we will first address the boundary lifting on the segments 
followed by the discussion on the rectangular computational domain $\Omega_i$.
Set $\field{I}=(-1,1)$ and let $L_n(y)$ denote the Legendre polynomial of degree 
$n$. Define the polynomial space 
\begin{equation}
\fs{P}_N = \Span\left\{L_p(y)|\;p=0,1,\ldots,N,\;y\in\ovl{\field{I}}\right\}.
\end{equation}

\subsubsection{Linear system: 1D}\label{sec:linear-sys-1D}
We observed that in order to implement the DtN maps for the interior field, we 
needed to solve a one-dimensional Schr\"odinger equation in terms of certain 
auxiliary functions with boundary segments as the underlying spatial domain 
together with the accompanying transparent boundary conditions.  
The recipe comprised two time-like variables, namely, $\tau_1$ and $\tau_2$, 
which played the role of the temporal dimension one at a time depending on 
the boundary segment involved. Therefore, it is worthwhile to consider a model 
IVP corresponding to the Schr\"odinger equation in 1D to represent each of the 
IVPs involved at the boundary segments.

Let $\field{S}_i=(x_-,x_+)$ denote the interior domain in this section. Let us 
state our model IVP in 1D with transparent boundary conditions as follows:
\begin{equation}\label{eq:1D-SE-CT}
\left\{\begin{aligned}
 &i\partial_tu+\partial^2_{x} u=0,\quad ({x},t)\in\field{S}_i\times\field{R}_+,\\
&u({x},0)=u_0({x})\in \fs{L}^2(\field{S}_i),\quad\supp\,u_0\subset\field{S}_i,\\
&\partial_{n}{u}+e^{-i\pi/4}\partial_{t}^{1/2}u=0,
\quad {x}\in\{x_-,x_+\},\quad t>0.
\end{aligned}\right.
\end{equation}
For the computational domain $\field{S}_i$, the reference domain is 
$\field{I}$. Let $y\in\field{I}$ such that
\begin{equation}
x = J y+\bar{x},\quad J = \frac{1}{2}(x_+-x_-),\quad
\bar{x}=\frac{1}{2}(x_-+x_+),\quad \beta=J^{-2}.
\end{equation}
Let $\Delta t$ denote the time-step. 
Following the earlier convention, $u^j(y)$ is taken to approximate 
$u(x,j\Delta t)$ for $j=0,1,2,\ldots,N_t-1$ for the sake of brevity 
of presentation. Also, we formulate the linear system in context of BDF1 
time-stepping for the interior problem as well as for discretizing the 
boundary conditions. The results can be easily extended to the case of TR. 
Setting $\rho=1/\Delta t$ and $\alpha=\sqrt{\rho/\beta}e^{-i\pi/4}$, 
then the BDF1 based discretization of~\eqref{eq:1D-SE-CT} becomes
\begin{equation}\label{eq:1D-SE}
\alpha^{-2}\partial^2_{y} u^{j+1}+ u^{j+1} = u^j,\quad
\left.(\partial_{y}-\kappa)u^{j+1}\right|_{y=-1}=+\alpha\mathcal{B}^{j+1}_-,\quad
\left.(\partial_{y}+\kappa)u^{j+1}\right|_{y=+1}=-\alpha\mathcal{B}^{j+1}_+,
\end{equation}
where $\kappa$ and history functions, $\mathcal{B}^{j+1}_{\pm}$, in the 
Robin-type boundary condition are specific to the method chosen, 
namely, CQ and NP, for discretizing the boundary conditions. 
In order to enforce the boundary conditions exactly, we introduce 
the space of boundary adapted basis functions given by
\begin{equation}
\fs{X}_N=\left\{u\in\fs{P}_N\left|\;
\begin{aligned}
&\left.(\partial_{y}-\kappa)u(y)\right|_{y=-1}=0,\\
&\left.(\partial_{y}+\kappa)u(y)\right|_{y=+1}=0
\end{aligned}\right.
\right\}.
\end{equation}
Let the index set $\{0,1,\ldots,N-2\}$ be denoted by $\field{J}$. 
Next, we introduce a function 
$\chi$ to convert the discrete TBCs to a homogeneous form by stating the original field as
$u^j=w^j+\chi^{j}(y)$ where $w^j\in\fs{X}_N$ so that
\begin{equation}\label{eq:1D-SE-homogenized}
\alpha^{-2}\partial^2_{y} w^{j+1}+ w^{j+1} = u^j-\chi^{j+1},\quad
\left.(\partial_{y}-\kappa)w^{j+1}\right|_{y=-1}=0,\quad
\left.(\partial_{y}+\kappa)w^{j+1}\right|_{y=+1}=0,
\end{equation}
and the field $\chi^j(y)$ is forced to satisfy the constraints below
\begin{equation}
\left.(\partial_{y}-\kappa)\chi^{j}(y)\right|_{y=-1}=+\alpha\mathcal{B}^{j}_-,\quad
\left.(\partial_{y}+\kappa)\chi^{j}(y)\right|_{y=+1}=-\alpha\mathcal{B}^{j}_+.
\end{equation}
The form of the equations suggest the use of an ansatz of the form 
$\chi^j(y)=c^j_-\chi_-+c^j_+\chi_+$ with $c^j_-=\alpha\mathcal{B}^j_-$ and 
$c^j_+=-\alpha\mathcal{B}^j_+$ so that the lifting functions satisfy
\begin{equation}
\left.(\partial_{y}-\kappa)
\begin{pmatrix}
\chi_-\\
\chi_+
\end{pmatrix}
\right|_{y=-1}=
\begin{pmatrix}
1\\
0
\end{pmatrix},\quad
\left.(\partial_{y}+\kappa)
\begin{pmatrix}
\chi_-\\
\chi_+
\end{pmatrix}
\right|_{y=+1}=
\begin{pmatrix}
0\\
1
\end{pmatrix}.
\end{equation}
From here, the degree of the lifting functions $\chi_{\pm}$ can be
inferred to be utmost $1$. Exploiting the fact 
that we can expand them in terms of Legendre polynomials, we have
\begin{equation}
\begin{pmatrix}
	\chi_-(y)\\
	\chi_+(y)
\end{pmatrix}=A_0
\begin{pmatrix}
	L_0(y)\\
	L_1(y)
\end{pmatrix},
\end{equation}
where $A_0$ is the unknown transformation matrix. Solving the linear system for $A_0$ yields
\begin{equation}
\left\{\begin{aligned}
\chi_-(y) &=-\frac{1}{2\kappa}L_0(y)+\frac{1}{2(\kappa+1)}L_1(y),\\
\chi_+(y) &=+\frac{1}{2\kappa}L_0(y)+\frac{1}{2(\kappa+1)}L_1(y).
\end{aligned}\right.
\end{equation}
The lifted field then takes the form
\begin{equation}
u^j(y)=w^j(y)-\frac{\alpha}{2\kappa}
\left[{\mathcal{B}}^j_{+}+{\mathcal{B}}^j_{-}\right]L_0(y)
-\frac{\alpha}{2(1+\kappa)}
\left[{\mathcal{B}}^j_{+}-{\mathcal{B}}^j_{-}\right]L_1(y).
\end{equation}
Next, we address the construction a basis $\{\phi_{p}|\,p\in\field{J}\}$ 
for $\fs{X}_N$ using the ansatz $\phi_{p}(y)=L_p(y)+a_pL_{p+1}(y)+b_pL_{p+2}(y)$.
Imposing the boundary conditions in the definition of $\fs{X}_N$, the 
sequences $(a_p)_{p\in\field{J}}$ and $(b_p)_{p\in\field{J}}$ work out to be
\begin{equation}
a_p=0,\quad
b_p=-\frac{\kappa+\frac{1}{2}p(p+1)}{\kappa+\frac{1}{2}(p+2)(p+3)},
\quad p\in\field{J}.
\end{equation}
Invoking the orthogonality property of the Legendre polynomials, 
system matrices can be computed easily. Denoting the stiffness matrix 
and the mass matrix for the boundary adapted Legendre basis by 
$S = (\melem{s}_{j,k})_{j,k\in\field{J}}$ and 
$M = (\melem{m}_{j,k})_{j,k\in\field{J}}$, respectively, where
\begin{equation}\label{eq:sys-mat-1d}
\melem{s}_{j,k}=-(\phi_j,\phi''_k)_{\field{I}}
=\begin{cases}
-2(2k+3)b_k,&j=k,\\
0,&\mbox{otherwise},
\end{cases}
\quad
\melem{m}_{j,k}=(\phi_j,\phi_k)_{\field{I}}
=\begin{cases}
\frac{2b_{k-2}}{2k+1},&j=k-2,\\
\frac{2}{2k+1}+\frac{2b_k^2}{2k+5},&j=k,\\
\frac{2b_{k}}{2k+5},&j=k+2,\\
0&\mbox{otherwise},
\end{cases}.
\end{equation}
Therefore, it is evident that the mass matrix is a pentadiagonal symmetric
matrix while the stiffness matrix is diagonal. 

Within the variational formulation, our goal is to find the 
approximate solution $u_N\in\fs{X}_N$ to~\eqref{eq:1D-SE} such that
\begin{equation}\label{eq:vform-SE1D}
\alpha^{-2}(\partial^2_{y} u_N^{j+1},\phi_p)_{\field{I}} 
+ (u_N^{j+1},\phi_p)_{\field{I}} = (u_N^j,\phi_p)_{\field{I}}, 
\quad \phi_p\in\fs{X}_N,\quad p\in\field{J},
\end{equation}
where $(u,\phi)_{\field{I}}=\int_{\field{I}}u\phi dy$ is the scalar product in 
$\fs{L}^2(\field{I})$. We approximate $u$ by polynomial interpolation over the 
Legendre-Gauss-Lobatto (LGL) nodes in the context of discrete Legendre 
transforms. 
The variational formulation~\eqref{eq:vform-SE1D} in terms of the 
unknown field $w^{j+1}$ (dropping the subscript `$N$' for the sake of brevity) 
can now be stated as
\begin{equation}\label{eq:vform-1d}
-\alpha^{-2}\left(\partial^2_{y}w^{j+1},\phi_{p}\right)_{\field{I}}
+(w^{j+1},\phi_{p})_{\field{I}}=(u^j,\phi_{p})_{\field{I}}
+\alpha\left(\mathcal{B}^{j+1}_{+}\chi_+-\mathcal{B}^{j+1}_{-}\chi_-,\phi_{p}\right)_{\field{I}}
\end{equation}
Let $\wtilde{u}_p$ denote the expansion coefficients of the field in 
Legendre basis and $\what{w}_p$ denote the expansion coefficients in the 
boundary adapted basis
\begin{equation}
w^{j+1}(y) = \sum_{p'\in\field{J}}\what{w}^{j+1}_{p'}\phi_{p'}(y),\quad
u^{j+1}(y)= \sum_{p'=0}^N\wtilde{u}^{j+1}_{p'}L_{p'}(y).
\end{equation}
where $\wtilde{\vv{u}}=(\wtilde{u}_0,\wtilde{u}_1,\ldots,\wtilde{u}_N)^{\tp}$
and $\what{\vv{w}}=(\what{w}_0,\what{w}_1,\ldots,\what{w}_{N-2})^{\tp}$. 
We can compute the inner product terms on LHS in~\eqref{eq:vform-1d} as
\begin{equation}
-\left(\partial^2_{y}w^{j+1},\phi_{p}\right)_{\field{I}}
=\sum_{p'}\melem{s}_{pp'}\what{w}^{j+1}_{p'}
=\left(S\what{\vv{w}}^{j+1}\right)_{p}
,\quad (w^{j+1},\phi_{p})_{\field{I}}
=\sum_{p'}\melem{m}_{pp'}\what{w}^{j+1}_{p'}
=\left(M\what{\vv{w}}^{j+1}\right)_{p}
\end{equation}
where $M$ and $S$ denote the mass and stiffness matrices, respectively, 
defined in~\eqref{eq:sys-mat-1d}. In matrix form, the expansion coefficients 
for the fields in Legendre basis satisfy
\begin{equation}
\wtilde{\vv{u}}^j=\wtilde{\vv{w}}^j-\frac{\alpha}{2\kappa}
\left[{\mathcal{B}}^j_{+}+{\mathcal{B}}^j_{-}\right]\vv{e}_0
-\frac{\alpha}{2(1+\kappa)}
\left[{\mathcal{B}}^j_{+}-{\mathcal{B}}^j_{-}\right]\vv{e}_1,
\end{equation}
where 
$\vv{e}_0 = (1,0,\ldots,0)\in\field{R}^{(N+1)\times1}$ and 
$\vv{e}_1 = (0,1,\ldots,0)\in\field{R}^{(N+1)\times1}$. We seek to 
further simplify the expression for the linear system by collecting the 
right hand side of~\eqref{eq:vform-1d} into one term. To this end, 
let $f^{j}(y)$ be such that
\begin{equation}
(u^j,\phi_{p})_{\field{I}}
+\alpha\left[
\mathcal{B}^{j+1}_{+}(\chi_+,\phi_{p})_{\field{I}}
-\mathcal{B}^{j+1}_{-}(\chi_-,\phi_{p})_{\field{I}}\right]
=(f^j,\phi_{p})_{\field{I}},
\end{equation}
then it follows that
\begin{equation}
f^j(y)
=u^j(y)+\frac{\alpha}{2\kappa}
\left(\mathcal{B}^{j+1}_{+}+\mathcal{B}^{j+1}_{-}\right)
+\frac{\alpha}{2(1+\kappa)}
\left(\mathcal{B}^{j+1}_{+}-\mathcal{B}^{j+1}_{-}\right).
\end{equation}
In matrix form, the expansion coefficients for the function $f^j(y)$ in the Legendre basis
satisfy
\begin{equation*}
\wtilde{\vv{f}}^j=\wtilde{\vv{u}}^j+\frac{\alpha}{2\kappa}
\left[{\mathcal{B}}^{j+1}_{+}+{\mathcal{B}}^{j+1}_{-}\right]\vv{e}_0
+\frac{\alpha}{2(1+\kappa)}
\left[{\mathcal{B}}^{j+1}_{+}-{\mathcal{B}}^{j+1}_{-}\right]\vv{e}_1.
\end{equation*}
Introducing the normalization factors 
$\gamma_k=\|L_k\|^2_2=2/(2k+1),\,k\in\field{N}_0$, we define the 
normalization matrix as $\Gamma=\diag(\gamma_0,\gamma_1,\ldots,\gamma_{N})$. 
The inner product $(f^j,\phi_{p})_{\field{I}}$ can be computed from the 
knowledge of the Legendre coefficient $\wtilde{\vv{f}}^j$ using the 
specific form of the boundary adapted basis. In~\ref{app:ip}, it is 
shown that the aforementioned operation can be 
achieved via a \emph{quadrature matrix} ($Q$) given by
\begin{equation}
Q=
\begin{pmatrix}
1   &  0  &  b_0 &      &      &      &\\
    &  1  &   0  &  b_1 &      &      &\\
    &     &   1  &   0  &  b_2 &      &\\
    &     &      &\ddots&\ddots&\ddots&\\
    &     &      &      &   1  & 0    &b_{N-2}\\
\end{pmatrix}\in\field{C}^{(N-1)\times (N+1)}\quad\text{such that}\quad
(f^j,\phi_{p})_{\field{I}}
=\left(Q\Gamma\wtilde{\vv{f}}^{j}\right)_p.
\end{equation}
The linear system for~\eqref{eq:vform-1d} then becomes
\begin{equation}\label{eq:linear-1d-system}
\alpha^{-2}S\what{\vv{w}}^{j+1}+M\what{\vv{w}}^{j+1}
=Q\Gamma\wtilde{\vv{f}}^{j}.
\end{equation}
In our implementation, we have found that, for the space depended variables, 
it is convenient to store them as samples over the LGL-nodes or as coefficients 
of their discrete Legendre transform. It is easy to work out that the transpose 
of the quadrature matrix which transforms Legendre coefficients to that of the 
boundary adapted basis (see~\ref{app:ip}), therefore, storage in terms of the
latter is not necessary. This linear system is 
then solved using LU-decomposition method. The matrix 
$\alpha^{-2}S+M$ is pentadiagonal. The LU decomposition of this
matrix inherits its bandedness provided pivoting is not required~\cite{Higham2002} 
and the resulting complexity of solving the linear system becomes $\bigO{N}$. A thorough
study of this issue is beyond the scope of this paper. 
Let us note that, in the realization of the TBCs, the linear system considered 
here plays a significant role. For a range of values of $\Delta t$, we note 
that MATLAB arrives at a banded $L$ and $U$ matrix consistent with lack of pivoting. 

\subsubsection{Linear system: 2D}
We formulate the linear system for solving the IBVP~\eqref{eq:2D-SE-CT} in the context
of the BDF1 time-stepping for the interior problem as well as for temporal 
discretization of the boundary conditions. The results can be easily extended to
the case of TR. For the 2D problem, we consider the tensor product space 
$\fs{P}_{N_1}\otimes\fs{P}_{N_2}$ of the Legendre polynomial spaces
$\fs{P}_{N_j},\,j=1,2$. In order to enforce 
the boundary conditions exactly, we introduce the space of boundary adapted basis given by
\begin{equation*}
\fs{X}_{N_1}\otimes\fs{X}_{N_2}=\left\{u\in\fs{P}_{N_1}\otimes\fs{P}_{N_2}\left|\;
\begin{aligned}
&\left.(\partial_{y_1}-\kappa_1)u(y_1,y_2)\right|_{y_1=-1}=0,\\
&\left.(\partial_{y_1}+\kappa_1)u(y_1,y_2)\right|_{y_1=+1}=0,\\
&\left.(\partial_{y_2}-\kappa_2)u(y_1,y_2)\right|_{y_2=-1}=0,\\
&\left.(\partial_{y_2}+\kappa_2)u(y_1,y_2)\right|_{y_2=+1}=0
\end{aligned}\right.
\right\},
\end{equation*}
where $\kappa_1$ and $\kappa_2$ are determined by the Robin-type formulation of 
TBCs specific to the cases, CQ and NP, for temporal discretization of the boundary conditions. 
\begin{rem}
The form of $\kappa$ follows from the Robin-type formulation of the discretized
version of the DtN-maps as discussed in Sec.~\ref{sec:dtn-maps}.
For CQ methods $\kappa_j=\alpha_j=\sqrt{\rho/\beta_j}\exp(-i\pi/4)$ while 
for the NP methods, we have $\kappa_j=\alpha_j\varpi$ where
$\varpi=\sqrt{\rho}R_M(\rho)$ as defined in~\eqref{eq:varpi-defn} for $j=1,2$.
\end{rem}

Let the index set $\{0,1,\ldots,N_j-2\}$ be denoted by $\field{J}_j$ for
$j=1,2$. We once again introduce a function $\chi$ in order to homogenize the 
Robin-type formulation of the discrete TBCs such that 
$u^{j+1}=w^{j+1}+\chi^{j+1}(y_1,y_2)$ 
where $w^j\in\fs{X}_{N_1}\otimes\fs{X}_{N_2}$ so that the interior problem 
in~\eqref{eq:2d-se-bdf1} reads as
\begin{equation}
-\left(\alpha_1^{-2}\partial^2_{y_1}+\alpha_2^{-2}\partial^2_{y_2}\right)w^{j+1}+w^{j+1}
=u^j-\chi^{j+1}+\left(\alpha_1^{-2}\partial^2_{y_1}+\alpha_2^{-2}\partial^2_{y_2}\right)\chi^{j+1},
\end{equation}
now with the homogeneous TBCs given by
\begin{equation}
\begin{split}
&\left.(\partial_{y_1}-\kappa_1)w^{j+1}(y_1,y_2)\right|_{y_1=-1}=0,\quad
\left.(\partial_{y_2}-\kappa_2)w^{j+1}(y_1,y_2)\right|_{y_2=-1}=0,\\
&\left.(\partial_{y_1}+\kappa_1)w^{j+1}(y_1,y_2)\right|_{y_1=+1}=0,\quad
\left.(\partial_{y_2}+\kappa_2)w^{j+1}(y_1,y_2)\right|_{y_2=+1}=0.
\end{split}
\end{equation}
The field $\chi^{j+1}(y_1,y_2)$ is forced to satisfy the following 
constraints on the segments 
\begin{equation}\label{eq:lift1}
\begin{split}
&\left.(\partial_{y_1}+\kappa_1)\chi^{j+1}(y_1,y_2)\right|_{y_1=y_{r}}=-\alpha_1\mathcal{B}^{j+1}_r(y_2),\quad
\left.(\partial_{y_2}+\kappa_2)\chi^{j+1}(y_1,y_2)\right|_{y_2=y_{t}}=-\alpha_2\mathcal{B}^{j+1}_t(y_1),\\
&\left.(\partial_{y_1}-\kappa_1)\chi^{j+1}(y_1,y_2)\right|_{y_1=y_{l}}=+\alpha_1\mathcal{B}^{j+1}_l(y_2),\quad
\left.(\partial_{y_2}-\kappa_2)\chi^{j+1}(y_1,y_2)\right|_{y_2=y_{b}}=+\alpha_2\mathcal{B}^{j+1}_b(y_1),
\end{split}
\end{equation}
together with the constraints on the corners below
\begin{equation}
\begin{split}
&\left.(\partial_{y_1}+\kappa_1)(\partial_{y_2}+\kappa_2)\chi^{j+1}(y_1,y_2)\right|_{(y_r,y_t)}=+\alpha_1\alpha_2\mathcal{C}^{j+1}_{rt},\quad
\left.(\partial_{y_1}-\kappa_1)(\partial_{y_2}-\kappa_2)\chi^{j+1}(y_1,y_2)\right|_{(y_l,y_b)}=+\alpha_1\alpha_2\mathcal{C}^{j+1}_{lb},\\
&\left.(\partial_{y_1}+\kappa_1)(\partial_{y_2}-\kappa_2)\chi^{j+1}(y_1,y_2)\right|_{(y_r,y_b)}=-\alpha_1\alpha_2\mathcal{C}^{j+1}_{rb},\quad
\left.(\partial_{y_1}-\kappa_1)(\partial_{y_2}+\kappa_2)\chi^{j+1}(y_1,y_2)\right|_{(y_l,y_t)}=-\alpha_1\alpha_2\mathcal{C}^{j+1}_{lt}.
\end{split}
\end{equation}
These boundary condition are taken from Sec.~\ref{sec:dtn-maps} where precise definition
of the quantities are provided. We can express these operators compactly by introducing the notation
$\partial_1^{\pm}\equiv (\partial_{y_1}\pm\alpha_1),\; \partial_2^{\pm}\equiv (\partial_{y_2}\pm\alpha_2)$.
Therefore,
\begin{equation}
\begin{split}
&\left.\partial^+_1\partial^+_2\chi^{j+1}(y_1,y_2)\right|_{(y_r,y_t)}
=+\alpha_1\alpha_2\mathcal{C}^{j+1}_{rt},\quad
\left.\partial_1^-\partial_2^-\chi^{j+1}(y_1,y_2)\right|_{(y_l,y_b)}
=+\alpha_1\alpha_2\mathcal{C}^{j+1}_{lb},\\
&\left.\partial_1^+\partial_2^-\chi^{j+1}(y_1,y_2)\right|_{(y_r,y_b)}
=-\alpha_1\alpha_2\mathcal{C}^{j+1}_{rb},\quad
\left.\partial_1^-\partial_2^+\chi^{j+1}(y_1,y_2)\right|_{(y_l,y_t)}
=-\alpha_1\alpha_2\mathcal{C}^{j+1}_{lt}.
\end{split}
\end{equation}
The form of the equations suggest an ansatz of the form
\begin{equation}
\chi^{j+1}=\eta^{j+1} 
+\alpha_1\alpha_2\mathcal{C}^{j+1}_{rt}\chi_{rt}
+\alpha_1\alpha_2\mathcal{C}^{j+1}_{lb}\chi_{lb}
-\alpha_1\alpha_2\mathcal{C}^{j+1}_{rb}\chi_{rb}
-\alpha_1\alpha_2\mathcal{C}^{j+1}_{lt}\chi_{lt},
\end{equation}
where $\eta^{j+1}$ satisfy 
\begin{equation}
\begin{split}
&\left. \partial^+_1\partial^+_2\eta^{j+1}(y_1,y_2)\right|_{(y_r,y_t)}=0,\quad
\left. \partial_1^-\partial_2^-\eta^{j+1}(y_1,y_2)\right|_{(y_l,y_b)}=0,\\
&\left. \partial_1^+\partial_2^-\eta^{j+1}(y_1,y_2)\right|_{(y_r,y_b)}=0,\quad
\left. \partial_1^-\partial_2^+\eta^{j+1}(y_1,y_2)\right|_{(y_l,y_t)}=0.
\end{split}
\end{equation}
Here, we suppress the temporal index $j$ for the sake of brevity of presentation.
Following along the similar lines as that of boundary lifting in 1D, we can set up the constraint 
equations for the lift functions as following:
\begin{equation}
\left.\partial^+_1\partial^+_2
\begin{pmatrix}
\chi_{rt}\\
\chi_{lb}\\
\chi_{rb}\\
\chi_{lt}
\end{pmatrix}
\right|_{(y_r,y_t)}=
\begin{pmatrix}
1\\
0\\
0\\
0
\end{pmatrix},\quad
\left.\partial^-_1\partial^-_2
\begin{pmatrix}
\chi_{rt}\\
\chi_{lb}\\
\chi_{rb}\\
\chi_{lt}
\end{pmatrix}
\right|_{(y_l,y_b)}=
\begin{pmatrix}
0\\
1\\
0\\
0
\end{pmatrix},\quad
\left.\partial^+_1\partial^-_2
\begin{pmatrix}
\chi_{rt}\\
\chi_{lb}\\
\chi_{rb}\\
\chi_{lt}
\end{pmatrix}
\right|_{(y_r,y_b)}=
\begin{pmatrix}
0\\
0\\
1\\
0
\end{pmatrix},\quad
\left.\partial^-_1\partial^+_2
\begin{pmatrix}
\chi_{rt}\\
\chi_{lb}\\
\chi_{rb}\\
\chi_{lt}
\end{pmatrix}
\right|_{(y_l,y_t)}=
\begin{pmatrix}
0\\
0\\
0\\
1
\end{pmatrix}.
\end{equation}
From here, the degree of the lifting functions can be inferred to be utmost 
$1$ in terms of each of the independent variables. Exploiting 
the fact that we can expand them in terms of Legendre polynomials, we have
\begin{equation}
\begin{pmatrix}
\chi_{rt}\\
\chi_{lb}\\
\chi_{rb}\\
\chi_{lt}\\
\end{pmatrix}
=C
\begin{pmatrix}
L_0(y_1)L_0(y_2)\\
L_0(y_1)L_1(y_2)\\
L_1(y_1)L_0(y_2)\\
L_1(y_1)L_1(y_2)
\end{pmatrix},
\end{equation}
where $C$ is the unknown transformation matrix. Solving the linear system
for $C$ yields
\begin{equation}
\begin{split}
&\chi_{rt}(y_1,y_2)=\chi_r(y_1)\chi_t(y_2),\quad
\chi_{lb}(y_1,y_2)=\chi_l(y_1)\chi_b(y_2),\\
&\chi_{rb}(y_1,y_2)=\chi_r(y_1)\chi_b(y_2),\quad
\chi_{lt}(y_1,y_2)=\chi_l(y_1)\chi_t(y_2),
\end{split}
\end{equation}
where
\begin{equation}
\left\{\begin{aligned}
\chi_r(y_1) &=+\frac{1}{2\kappa_1}L_0(y_1)+\frac{1}{2(\kappa_1+1)}L_1(y_1),\\
\chi_l(y_1) &=-\frac{1}{2\kappa_1}L_0(y_1)+\frac{1}{2(\kappa_1+1)}L_1(y_1),\\
\end{aligned}\right.\quad\text{and}\quad
\left\{\begin{aligned}
\chi_t(y_2) &=+\frac{1}{2\kappa_2}L_0(y_2)+\frac{1}{2(\kappa_2+1)}L_1(y_2),\\
\chi_b(y_2) &=-\frac{1}{2\kappa_2}L_0(y_2)+\frac{1}{2(\kappa_2+1)}L_1(y_2).
\end{aligned}\right.
\end{equation}
Plugging-in the form of $\chi$ in~\eqref{eq:lift1} and favouring $\eta$, 
we get 
\begin{equation}
\begin{split}
&\left.\partial^+_1\eta(y_1,y_2)\right|_{y_1=y_r}
=-\alpha_1\left[\mathcal{B}_r(y_2)-\alpha_2\mathcal{C}_{rb}\chi_b(y_2)
+\alpha_2\mathcal{C}_{rt}\chi_t(y_2)\right],\\
&\left.\partial^-_1\eta(y_1,y_2)\right|_{y_1=y_l}
=+\alpha_1\left[\mathcal{B}_l(y_2)-\alpha_2\mathcal{C}_{lb}\chi_b(y_2)
+\alpha_2\mathcal{C}_{lt}\chi_t(y_2)\right],\\
&\left.\partial^+_2\eta(y_1,y_2)\right|_{y_2=y_t}
=-\alpha_2\left[\mathcal{B}_t(y_1)-\alpha_1\mathcal{C}_{lt}\chi_l(y_1)
+\alpha_1\mathcal{C}_{rt}\chi_r(y_1)\right],\\
&\left.\partial^-_2\eta(y_1,y_2)\right|_{y_2=y_b}
=+\alpha_2\left[\mathcal{B}_b(y_1)-\alpha_1\mathcal{C}_{lb}\chi_l(y_1)
+\alpha_1\mathcal{C}_{rb}\chi_r(y_1)\right].
\end{split}
\end{equation}
We can rewrite the set of equations above by identifying the RHS as:
\begin{equation}
\begin{split}
&\left.\partial^+_1\eta(y_1,y_2)\right|_{y_1=y_r}=-\alpha_1\breve{\mathcal{B}}_r(y_2),\quad
\left.\partial^+_2\eta(y_1,y_2)\right|_{y_2=y_t}=-\alpha_2\breve{\mathcal{B}}_t(y_1),\\
&\left.\partial^-_1\eta(y_1,y_2)\right|_{y_1=y_l}=+\alpha_1\breve{\mathcal{B}}_l(y_2),\quad
\left.\partial^-_2\eta(y_1,y_2)\right|_{y_2=y_b}=+\alpha_2\breve{\mathcal{B}}_b(y_1).
\end{split}
\end{equation}
Restoring the superscripts, the system above can be solved in a similar manner 
as that of 1D case yielding
\begin{equation}
\eta^{j+1}=
-\alpha_1\chi_r(y_1)\breve{\mathcal{B}}^{j+1}_r(y_2)
+\alpha_1\chi_l(y_1)\breve{\mathcal{B}}^{j+1}_l(y_2)
-\alpha_2\breve{\mathcal{B}}^{j+1}_t(y_1)\chi_t(y_2)
+\alpha_2\breve{\mathcal{B}}^{j+1}_b(y_1)\chi_b(y_2).
\end{equation}
The lifted field $u^{j+1}=w^{j+1}+\chi^{j+1}(y_1,y_2)$ now expands to
\begin{equation}\label{eq:boundary-corner-lifting}
\begin{split}
u^{j+1} &= w^{j+1}
-\alpha_1\chi_r(y_1){\mathcal{B}}^{j+1}_r(y_2)
+\alpha_1\chi_l(y_1){\mathcal{B}}^{j+1}_l(y_2)
-\alpha_2{\mathcal{B}}^{j+1}_t(y_1)\chi_t(y_2)
+\alpha_2{\mathcal{B}}^{j+1}_b(y_1)\chi_b(y_2)\\
&\quad
-\alpha_1\alpha_2\mathcal{C}^{j+1}_{rt}\chi_r(y_1)\chi_t(y_2)
-\alpha_1\alpha_2\mathcal{C}^{j+1}_{lb}\chi_l(y_1)\chi_b(y_2)
+\alpha_1\alpha_2\mathcal{C}^{j+1}_{rb}\chi_r(y_1)\chi_b(y_2)
+\alpha_1\alpha_2\mathcal{C}^{j+1}_{lt}\chi_l(y_1)\chi_t(y_2).
\end{split}
\end{equation}
Next, we address the construction of the basis for the boundary adapted space 
$\fs{X}_{N_1}\otimes\fs{X}_{N_2}$ by using the ansatz
\begin{equation}
\theta_{p_1,p_2}(y_1,y_2)=\phi^{(1)}_{p_1}(y_1)\phi^{(2)}_{p_2}(y_2),
\end{equation}
where
\begin{equation}
\phi^{(j)}_{p_j}(y_j)=L_{p_j}(y_j)+a^{(j)}_{p_j}L_{p_j+1}(y_j)+b^{(j)}_{p_j}L_{p_j+2}(y_j),
\quad p_1\in\field{J}_1,\;p_2\in\field{J}_2.
\end{equation}
Imposing the boundary conditions in the definition of $\fs{X}_{N_1}\otimes\fs{X}_{N_2}$, the 
sequences $\{a^{(j)}_{p_j}\}$ and $\{b^{(j)}_{p_j}\}$ work out to be
\begin{equation}
a^{(j)}_{p_j}=0,\quad
b^{(j)}_{p_j}=-\frac{\kappa_j+\frac{1}{2}p_j(p_j+1)}{\kappa_j+\frac{1}{2}(p_j+2)(p_j+3)},
\quad p_j\in\field{J}_j,\;j=1,2.
\end{equation}
The stiffness matrix and the mass matrix for the boundary adapted Legendre basis are
denoted by $S_p = (\melem{s}^{(p)}_{kj})_{k,j\in\field{J}_p}$ and
$M_p = (\melem{m}^{(p)}_{kj})_{k,j\in\field{J}_p},\;p=1,2$, respectively, where
\begin{equation}\label{eq:sys-mat-2d}
\melem{s}^{(p)}_{jk}=-\left(\phi^{(p)}_j,\partial^2_{y_l}\phi^{(p)}_k\right)_{\field{I}}
=\begin{cases}
-2(2k+3)b^{(p)}_k,&j=k,\\
0,&\mbox{otherwise},
\end{cases}\quad
\melem{m}^{(p)}_{jk}=\left(\phi^{(p)}_j,\phi^{(p)}_k\right)_{\field{I}}
=\begin{cases}
\frac{2b^{(p)}_{k-2}}{2k+1},&j=k-2,\\
\frac{2}{2k+1}+\frac{2(b^{(p)}_k)^2}{2k+5},&j=k,\\
\frac{2b^{(p)}_{k}}{2k+5},&j=k+2,\\
0&\mbox{otherwise}.
\end{cases}
\end{equation}
It is interesting to note that mass matrices are pentadiagonal and 
stiffness matrices are simply diagonal. 

Within the variational formulation, our goal is to find the approximate solution 
$u_{N_1,N_2}\in\fs{X}_{N_1}\otimes\fs{X}_{N_2}$ for the interior problem in~\eqref{eq:2d-se-bdf1}:
\begin{equation}
\begin{split}
&-\alpha_1^{-2}\left(\partial^2_{y_1}u_{N_1,N_2}^{j+1},\theta_{p_1,p_2}\right)_{\Omega_i^{\text{ref}}}
 -\alpha^{-2}_2\left(\partial^2_{y_2}u_{N_1,N_2}^{j+1},\theta_{p_1,p_2}\right)_{\Omega_i^{\text{ref}}}
 +\left(u_{N_1,N_2}^{j+1},\theta_{p_1,p_2}\right)_{\Omega_i^{\text{ref}}}
=\left(u_{N_1,N_2}^j,\theta_{p_1,p_2}\right)_{\Omega_i^{\text{ref}}},
\quad \theta_{p_1,p_2}\in\fs{X}_{N_1}\otimes\fs{X}_{N_2},
\end{split}
\end{equation}
where $(u,\theta)_{\Omega_i^{\text{ref}}}=\int_{\Omega_i^{\text{ref}}}u\theta d^2\vv{y}$ is the 
scalar product in $\fs{L}^2(\Omega_i^{\text{ref}})$.
The variational formulation defined above can be expressed in terms of
unknown field $w^{j+1}$ (dropping the subscripts `$N_1$' and `$N_2$' for the sake of brevity) as
\begin{equation}\label{eq:vform-2d}
\begin{split}
-\alpha_1^{-2}\left(\partial^2_{y_1}w^{j+1},\theta_{p_1,p_2}\right)_{\Omega_i^{\text{ref}}}
 -\alpha^{-2}_2\left(\partial^2_{y_2}w^{j+1},\theta_{p_1,p_2}\right)_{\Omega_i^{\text{ref}}}
& +\left(w^{j+1},\theta_{p_1,p_2}\right)_{\Omega_i^{\text{ref}}}\\
 =\left(u^j-\chi^{j+1},\theta_{p_1,p_2}\right)_{\Omega_i^{\text{ref}}}
&+\alpha_1^{-2}\left(\partial^2_{y_1}\chi^{j+1},\theta_{p_1,p_2}\right)_{\Omega_i^{\text{ref}}}
+\alpha^{-2}_2\left(\partial^2_{y_2}\chi^{j+1},\theta_{p_1,p_2}\right)_{\Omega_i^{\text{ref}}}.
\end{split}
\end{equation}
Next, we simplify the second order derivative terms on the right hand side as:
\begin{equation}
\begin{split}
&\partial^2_{y_1}\chi^{j+1}=
-\alpha_2\partial^2_{y_1}{\mathcal{B}}^{j+1}_t(y_1)\chi_t(y_2)
+\alpha_2\partial^2_{y_1}{\mathcal{B}}^{j+1}_b(y_1)\chi_b(y_2),\\
&\partial^2_{y_2}\chi^{j+1}=
-\alpha_1\partial^2_{y_2}{\mathcal{B}}^{j+1}_r(y_2)\chi_r(y_1)
+\alpha_1\partial^2_{y_2}{\mathcal{B}}^{j+1}_l(y_2)\chi_l(y_1).
\end{split}
\end{equation}
Let $\wtilde{u}_{p,p'}$ denote the expansion coefficients of the field in Legendre basis 
and $\what{w}_{p,p'}$ denote the expansion coefficients in the boundary adapted basis as
\begin{equation}
u^{j+1}(y_1,y_2)=\sum_{p'=0}^{N_1}\sum_{q'=0}^{N_2}\wtilde{u}^{j+1}_{p',q'}L_{p'}(y_1)L_{q'}(y_2),
\quad\text{and}\quad
w^{j+1}(y_1,y_2)=\sum_{p'\in\field{J}_1}\sum_{q'\in\field{J}_2}\what{w}^{j+1}_{p',q'}
\phi^{(1)}_{p'}(y_1)\phi^{(2)}_{q'}(y_2),
\end{equation}
respectively. Let us introduce the matrices 
$\what{W}=\left(\what{w}_{p',q'}\right)_{p'\in\field{J}_1,\;q'\in\field{J}_2}$
and $\wtilde{U}=\left(\wtilde{u}_{p',q'}\right)_{0\leq p'\leq N_1,\;0\leq q'\leq N_2}$ 
for convenience. In matrix form, the expansion coefficients for the 
fields introduced in~\eqref{eq:boundary-corner-lifting} in the Legendre basis
assume the form
\begin{equation}
\begin{split}
\wtilde{U}^{j+1}
= \wtilde{W}^{j+1}&
-\frac{\alpha_1}{2\kappa_1}\vv{e}^{(1)}_0\otimes
\left(\wtilde{\vs{\mathcal{B}}}^{j+1}_r+\wtilde{\vs{\mathcal{B}}}^{j+1}_l\right)^{\tp}
-\frac{\alpha_1}{2(\kappa_1+1)}\vv{e}^{(1)}_1\otimes
\left(\wtilde{\vs{\mathcal{B}}}^{j+1}_r-\wtilde{\vs{\mathcal{B}}}^{j+1}_l\right)^{\tp}\\
&-\frac{\alpha_2}{2\kappa_2}\left(\wtilde{\vs{\mathcal{B}}}^{j+1}_t+\wtilde{\vs{\mathcal{B}}}^{j+1}_b\right)
\otimes\left(\vv{e}_0^{(2)}\right)^{\tp}
-\frac{\alpha_2}{2(\kappa_2+1)}
\left(\wtilde{\vs{\mathcal{B}}}^{j+1}_t-\wtilde{\vs{\mathcal{B}}}^{j+1}_b\right)
\otimes\left(\vv{e}_1^{(2)}\right)^{\tp}\\
&-\frac{\alpha_1\alpha_2}{4\kappa_1\kappa_2}\left(\mathcal{C}^{j+1}_{rt}+\mathcal{C}^{j+1}_{lb}
+\mathcal{C}^{j+1}_{rb}+\mathcal{C}^{j+1}_{lt}\right)
\vv{e}^{(1)}_0\otimes\left(\vv{e}^{(2)}_0\right)^{\tp}\\
&-\frac{\alpha_1\alpha_2}{4\kappa_1(\kappa_2+1)}
\left(\mathcal{C}^{j+1}_{rt}-\mathcal{C}^{j+1}_{lb}
-\mathcal{C}^{j+1}_{rb}+\mathcal{C}^{j+1}_{lt}\right)
\vv{e}^{(1)}_0\otimes\left(\vv{e}^{(2)}_1\right)^{\tp}\\
&-\frac{\alpha_1\alpha_2}{4\kappa_2(\kappa_1+1)}
\left(\mathcal{C}^{j+1}_{rt}-\mathcal{C}^{j+1}_{lb}
+\mathcal{C}^{j+1}_{rb}-\mathcal{C}^{j+1}_{lt}\right)
\vv{e}^{(1)}_1\otimes\left(\vv{e}^{(2)}_0\right)^{\tp}\\
&-\frac{\alpha_1\alpha_2}{4(\kappa_1+1)(\kappa_2+1)}
\left(\mathcal{C}^{j+1}_{rt}+\mathcal{C}^{j+1}_{lb}
-\mathcal{C}^{j+1}_{rb}-\mathcal{C}^{j+1}_{lt}\right)
\vv{e}^{(1)}_1\otimes\left(\vv{e}^{(2)}_1\right)^{\tp},
\end{split}
\end{equation}
where the hat and tilde convention is extended to other quantities as well.
We seek to further simplify the expression for the linear system by collecting the 
right hand side of~\eqref{eq:vform-2d} into one term. To this end, let $f^j(y_1,y_2)$ 
be such that
\begin{equation}
\left(f^j,\theta_{p_1,p_2}\right)
=\left(u^j-\chi^{j+1},\theta_{p_1,p_2}\right)
+\alpha_1^{-2}\left(\partial^2_{y_1}\chi^{j+1},\theta_{p_1,p_2}\right)
+\alpha^{-2}_2\left(\partial^2_{y_2}\chi^{j+1},\theta_{p_1,p_2}\right).
\end{equation}
In matrix form, the expansion coefficients for the function $f^j(y_1,y_2)$ in the Legendre basis
satisfy
\begin{equation}
\begin{split}
\wtilde{F}^j
= \wtilde{U}^{j}
&
+\frac{\alpha_1}{2\kappa_1}\vv{e}^{(1)}_0\otimes
\left(\wtilde{\vs{\mathcal{B}}}^{j+1}_r+\wtilde{\vs{\mathcal{B}}}^{j+1}_l\right)^{\tp}\mathcal{D}^{\tp}_2
+\frac{\alpha_1}{2(\kappa_1+1)}\vv{e}^{(1)}_1\otimes
\left(\wtilde{\vs{\mathcal{B}}}^{j+1}_r-\wtilde{\vs{\mathcal{B}}}^{j+1}_l\right)^{\tp}\mathcal{D}^{\tp}_2\\
&
+\frac{\alpha_2}{2\kappa_2}\mathcal{D}_1\left(\wtilde{\vs{\mathcal{B}}}^{j+1}_t+\wtilde{\vs{\mathcal{B}}}^{j+1}_b\right)
\otimes\left(\vv{e}_0^{(2)}\right)^{\tp}
+\frac{\alpha_2}{2(\kappa_2+1)}
\mathcal{D}_1\left(\wtilde{\vs{\mathcal{B}}}^{j+1}_t-\wtilde{\vs{\mathcal{B}}}^{j+1}_b\right)
\otimes\left(\vv{e}_1^{(2)}\right)^{\tp}+\wtilde{G}^{j+1},
\end{split}
\end{equation}
where $\mathcal{D}_k:1-\alpha_k^{-2}\partial^2_{y_k},\;k=1,2,$ are symbolic
representations of the second order differential 
operators\footnote{We would like to mention that the second order derivative of the history 
functions will be handled by first expanding them in the Legendre basis and 
then using the backward recursion in the frequency space to obtain the expansion 
coefficients for $\partial^2_{y_k}{\mathcal{B}}^{j+1}(y_k),\;k=1,2$. For more details the
reader is referred to Shen~\cite[Chap.~3]{Shen2011}.} and 
\begin{equation}
\wtilde{G}^{j+1}=(\wtilde{g}^{j+1}_{p_1,p_2})
=\begin{cases}
\frac{\alpha_1\alpha_2}{4\kappa_1(\kappa_2+1)} \left(\mathcal{C}^{j+1}_{rt}+\mathcal{C}^{j+1}_{lb}
+\mathcal{C}^{j+1}_{rb}+\mathcal{C}^{j+1}_{lt}\right), & p_1=0,p_2=0,\\
\frac{\alpha_1\alpha_2}{4\kappa_1(\kappa_2+1)}
\left(\mathcal{C}^{j+1}_{rt}-\mathcal{C}^{j+1}_{lb}
-\mathcal{C}^{j+1}_{rb}+\mathcal{C}^{j+1}_{lt}\right),& p_1=0,p_2=1,\\
\frac{\alpha_1\alpha_2}{4\kappa_2(\kappa_1+1)}
\left(\mathcal{C}^{j+1}_{rt}-\mathcal{C}^{j+1}_{lb}
+\mathcal{C}^{j+1}_{rb}-\mathcal{C}^{j+1}_{lt}\right),& p_1=1,p_2=0,\\
\frac{\alpha_1\alpha_2}{4(\kappa_1+1)(\kappa_2+1)}
\left(\mathcal{C}^{j+1}_{rt}+\mathcal{C}^{j+1}_{lb}
-\mathcal{C}^{j+1}_{rb}-\mathcal{C}^{j+1}_{lt}\right),& p_1=1,p_2=1,\\
0, & p_1=2,\ldots,N_1,\;p_2=2,\ldots,N_2, 
\end{cases}
\end{equation}
The variational formulation in~\eqref{eq:vform-2d} now takes the form
\begin{equation}
\begin{split}
&-\alpha_1^{-2}\left(\partial^2_{y_1}w^{j+1},\theta_{p_1,p_2}\right)_{\Omega_i^{\text{ref}}}
-\alpha^{-2}_2\left(\partial^2_{y_2}w^{j+1},\theta_{p_1,p_2}\right)_{\Omega_i^{\text{ref}}}
+\left(w^{j+1},\theta_{p_1,p_2}\right)_{\Omega_i^{\text{ref}}}
=\left(f^j,\theta_{p_1,p_2}\right)_{\Omega_i^{\text{ref}}} .
\end{split}
\end{equation}
The inner product $(f^j,\theta_{p_1,p_2})_{\Omega_i^{\text{ref}}}$ can be computed 
from the knowledge of the Legendre coefficients $\wtilde{F}^j$ using the specific 
form of the boundary adapted basis. In~\ref{app:ip}, it is 
shown that the aforementioned operation can be 
achieved via \emph{quadrature matrices} ($Q_j$) given by
\begin{equation}
Q_j=B_j^{\tp}=
\begin{pmatrix}
1   &  0  &  b^{(j)}_0 &      &      &      &\\
    &  1  &   0  &  b^{(j)}_1 &      &      &\\
    &     &   1  &   0  &  b^{(j)}_2 &      &\\
    &     &      &\ddots&\ddots&\ddots&\\
    &     &      &      &   1  & 0    &b^{(j)}_{N-2}\\
\end{pmatrix}\in\field{C}^{(N-1)\times (N+1)},\quad j=1,2.
\end{equation}
The linear system thus becomes
\begin{equation}\label{eq:linear-2d-system}
\begin{split}
&\alpha_1^{-2}S_1^{\tp}\what{W}^{j+1}M_2+\alpha^{-2}_2M_1^{\tp}\what{W}^{j+1}S_2
+M_1^{\tp}\what{W}^{j+1}M_2=Q_1\Gamma\wtilde{F}^{j}\Gamma Q_2^{\tp}\equiv \what{F},\\
&\left(\alpha_1^{-2}M_2^{\tp}\otimes S_1^{\tp}+\alpha_2^{-2}S_2^{\tp}\otimes M_1^{\tp}
+M_2^{\tp}\otimes M_1^{\tp}\right)\what{\vv{w}}^{j+1}=\what{\vv{f}}^j, 
\end{split}
\end{equation}
where $\otimes$ denotes the tensor product and $\what{\vv{w}},\;\what{\vv{f}}$
denote the column vectors obtained by stacking columns of matrices 
$\what{W},\;\what{F}$, respectively, below one another. The mass and stiffness 
matrices are defined in~\eqref{eq:sys-mat-2d} and the linear system is then 
solved using LU-decomposition method. In order to understand the sparsity pattern 
of the system matrix, we proceed as follows: Taking into account the symmetric 
nature of the mass matrices (noting stiffness matrices are already diagonal), 
we have
\begin{equation}
\left(\alpha_1^{-2}M_2^{\tp}\otimes S_1^{\tp}+\alpha_2^{-2}S_2^{\tp}\otimes M_1^{\tp}
+M_2^{\tp}\otimes M_1^{\tp}\right)= 
\left(\alpha_2^{-2}S_2+M_2\right)\otimes\left(\alpha_1^{-2}S_1+M_1\right)
-\alpha_2^{-2}\alpha_1^{-2}\left(S_2\otimes S_1\right).
\end{equation}
From the property of Kronecker product, it follows that the system matrix is 
symmetric block pentadiagonal because the individual matrices in the Kronecker 
product are symmetric pentadiagonal in nature. The upper and lower bandwidth
works out to be $2N_1$. The cost of LU decomposition can be worked out to be
$\bigO{N^2_1N_{\text{dim}}}$ where $N_{\text{dim}}=(N_1-1)(N_2-1)$. The lower and upper
triangular matrices inherit the bandedness of the system matrix 
provided pivoting is not required~\cite{Higham2002}. Under this assumption, 
the resulting complexity of solving the linear system (excluding the cost of LU
decomposition) becomes $\bigO{N_1N_{\text{dim}}}$. A thorough study of this 
issue is beyond the scope of this paper. For a range of values of $\Delta t$, we note 
that MATLAB arrives at banded lower and upper triangular matrices which is 
consistent with lack of pivoting.

\subsection{Algorithmic summary}\label{sec:algo}
In this section, we present a brief algorithmic summary for the CQ and NP methods
in context of BDF1 based temporal discretization of the boundary maps and also the
interior problem. For the TR case, the algorithmic steps are largely similar.

\subsubsection{CQ--BDF1}
Let us recall that, on the boundary segments $\Gamma_{a_1}$ and $\Gamma_{a_2}$, 
the auxiliary function is denoted by 
$\varphi_{a_1}(y_2,\tau_1,\tau_2)$ and $\varphi_{a_2}(y_1,\tau_1,\tau_2)$, 
respectively, while the auxiliary function at the corners is denoted by 
$\varphi_{a_1a_2}(\tau_1,\tau_2)$ where $a_1\in\{l,r\},\;a_2\in\{b,t\}$. Let 
$N$ be the number of spatial nodes along $y_1$ as well as $y_2$, i.e.,
$N_1=N_2=N$. Let the number of time-steps be denoted by $N_t$. The temporal
variables $t$, $\tau_1$ and $\tau_2$ are discretized with the same step-size
$\Delta t$ so that the maximum time $T=N_t\Delta t$. The storage of the 
auxiliary field on the boundary segments require four 
3D-arrays of size $N\times N_t\times N_t$ and the storage of auxiliary field 
at corner points require four 2D-arrays of size $N_t\times N_t$ which 
corresponds to $\varphi_{a_1a_2}(\tau_1,\tau_2)$. The algorithmic steps for the
solving the IBVP are enumerated below.
\begin{enumerate}[%
,leftmargin=*
,label={\bfseries Step \arabic*:}
]
\item Advance the following IVPs corresponding to auxiliary functions 
by one-step, say, $j\rightarrow (j+1)$ using the previously computed values 
of the auxiliary field on the segments (as described in
Fig.~\ref{fig:IVP-auxi}):
\begin{equation*}
\begin{split}
&-\alpha^{-2}_1\partial_{y_1}^2\varphi^{j+1,q}_{a_2}+\varphi^{j+1,q}_{a_2}
=\varphi^{j,q}_{a_2},\quad q=0,1,\ldots,j,\\
&-\alpha^{-2}_2\partial_{y_2}^2\varphi^{m,j+1}_{a_1}+\varphi^{m,j+1}_{a_1}
=\varphi^{m,j}_{a_1},\quad m=0,1,\ldots,j,
\end{split}
\end{equation*}
using the TBCs at the corner points given by 
\begin{equation*}
(\partial_{n_1}+\alpha_1)\varphi^{j+1,q}_{a_2}(y_{a_1})
=-\alpha_1\mathcal{B}^{j+1,q}_{a_2,a_1},\quad
(\partial_{n_2}+\alpha_2)\varphi^{m,j+1}_{a_1}(y_{a_2})
=-\alpha_2\mathcal{B}^{m,j+1}_{a_1,a_2},
\end{equation*}
where the history functions are computed as
\begin{equation*}
\mathcal{B}^{j+1,q}_{a_2,a_1} =\sum_{k=1}^{j+1}\omega_{k}\varphi^{j+1-k,q}_{a_1a_2},\quad
\mathcal{B}^{m,j+1}_{a_1,a_2} =\sum_{k=1}^{j+1}\omega_{k}\varphi^{m,j+1-k}_{a_1a_2}.
\end{equation*}
Let the cost of solving the linear 1D system of the 
kind~\eqref{eq:linear-1d-system} be denoted by $\Theta_1(N_{\text{dim}})$ where the 
dimension is $N_{\text{dim}}=N-1$. The cost of the steps outlined here works out to be 
$4(j+1)\Theta_1(N-1)+4(j+1)^2$ where the second part is the cost of computing
the history functions. The memory requirement is $4(j+1)^2N+4(j+1)^2$ size array.
Therefore, the complexity scales with $j$. Moreover, 
for large number of time-steps, i.e., 
$N_t\gg N$, the term $(j+1)^2$ starts dominating after a point.

\item Compute the history functions on the boundary segments of the domain $\Omega_i$ as
\begin{equation*}
\mathcal{B}^{j+1}_{a_1}(y_2)
=\sum_{k=1}^{j+1}\omega_{k}\varphi^{j+1-k,j+1}_{a_1}(y_2),\quad
\mathcal{B}^{j+1}_{a_2}(y_1)
=\sum_{k=1}^{j+1}\omega_{k}\varphi^{j+1,j+1-k}_{a_2}(y_1).
\end{equation*}
The cost of this step is $\bigO{(j+1)N}$.

\item Solve the linear system for the interior problem to compute the solution
at $(j+1)$-th time step 
\begin{equation*}
-\left(\alpha^{-2}_1\partial^2_{y_1}
+\alpha^{-2}_2\partial^2_{y_2}\right)u^{j+1}+u^{j+1}
=u^j,\quad \rho = 1/\Delta t.
\end{equation*}
using the following Robin-type TBCs at discrete level
\begin{equation*}
\partial_{n_{1}}u^{j+1}+\alpha_1 u^{j+1}=-\alpha_1\mathcal{B}^{j+1}_{a_1}(y_2),\quad
\partial_{n_{2}}u^{j+1}+\alpha_2 u^{j+1}=-\alpha_2\mathcal{B}^{j+1}_{a_2}(y_2).
\end{equation*}
Let the cost of solving the linear 2D system of the 
kind~\eqref{eq:linear-2d-system} be denoted by $\Theta_2(N_{\text{dim}})$ where 
the dimension is $N_{\text{dim}}=(N-1)^2$. The cost of solving the linear system
in this step works out to
be $\Theta_2(N_{\text{dim}})=\bigO{N^3}$ provided pivoting is not needed.
\item Update the 2D and 3D arrays storing the values of the auxiliary field with 
the $(j+1)$-th time step values using the fact that 
$\varphi^{j+1,j+1}(y_1,y_2)= u^{j+1}(y_1,y_2)$.
\end{enumerate}
For moderate value of $N_t$, i.e., $N_t=\bigO{N}$, the cost of computation at
each time-step is dominated by $\Theta_2(N_{\text{dim}})=\bigO{N^3}$ which 
corresponds to the linear system for the interior problem. The memory
requirement for the auxiliary functions in this case scales as
$\bigO{N^3}$ which is quite prohibitive.

In the regime $N_t\gg N$, say, $N_t=\bigO{N^2}$, the cost of computation scales
as $\bigO{N^4}$ which corresponds to cost of computing the history functions 
while the memory requirement scales as $\bigO{N^5}$ .

\subsubsection{NP--BDF1}
Let us recall that, on the boundary segments $\Gamma_{a_1}$ and $\Gamma_{a_2}$, 
the auxiliary fields are denoted by $\varphi_{k,a_1}(y_2,\tau_1,\tau_2)$ and
$\varphi_{k,a_2}(y_1,\tau_1,\tau_2)$, respectively,  while the auxiliary field 
at the corners is denoted by $\psi_{k,k',a_1,a_2}(\tau_1,\tau_2)$ where 
$a_1\in\{l,r\},\;a_2\in\{b,t\}$ and $k,k'=1,2,\ldots,M$. Let $N$ be the number of spatial nodes along 
$y_1$ as well as $y_2$, i.e., $N_1=N_2=N$. The storage of the auxiliary fields 
on the boundary segments require four 2D-arrays of size $M\times N$ and the 
storage of auxiliary field at corner points require four 2D-arrays of size 
$M\times M$  where $M (<N)$ is the order of the diagonal Pad\'e approximants used
in the NP method. The evolution of these fields are effectively local in time 
which circumvents the need to store the history of these fields in time.
The algorithmic steps for the solving the IBVP are enumerated below.
\begin{enumerate}[%
,leftmargin=*
,label={\bfseries Step \arabic*:}
]
\item Advance the following IVPs corresponding to auxiliary fields by one step
, say, $j\rightarrow (j+1)$ using the previously computed values of the auxiliary 
field on the segments (as described in Fig.~\ref{fig:IVP-auxi-pade}): 
\begin{equation*}
-\alpha^{-2}_1\partial_{y_1}^2\varphi^{j+1,j}_{k,a_2}+\varphi^{j+1,j}_{k,a_2}
=\varphi^{j,j}_{k,a_2},\quad
-\alpha^{-2}_1\partial_{y_2}^2\varphi^{j,j+1}_{k,a_1}+\varphi^{j,j+1}_{k,a_1}
=\varphi^{j,j}_{k,a_1}, \quad k=1,2,\ldots,M,
\end{equation*}
together with the boundary conditions 
\begin{equation*}
\partial_{n_2}\varphi_{k,a_1}^{j,j+1}+\alpha_2\varpi\varphi_{k,a_1}^{j,j+1}
=-\alpha_2\mathcal{B}^{j+1}_{k,a_1,a_2},\quad
\partial_{n_1}\varphi_{k,a_2}^{j+1,j}+\alpha_1\varpi\varphi_{k,a_2}^{j+1,j}
=-\alpha_1\mathcal{B}^{j+1}_{k,a_2,a_1},
\end{equation*}
where the history functions are given by
\begin{equation*}
\mathcal{B}^{j+1}_{k,a_1,a_2}=\sum_{k'=1}^M \Gamma_{k'}\psi_{k,k',a_1,a_2}^{j,j},
\quad \mathcal{B}^{j+1}_{k,a_2,a_1}=\sum_{k'=1}^M \Gamma_{k'}\psi_{k,k',a_2,a_1}^{j,j}
=\sum_{k'=1}^M \Gamma_{k'}\psi_{k',k,a_1,a_2}^{j,j}.
\end{equation*}
These boundary condition are taken from Sec.~\ref{sec:dtn-maps} where precise definition
of the quantities are provided. The cost of the steps outlined here works out to be 
$4M\Theta_1(N-1)+4M^2$ where the second part is the cost of computing
the history functions. The memory requirement is $4MN+4M^2$ size array.

\item Compute the history functions on the boundary segments of the domain $\Omega_i$ as
\begin{equation*}
\mathcal{B}^{j+1}_{a_1}(y_2)=\sum_{k=1}^M \Gamma_{k}\varphi_{k,a_1}^{j,j+1},\quad 
\mathcal{B}^{j+1}_{a_2}(y_1)=\sum_{k=1}^M \Gamma_{k}\varphi_{k,a_2}^{j+1,j}.
\end{equation*}
The cost of this step is $\bigO{MN}$.
\item Solve the linear system for the interior problem to compute the solution
at $(j+1)$-th time step 
\begin{equation*}
-\left(\alpha^{-2}_1\partial^2_{y_1}
+\alpha^{-2}_2\partial^2_{y_2}\right)u^{j+1}+u^{j+1}
=u^j,\quad \rho = 1/\Delta t.
\end{equation*}
using the following Robin-type TBCs at discrete level 
\begin{equation*}
\partial_{n_1} u^{j+1}+ \alpha_1\varpi u^{j+1}=-\alpha_1\mathcal{B}^{j+1}_{a_1}(y_2),\quad
\partial_{n_2} u^{j+1}+ \alpha_2\varpi u^{j+1}=-\alpha_2\mathcal{B}^{j+1}_{a_2}(y_1).
\end{equation*}
The cost of solving the linear system in this step works out to
be $\Theta_2(N_{\text{dim}})=\bigO{N^3}$ provided pivoting is not needed.
\item Update the 2D arrays storing the values of the auxiliary fields with 
the $(j+1)$-th time step values using the following update relations
\begin{equation*}
\left\{
\begin{aligned}
&\varphi^{j+1,j+1}_{k,a_1} 
=\frac{1}{(1+\ovl{\eta}^2_k)}\varphi^{j,j+1}_{k,a_1}
+\frac{1/\rho}{(1+\ovl{\eta}^2_k)}u^{j+1},\\
&\varphi^{j+1,j+1}_{k,a_2} 
=\frac{1}{(1+\ovl{\eta}^2_k)}\varphi^{j+1,j}_{k,a_2}
+\frac{1/\rho}{(1+\ovl{\eta}^2_k)}u^{j+1},
\end{aligned}\right.
\quad\text{and}\quad
\left\{
\begin{aligned}
&\psi^{j,j+1}_{k,k',a_1a_2} 
=\frac{1}{(1+\ovl{\eta}^2_{k'})}\psi^{j,j}_{k,k',a_1,a_2}
+\frac{1/\rho}{(1+\ovl{\eta}^2_{k'})}\varphi_{k,a_1}^{j,j+1},\\
&\psi^{j+1,j+1}_{k,k',a_1a_2} 
 =\frac{1}{(1+\ovl{\eta}^2_{k})}\psi^{j,j+1}_{k,k',a_1,a_2}
 +\frac{1/\rho}{(1+\ovl{\eta}^2_{k})}\varphi_{k',a_2}^{j+1,j+1}.
\end{aligned}\right.
\end{equation*}
\end{enumerate}
From the discussion above, it is evident that the overall 
cost of computation is completely independent of the time-step and the DtN-maps
get implemented at the cost of $\bigO{MN}$. The memory
requirement for the auxiliary functions also remain independent of time-step
with an estimate of $\bigO{MN}$. Therefore, the overall complexity is 
dominated by $\Theta_2(N_{\text{dim}})=\bigO{N^3}$ which 
corresponds to the linear system for the interior problem. 

\begin{figure}[!ht]
\begin{center}
\def\myscale{1}
\includegraphics[width=\textwidth]{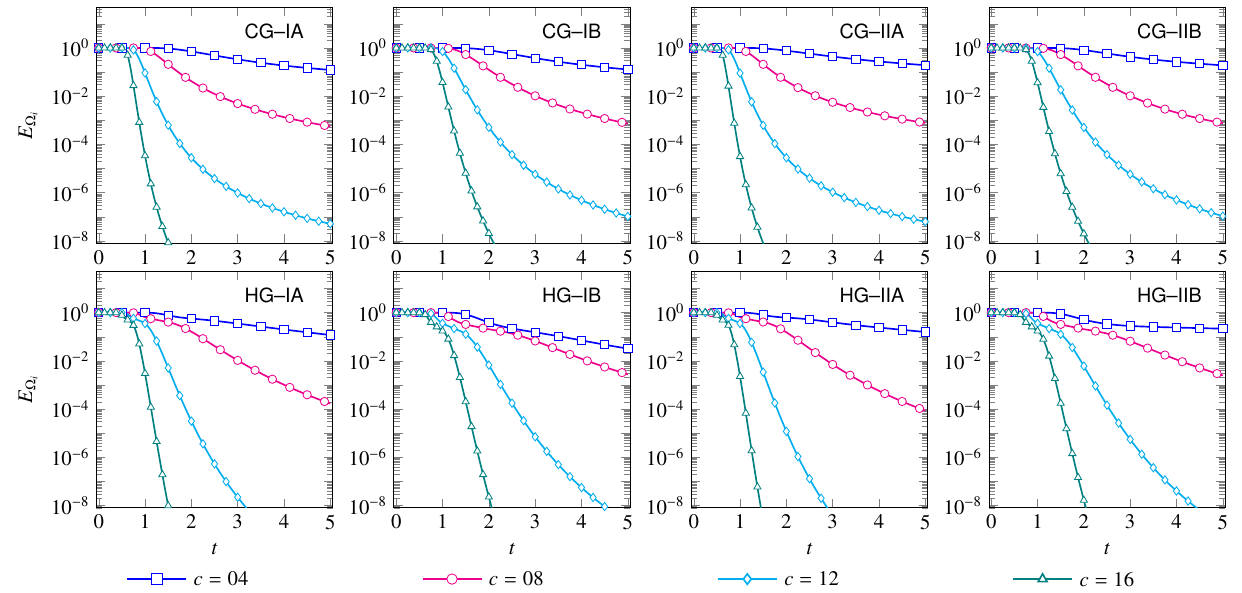}
\end{center}
\caption{\label{fig:wp2d-energy-content}The figure shows the evolution of the 
relative energy content as defined in~\eqref{eq:cg2d-energy-content} of the 
chirped-Gaussian and the Hermite-Gaussian profiles considered in 
Table~\ref{tab:cg2d} and Table~\ref{tab:hg2d}. Here the computational domain is 
$\Omega_i=(-10,10)^2$.}
\end{figure}
\section{Numerical Experiments}\label{sec:numerical-experiments}
In this section, we will carry out the extensive numerical tests to showcase the 
accuracy of the numerical schemes developed in this work to solve the 
IBVP in~\eqref{eq:2D-SE-CT}. We start by analyzing the behaviour of the exact solutions
for the IBVP under consideration followed by studying the error evolution behaviour
and convergence analysis of the numerical schemes.

\subsection{Exact solutions}\label{sec:exact-solution}
The exact solutions admissible for our numerical experiments are such that the 
initial profile must be effectively supported within the computational domain.
We primarily consider the wavepackets which are a modulation of a Gaussian
envelop so that the requirement of effective initial support can be easily met.
In this class of solutions, we consider the Chirped-Gaussian and the
Hermite-Gaussian profiles.

\def\arraystretch{2}
\setlength{\tabcolsep}{1mm}
\begin{table}[htb]
\centering
\caption{\label{tab:cg2d} Chirped-Gaussian profile with $A_0=2$ and 
$c_0\in\{4,8,12,16\}$.}
\begin{tabular}{m{10mm}m{5mm}m{30mm}m{20mm}m{40mm}}
\hline
\multicolumn{5}{c}{
$G\left(\vv{x},t\right) 
= A_0\sum_{j=1}^n G(\vv{x},t;\vv{a}_j,\vv{b}_j,\vv{c}_j),
\quad\vv{c}_j=c_0(\cos\theta_j,\sin\theta_j)$}\\
\hline
Type &$n$ & $(\vv{a}_j\in\field{R}^2)_{j=1}^n$
&$(\vv{b}_j\in\field{R}^2)_{j=1}^n$ & $\vs{\theta}\in\field{R}^n$\\
\hline
IA,\newline{IB} & $2$ & 
$\vv{a}_1=(1/2.5,1/2.4)$,\newline$\vv{a}_2=(1/2.3,1/2.2)$ & $\vv{b}_j=(1/2)\vv{1}$
&IA:~$\vs{\theta}_A=(0,\pi)$,\newline{IB}:~$\vs{\theta}_B=\vv{\theta}_A+(\pi/4)\vv{1}$\\
\hline
IIA,\newline{IIB} & $4$ & 
$\vv{a}_1 = (1/2.5,1/2.4)$,\newline 
$\vv{a}_2 = (1/2.3,1/2.2)$,\newline
$\vv{a}_3 = (1/2.7,1/2.6)$,\newline
$\vv{a}_4 = (1/2.2,1/2.5)$
& $\vv{b}_j=(1/2)\vv{1}$
&IIA:~$\vs{\theta}_A=(0,\pi/2,\pi,3\pi/2)$,\newline{IIB}:~$\vs{\theta}_B=\vv{\theta}_A+(\pi/4)\vv{1}$\\
\hline
\end{tabular}
\end{table}

\subsubsection{Chirped-Gaussian profile}
Using the function defined by
\begin{equation}
\mathcal{G}(x,t;a,b)=\frac{1}{\sqrt{1+4i(a+ib)t}}
\exp\left[-\frac{(a+ib)}{1+4i(a+ib)t}x^2\right],\quad a>0,\;b\in\field{R},
\end{equation}
one can define the family of solutions referred to as chirped-Gaussian profile by
\begin{equation}
G(\vv{x},t;\vv{a},\vv{b},\vv{c})
=\mathcal{G}(x_1-c_1t,t;a_1,b_1)\mathcal{G}(x_2-c_2t,t;a_2,b_2)
\exp\left(+i\frac{1}{2}\vv{c}\cdot\vv{x}-i\frac{1}{4}\vv{c}\cdot\vv{c}\,t\right),
\end{equation}
where $\vv{a}\in\field{R}^2_+$ determines the effective support of the 
profile at $t=0$, $\vv{b}\in\field{R}^2$ is the \emph{chirp} parameter and 
$\vv{c}\in\field{R}^2$ is the velocity of the profile. Using linear 
combination, one can further define a more general family of solutions with 
parameters $A_0,c_0\in\field{R}$, $\vs{\theta}\in\field{R}^n$, 
$(\vv{a}_j\in\field{R}^2_+)_{j=1}^n$ and 
$(\vv{b}_j\in\field{R}^2)_{j=1}^n$ given by
\begin{equation}
G\left(\vv{x},t;\;c_0,A_0;\;(\vv{a}_j)_{j=1}^n,(\vv{b}_j)_{j=1}^n,\vv{\theta}\right) 
= A_0\sum_{j=1}^n G(\vv{x},t;\vv{a}_j,\vv{b}_j,\vv{c}_j),
\quad\vv{c}_j=c_0(\cos\theta_j,\sin\theta_j).
\end{equation}
Let the constant vector $(1,1,\ldots)\in\field{R}^n$ be denoted by $\vv{1}$,
then the specific values of the parameters of the solutions used in the
numerical experiments can be summarized as in Table~\ref{tab:cg2d}. The energy 
content of the profile within the computational domain $\Omega_i$ 
over time is
\begin{equation}\label{eq:cg2d-energy-content}
E_{\Omega_i}(t)
=\left.{\int_{\Omega_i}|G(\vv{x},t)|^2d^2\vv{x}}\right/{\int_{\Omega_i}|G(\vv{x},0)|^2d^2\vv{x}},
\quad t\geq0.
\end{equation}
The profiles are chosen with non-zero speed $c_0$ so that the field hits the
boundary of $\Omega_i$, In case of a square computational domain, the type `A' class of 
solutions are directed to the segments of the boundary normally while type `B' 
class of solutions are directed to the corners. The behaviour of the 
$E_{\Omega_i}(t)$ for $t\in[0,5]$ is shown in Fig.~\ref{fig:wp2d-energy-content}.
\def\arraystretch{2}
\setlength{\tabcolsep}{1mm}
\begin{table}[htb!]
\centering
\caption{\label{tab:hg2d} Hermite-Gaussian profile with $A_0=2$ and 
$c_0\in\{4,8,12,16\}$.}
\begin{tabular}{m{10mm}m{5mm}m{30mm}m{20mm}m{40mm}}
\hline
\multicolumn{5}{c}{
$G\left(\vv{x},t\right) 
= A_0\sum_{j=1}^n G(\vv{x},t;\vv{m}_j,\vv{a}_j,\vv{c}_j),\;
\vv{c}_j=c_0(\cos\theta_j,\sin\theta_j)$
}\\
\hline
Type &$n$ & $(\vv{a}_j\in\field{R}^2)_{j=1}^n$
&$(\vv{m}_j\in\field{N}_0^2)_{j=1}^n$ & $\vs{\theta}\in\field{R}^n$\\
\hline
IA,\newline{IB} & $2$ & 
$\vv{a}_1=(1/2.5,1/2.4)$,\newline$\vv{a}_2=(1/2.3,1/2.2)$ & 
$\vv{m}_1 = (1,2)$,\newline$\vv{m}_2 = (2,1)$
&$\vs{\theta}_A=(0,\pi)$,\newline%
$\vs{\theta}_B=\vv{\theta}_A+(\pi/4)\vv{1}$\\
\hline
IIA,\newline{IIB} & $4$ & 
$\vv{a}_1 = (1/2.5,1/2.4)$,\newline 
$\vv{a}_2 = (1/2.3,1/2.2)$,\newline
$\vv{a}_3 = (1/2.7,1/2.6)$,\newline
$\vv{a}_4 = (1/2.2,1/2.5)$ &
$\vv{m}_1 = (1,2)$,\newline
$\vv{m}_2 = (2,1)$,\newline
$\vv{m}_3 = (2,1)$,\newline
$\vv{m}_4 = (1,2)$
&$\vs{\theta}_A=(0,\pi/2,\pi,3\pi/2)$,\newline%
$\vs{\theta}_B=\vv{\theta}_A+(\pi/4)\vv{1}$\\
\hline
\end{tabular}
\end{table}

\begin{figure}[!h]
\begin{center}
\includegraphics[width=\textwidth]{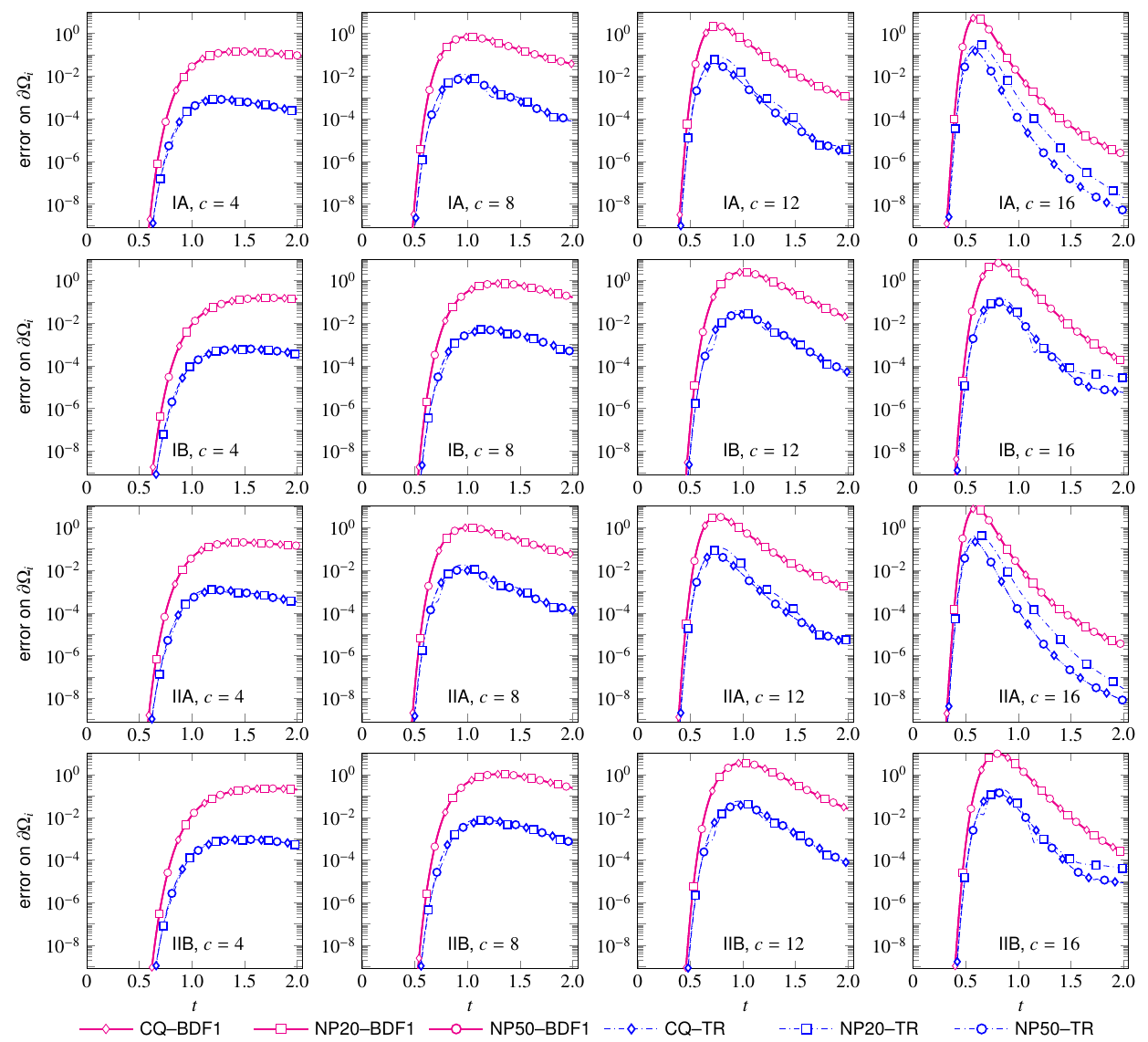}
\end{center}
\caption{\label{fig:wpcg2d-mt}The figure shows the behaviour of error involved in 
the discretization of the boundary map alone for the chirped-Gaussian profile with
different values of the speed `c' (see Table~\ref{tab:cg2d}). The numerical 
parameters and the labels are described in Sec.~\ref{sec:tests-mt} where the 
error is quantified by~\eqref{eq:error-dtn}.}
\end{figure}

\begin{figure}[!h]
\begin{center}
\includegraphics[width=\textwidth]{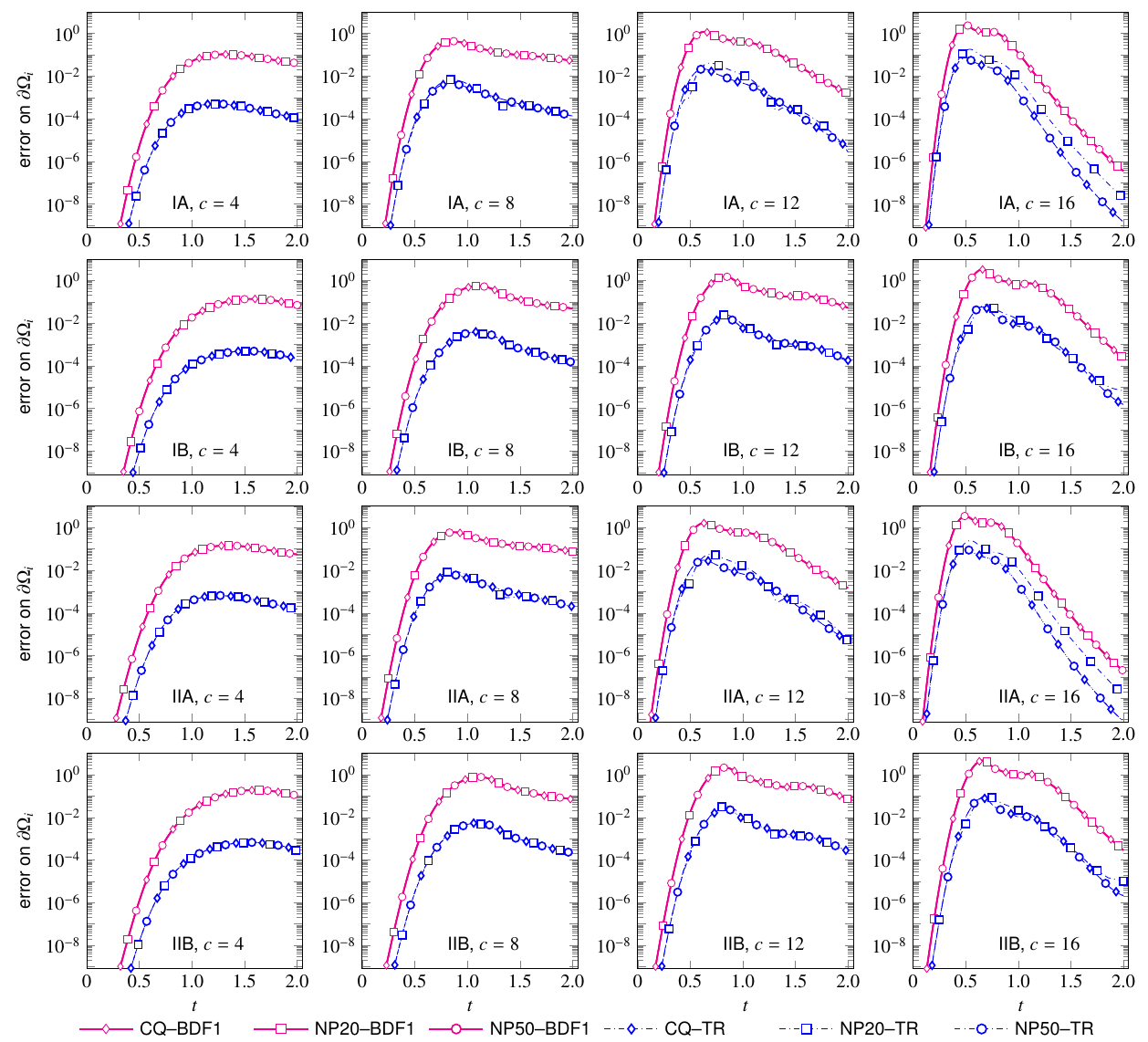}
\end{center}
\caption{\label{fig:wphg2d-mt}The figure shows the behaviour of error involved in 
the discretization of the boundary map alone for the Hermite-Gaussian profile 
with different values of the speed `c' (see Table~\ref{tab:hg2d}). The numerical 
parameters and the labels are described in Sec.~\ref{sec:tests-mt} where the 
error is quantified by~\eqref{eq:error-dtn}.}
\end{figure}

\begin{figure}[!h]
\begin{center}
\includegraphics[width=\textwidth]{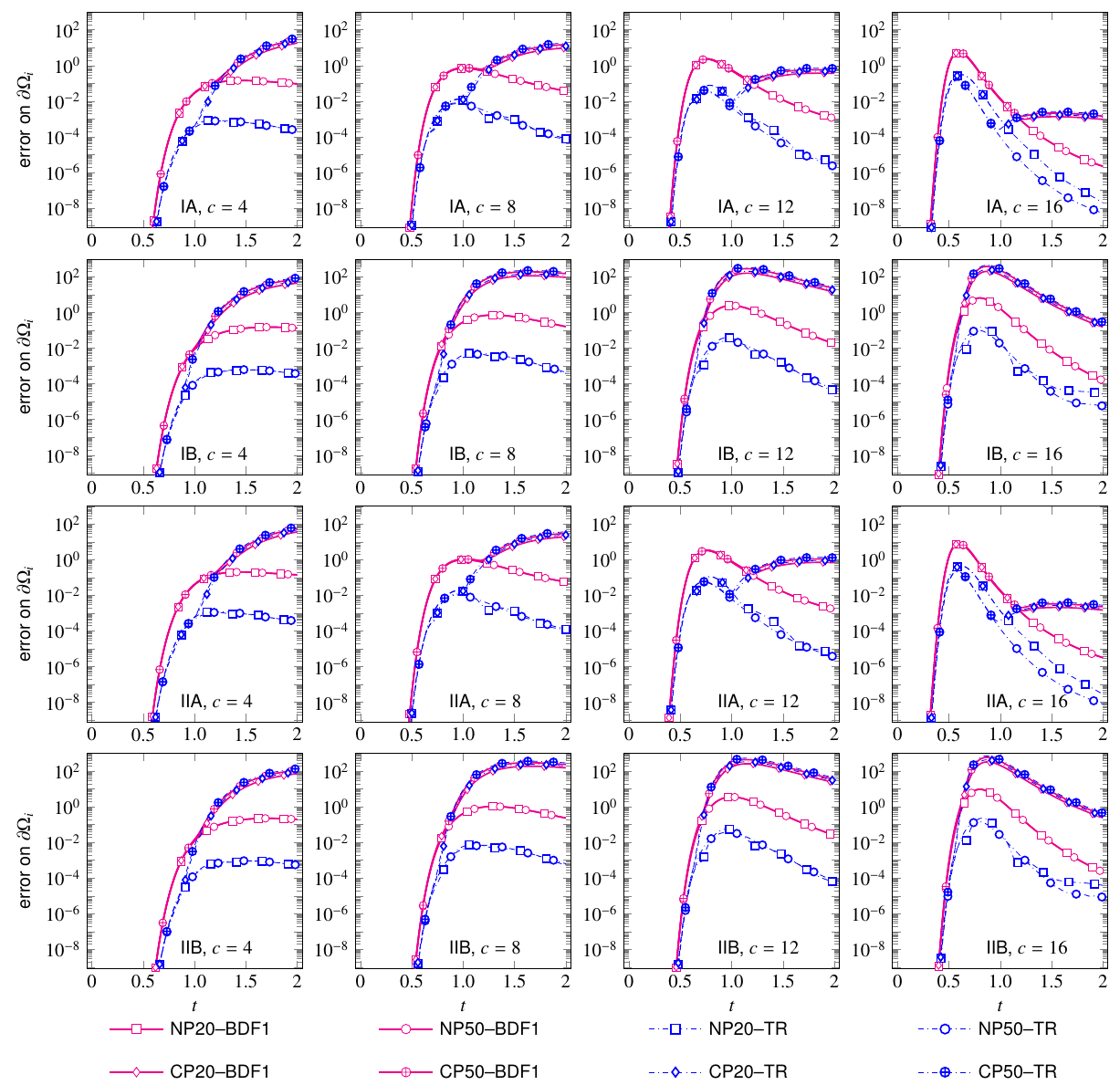}
\end{center}
\caption{\label{fig:wpcg2d-mt-menza}\color{myrem}The figure shows a comparison of the novel 
Pad\'e and Menza's approach towards realizing the DtN maps for the 
chirped-Gaussian profile with
different values of the speed `c' (see Table~\ref{tab:cg2d}). The numerical 
parameters and the labels are described in Sec.~\ref{sec:tests-mt} where the 
error is quantified by~\eqref{eq:error-dtn}.}
\end{figure}

\begin{figure}[!h]
\begin{center}
\includegraphics[width=\textwidth]{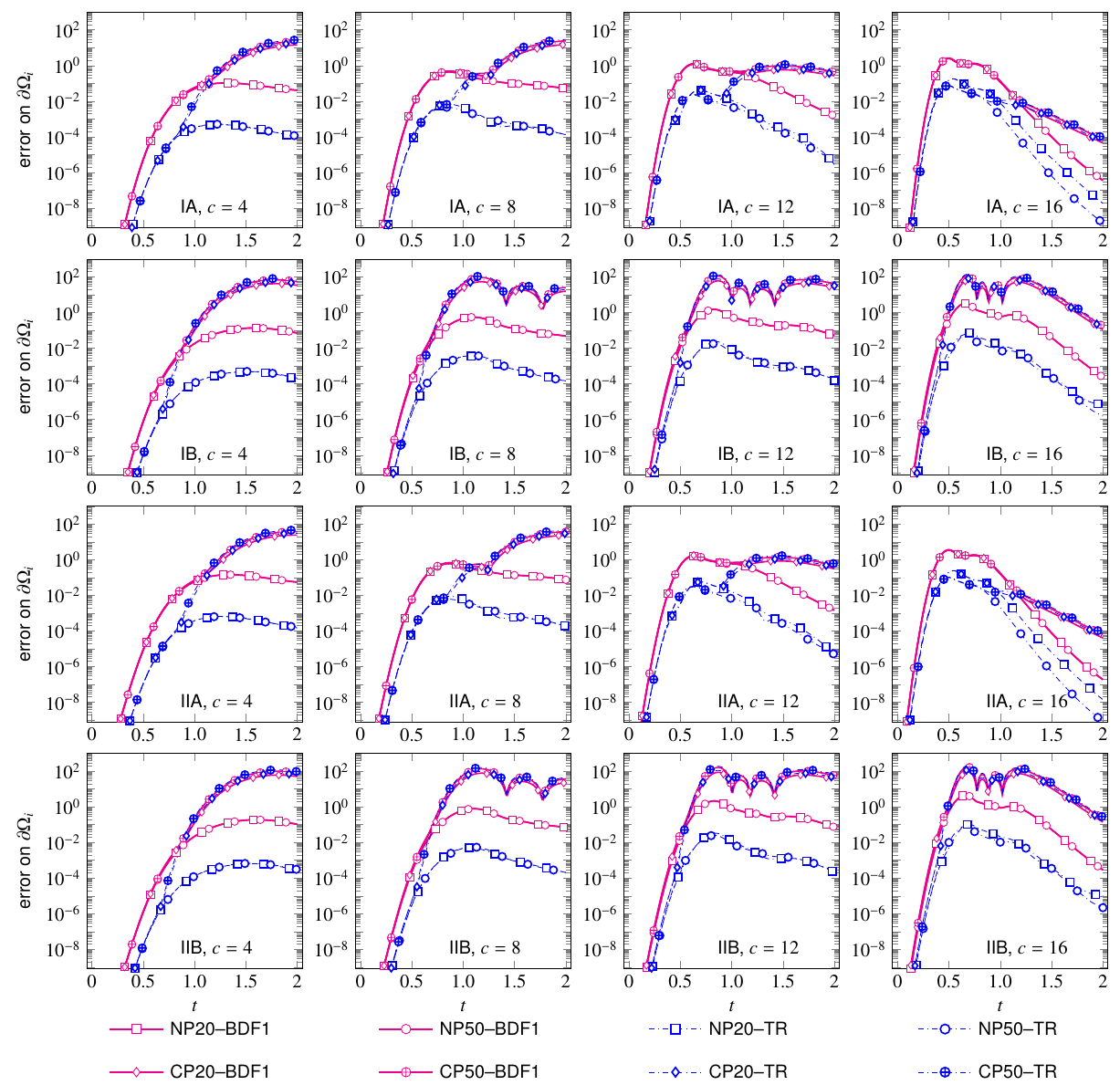}
\end{center}
\caption{\label{fig:wphg2d-mt-menza}\color{myrem}The figure shows a comparison of the 
novel Pad\'e and Menza's approach towards realizing the DtN maps for the 
Hermite-Gaussian profile with different values of the speed `c' 
(see Table~\ref{tab:hg2d}). The numerical parameters and the labels are 
described in Sec.~\ref{sec:tests-mt} where the error is quantified 
by~\eqref{eq:error-dtn}.}
\end{figure}
\begin{figure}[!htb]
\begin{center}
\includegraphics[width=\textwidth]{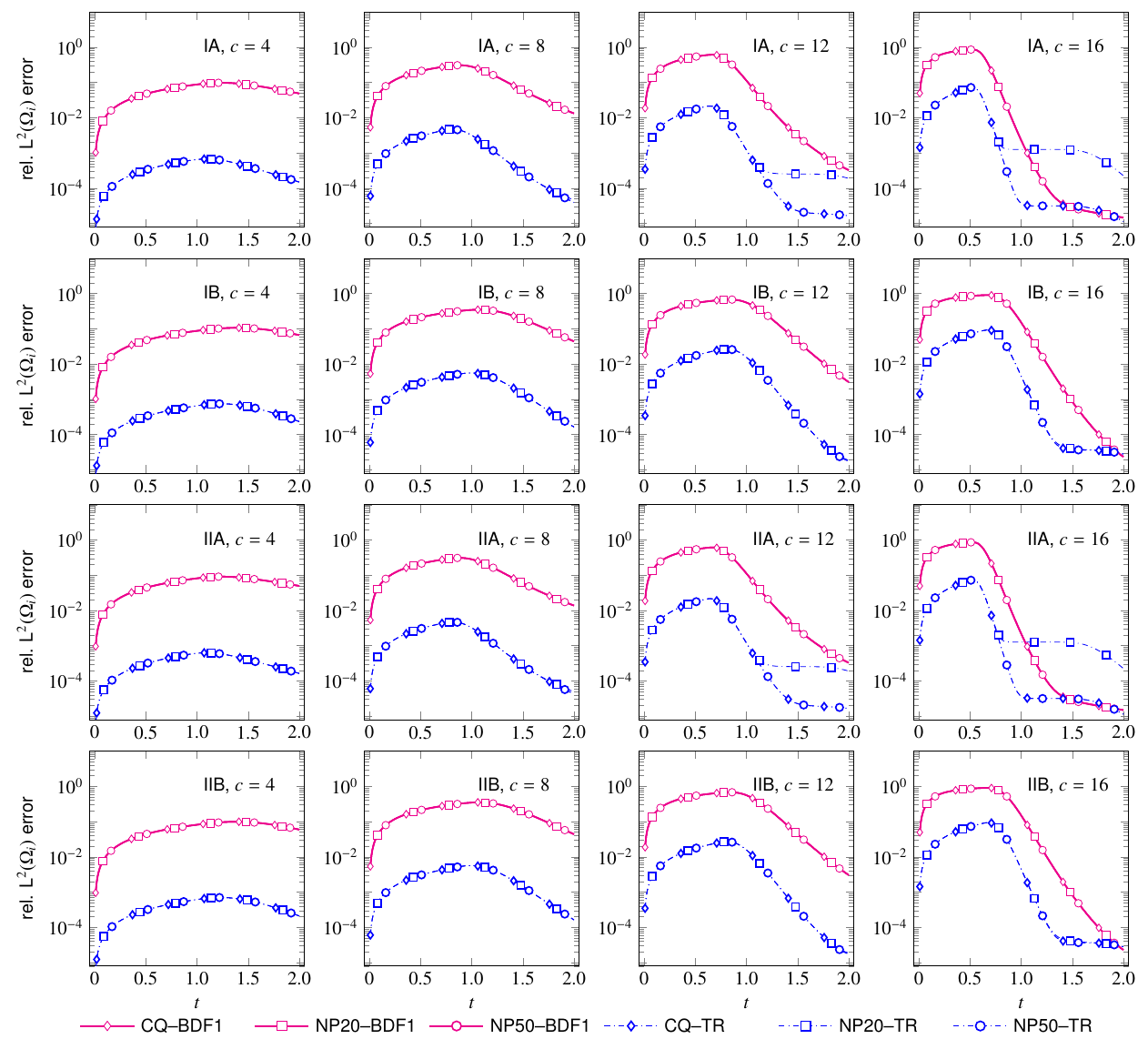}
\end{center}
\caption{\label{fig:wpcg2d-ee}The figure shows a comparison of evolution of error 
in the numerical solution of the IBVP~\eqref{eq:2D-SE-CT} with the various TBCs 
for the chirped-Gaussian profile with different values of the speed `c' 
(see Table~\ref{tab:cg2d}). The numerical parameters and the labels are described in 
Sec.~\ref{sec:tests-ee} where the error is quantified by~\eqref{eq:error-ibvp}.}
\end{figure}

\subsubsection{Hermite-Gaussian profile}
Consider the class of normalized Hermite-Gaussian functions defined by 
\begin{equation}
\mathcal{G}_{m}(x,t;a)=\gamma_m^{-1}H_{m}\left(\frac{\sqrt{2a}\,x }{w(t)}\right)
\sqrt{\frac{\mu(t)}{a}}
\exp\left[-\mu(t)x^2-i\,m\theta(t)\right],
\quad m\in\field{N}_0=\{0,1,\ldots\},
\end{equation}
where $a>0$,
\begin{equation}
w(t)=\sqrt{1+(4at)^2},\quad
\frac{1}{\mu(t)}=\frac{1}{a}+i{4t}=\frac{1}{a}w(t)\exp[i\theta(t)],
\end{equation}
and the normalization factor is given by 
$\gamma^2_{m}=2^m(m!)\sqrt{\pi}(2a)^{-1/2}$.

The Hermite polynomials are evaluated using the following relations
\begin{equation}
H_{n+1}(x)=2xH_n(x)-2nH_{n-1}(x),\quad H_{n+1}'(x)=2(n+1)H_{n}(x),
\end{equation}
with $H_0(x) = 1$ and $H_1(x)=2x$, These relations can be used to compute the
first spatial derivative of the Hermite-Gaussian functions:
\begin{equation}
\partial_x\mathcal{G}_{m}(x,t;a)=
-\sqrt{(m+1)a}\mathcal{G}_{m+1}(x,t;a)+\sqrt{ma}\mathcal{G}_{m-1}(x,t;a),
\quad m>0,
\end{equation}
with $\partial_x\mathcal{G}_{0}(x,t;a)=-\sqrt{a}\mathcal{G}_{1}(x,t;a)$. This is
required for testing our Dirichlet-to-Neumann maps. Using these functions, we
can define a family of solutions referred to as Hermite-Gaussian profile by
\begin{equation}
G(\vv{x},t;\vv{m},\vv{a},\vv{c})
=\mathcal{G}_{m_1}(x_1-c_1t,t;a_1)\mathcal{G}_{m_2}(x_2-c_2t,t;a_2)
\exp\left(+i\frac{1}{2}\vv{c}\cdot\vv{x}-i\frac{1}{4}\vv{c}\cdot\vv{c}\,t\right),
\end{equation}
where  $\vv{m}\in\field{N}_0^2$ is the order parameter, 
$\vv{a}\in\field{R}^2_+$ determines the effective support of the profile 
at $t=0$ and $\vv{c}\in\field{R}^2$ is the velocity of the profile. Once again,
using linear combination, one can further define a more general family of solutions with 
parameters $A_0,c_0\in\field{R}$, $\vs{\theta}\in\field{R}^n$, 
$(\vv{m}_j\in\field{N}^2_0)_{j=1}^n$ and
$(\vv{a}_j\in\field{R}^2_+)_{j=1}^n$ given by
\begin{equation}
G\left(\vv{x},t;\;c_0,A_0;\;(\vv{m}_j)_{j=1}^n,(\vv{a}_j)_{j=1}^n,\vv{\theta}\right) 
= A_0\sum_{j=1}^n G(\vv{x},t;\vv{m}_j,\vv{a}_j,\vv{c}_j),
\quad\vv{c}_j=c_0(\cos\theta_j,\sin\theta_j).
\end{equation}
As in the last case, we set $\vv{1}=(1,1,\ldots)\in\field{R}^n$, then the specific 
values of the parameters of the solutions used in the
numerical experiments can be summarized as in Table~\ref{tab:hg2d}. The energy 
content of the profile within the computational domain $\Omega_i$ 
over time is $E_{\Omega_i}(t)$ as defined by~\eqref{eq:cg2d-energy-content} with
the appropriate profile in the integrand. The profiles are chosen with non-zero 
speed $c_0$ so that the field hits the
boundary of $\Omega_i$. In case of a square computational domain, the type `A' 
class of solutions are directed to the segments of the boundary normally while 
type `B' class of solutions are directed to the corners. The behaviour of the 
$E_{\Omega_i}(t)$ for $t\in[0,5]$ is shown in Fig.~\ref{fig:wp2d-energy-content}.

\begin{figure}[!htb]
\begin{center}
\includegraphics[width=\textwidth]{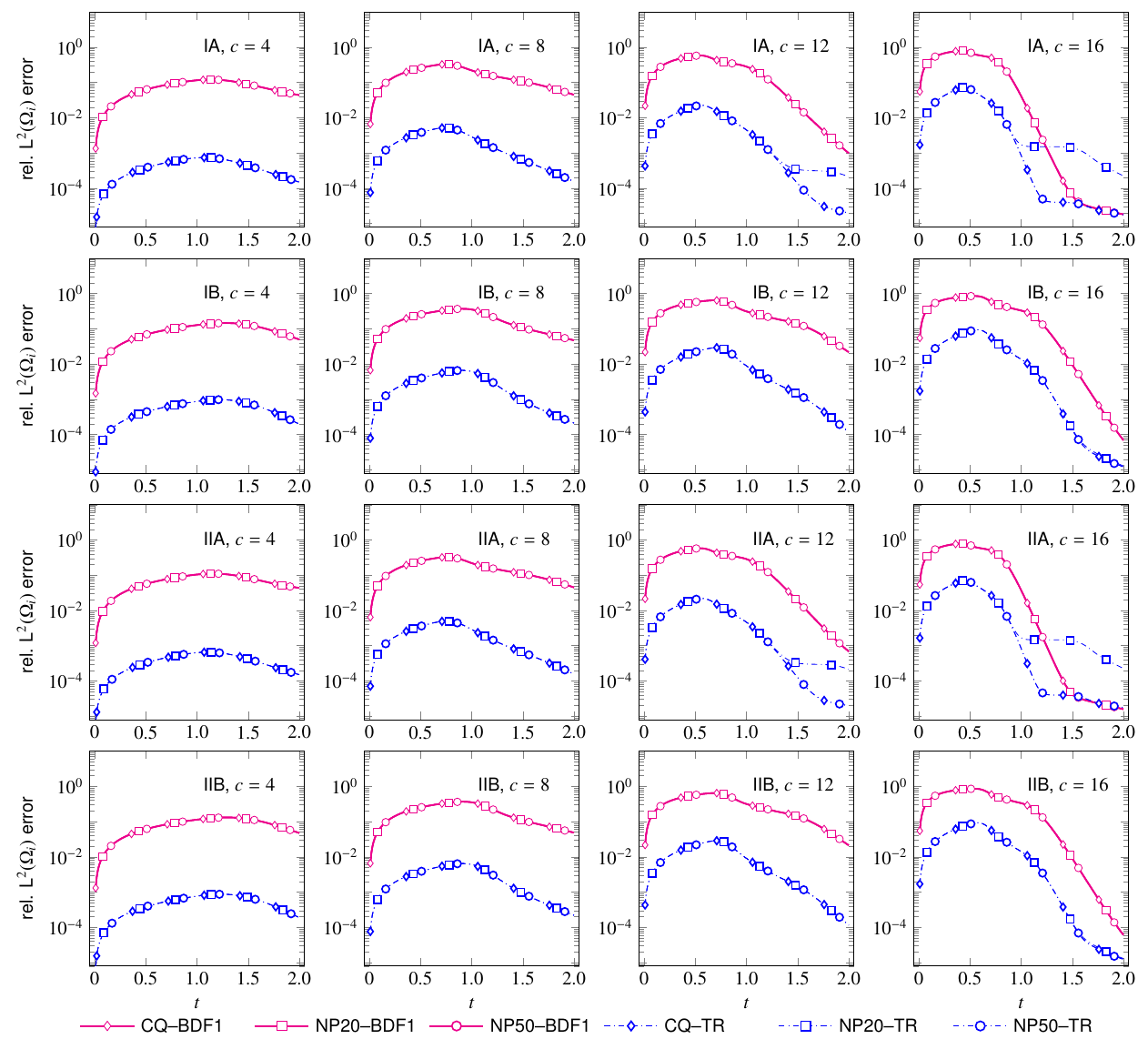}
\end{center}
\caption{\label{fig:wphg2d-ee}The figure shows a comparison of evolution of error 
in the numerical solution of the IBVP~\eqref{eq:2D-SE-CT} with the various TBCs 
for the Hermite-Gaussian profile with different values of the speed `c' 
(see Table~\ref{tab:hg2d}). The numerical parameters and the labels are described in 
Sec.~\ref{sec:tests-ee} where the error is quantified by~\eqref{eq:error-ibvp}.}
\end{figure}

\subsection{Tests for DtN maps}\label{sec:tests-mt}
Prior to the full implementation of the numerical scheme, one can test the 
accuracy of the DtN boundary maps where we quantify the error on
$\partial\Omega_i$ by
\begin{equation}\label{eq:error-dtn}
e(t_j)=
\left(\int_{\partial\Omega_i}
\left|\partial_{n}u(\vv{x},t_j)-[\partial_{n}u(\vv{x},t_j)]_{\text{num.}}
\right|^2d\sigma(\vv{x})\right)^{1/2},\quad t_j\in[0,T_{max}],
\end{equation}
where $d\sigma(\vv{x})$ denotes the measure along the boundary of the
computational domain $\Omega_i$. Along the spatial dimension, the integral
above will be approximated by Gauss quadrature over LGL-points~\footnote{The
quadrature nodes and the corresponding weights are computed using quad precision
in order to avoid any numerical artifacts due to round-off errors.}.
The quantity 
$[\partial_{n}u(\vv{x},t_j)]_{\text{num.}}$ refers to the numerical
approximation of the Neumann datum from the Dirichlet datum. All the intermediate 
quantities involved in this process is computed according to a numerical scheme
developed as part of the numerical realization of the DtN map.

In this work, we have developed two discrete versions of the DtN maps, namely, 
CQ and NP. Each of these methods have a variant determined by the choice of 
one-step method used in the temporal discretization so that the complete list of
methods to be tested can be labelled as CQ--BDF1, NP--BDF1 (corresponding to the
one-step method BDF1), and, CQ--TR, NP--TR (corresponding to the
one-step method TR). For the Pad\'e approximant based method `NP', we
distinguish the diagonal approximant of order $20$ from that of $50$ via the
labels `NP20' and `NP50'. The numerical parameters used in this section
are summarized in Table~\ref{tab:mt-params}. The exact solutions used are 
described in Sec.~\ref{sec:exact-solution}. Let us note that the accuracy of the 
DtN maps do not reflect the accuracy of the numerical solution computed with
the TBCs on account of the fact that the interior problem together with the discretized 
DtN maps may turn out to be unstable. 

\def\arraystretch{2}
\setlength{\tabcolsep}{1mm}
\begin{table}[!h]
\centering
\caption{\label{tab:mt-params} Numerical parameters for testing DtN maps}
\begin{tabular}{m{80mm}m{50mm}}\hline
Computational domain ($\Omega_i$) & $(-10,10)\times(-10,10)$\\\hline
Maximum time ($T_{max}$)          & $2$\\\hline
No. of time-steps ($N_t$)         & $1000+1$\\\hline
Time-step ($\Delta t$)            & $2\times 10^{-3}=T_{max}/(N_t-1)$\\\hline
Number of LGL-points ($N+1$)        & $200$\\\hline
\end{tabular}
\end{table}

The numerical results for error on $\partial\Omega_i$ corresponding to the 
chirped-Gaussian and the Hermite-Gaussian profiles are shown in Fig.~\ref{fig:wpcg2d-mt} 
and Fig.~\ref{fig:wphg2d-mt}, respectively. It can be seen that the diagonal 
Pad\'e approximant based methods, namely, NP20 and NP50, perform equally well as 
compared to that of the convolution quadrature based method. Moreover, the TR 
methods perform better than the BDF1 methods which is clear from the error peaks 
in figures~\ref{fig:wpcg2d-mt} and~\ref{fig:wphg2d-mt}. The accuracy deteriorates 
with increase in the value of $c$. For faster moving profiles corresponding to 
IA and IIA type solutions, NP-20 becomes slightly less accurate compared to 
CQ and NP-50. There is no significant difference in the results for error on 
$\partial\Omega_i$ corresponding to the chirped-Gaussian profile and the 
Hermite-Gaussian profile.  

{\color{myrem}Before concluding this section, we would like to present evidence of the phenomenon
discussed in Sec.~\ref{sec:CT-CP} with regard to deficiencies in the Menza's approach.
We can test the accuracy of the DtN boundary maps realized according to Menza's approach 
using the error metric defined in~\eqref{eq:error-dtn}. The discrete scheme 
for the boundary maps is labelled as `CP--' followed by the time-stepping method 
used for temporal derivatives in the realization of the boundary maps. 
Here, Menza's approach is labeled as `CP' which stands for the conventional Pad\'e approach 
. Moreover, we distinguish the diagonal Pad\'e approximant
of order $20$ from that of $50$ via the labels `CP20' and `CP50'. 
The numerical results showing a comparison of the novel Pad\'e and Menza's approach
corresponding to the chirped-Gaussian and the 
Hermite-Gaussian profiles are shown in Fig.~\ref{fig:wpcg2d-mt-menza} 
and Fig.~\ref{fig:wphg2d-mt-menza}, respectively. It can be seen that the NP 
methods perform superior to the CP methods which is clear from the error peaks 
in figures~\ref{fig:wpcg2d-mt-menza} and~\ref{fig:wphg2d-mt-menza}. 
Moreover, CP methods perform poorly in the case of IB and IIB type solutions
where the profiles are directed towards the corners of the computational domain.
This observation points towards the fact that Menza's approach failed to handle 
corners of the domain adequately. 
For faster moving profiles corresponding to IA and IIA type solutions 
(directed towards the segments), error in the CP methods agree upto some extent 
with that of the NP methods.
The error behaviour on $\partial\Omega_i$ with Menza's approach becomes slightly 
oscillatory in case of the Hermite-Gaussian profile as evident in Fig.~\ref{fig:wphg2d-mt-menza}.
On account of the poor accuracy of the Menza's approach as demonstrated here,
we do not consider the full implementation of the IBVP in~\eqref{eq:2D-SE-CT} 
within the context of Menza's approach.
}

\begin{figure}[!htb]
\begin{center}
\includegraphics[width=\textwidth]{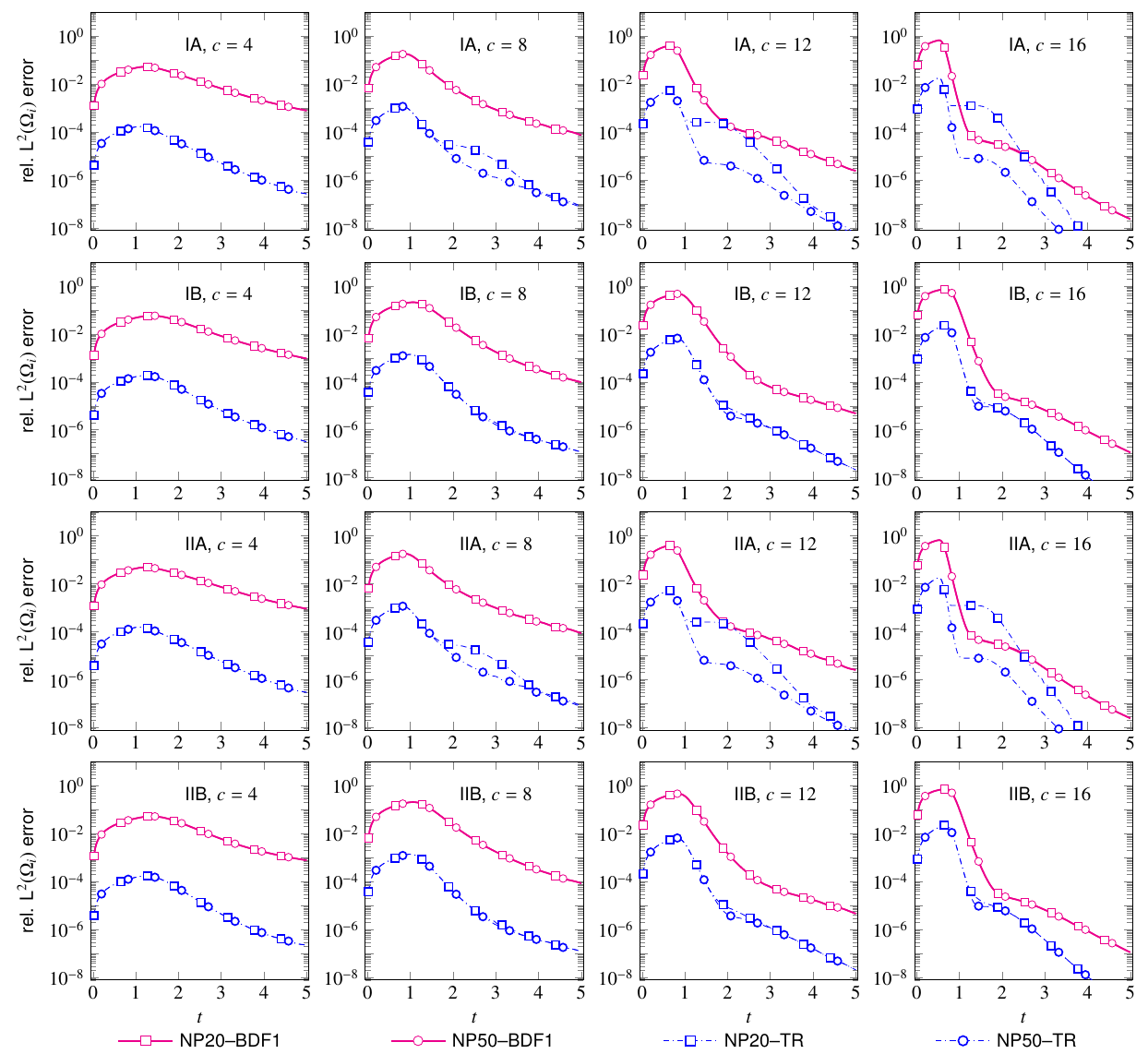}
\end{center}
\caption{\label{fig:wpcg2d-ee-np}The figure shows a comparison of evolution of error 
in the numerical solution of the IBVP~\eqref{eq:2D-SE-CT} with the TBCs developed using
novel Pad\'e approach for the chirped-Gaussian profile with different values of the speed `c' 
(see Table~\ref{tab:cg2d}). The numerical parameters and the labels are described in 
Sec.~\ref{sec:tests-ee} where the error is quantified by~\eqref{eq:error-ibvp}.}
\end{figure}

\begin{figure}[!htb]
\begin{center}
\includegraphics[width=\textwidth]{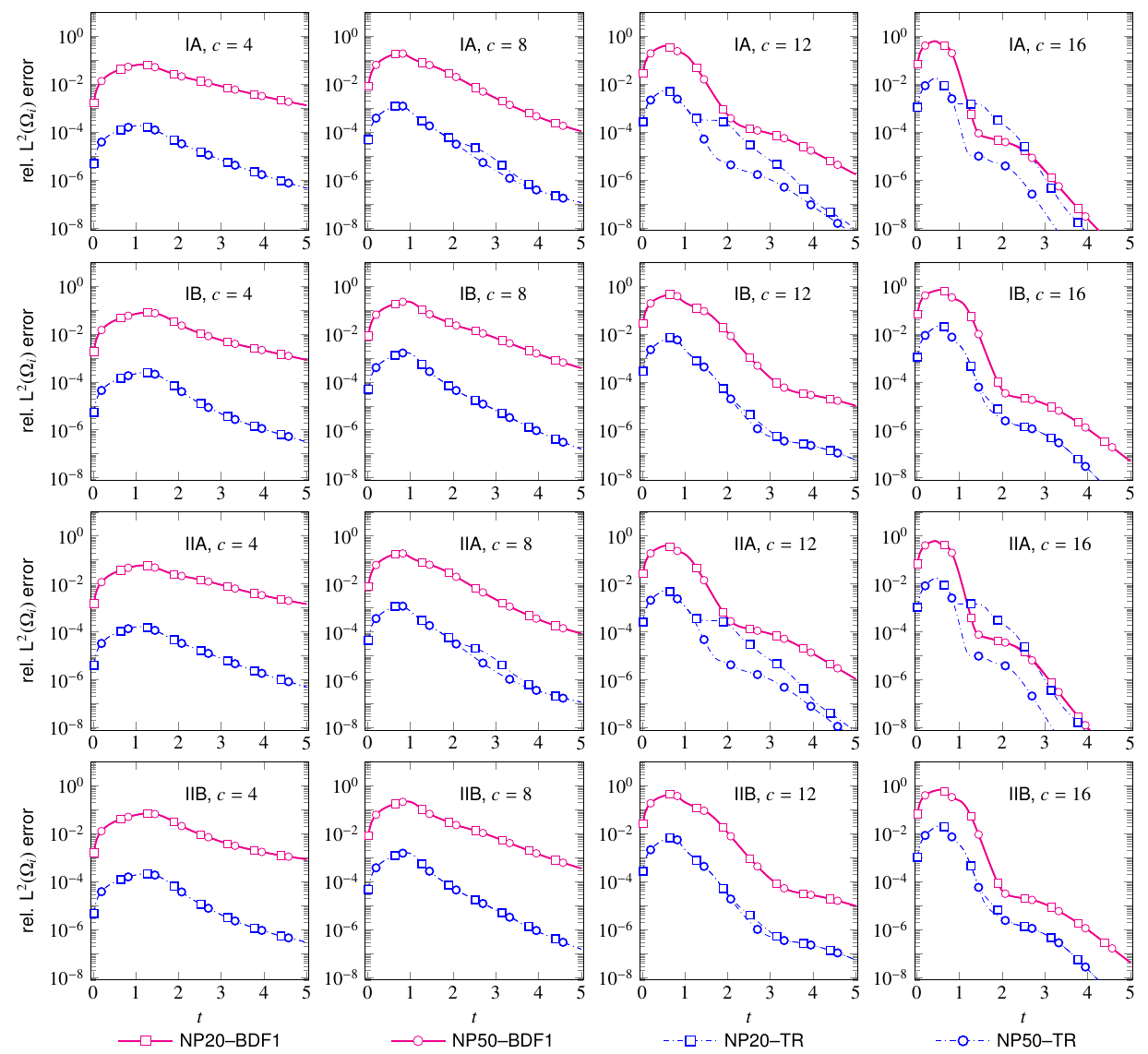}
\end{center}
\caption{\label{fig:wphg2d-ee-np}The figure shows a comparison of evolution of error 
in the numerical solution of the IBVP~\eqref{eq:2D-SE-CT} with the TBCs developed using
novel Pad\'e approach for the Hermite-Gaussian profile with different values of the speed `c' 
(see Table~\ref{tab:hg2d}). The numerical parameters and the labels are described in 
Sec.~\ref{sec:tests-ee} where the error is quantified by~\eqref{eq:error-ibvp}.}
\end{figure}

\subsection{Tests for evolution error}\label{sec:tests-ee}
In this section, we consider the IBVP in~\eqref{eq:2D-SE-CT} where the initial 
condition corresponds to the exact solutions described in 
Sec.~\ref{sec:exact-solution}. The numerical solution is labelled according to 
that of the boundary maps described in Sec.~\ref{sec:tests-mt}. The error in 
the evolution of the profile computed numerically is quantified by the 
relative $\fs{L}^2(\Omega_i)$-norm: 
\begin{equation}\label{eq:error-ibvp}
e(t_j)=
\left.\left(\int_{\Omega_i}
\left|u(\vv{x},t_j)-[u(\vv{x},t_j)]_{\text{num.}}
\right|^2d^2\vv{x}\right)^{1/2}\,\right/
\left(\int_{\Omega_i}\left|u_0(\vv{x})\right|^2d^2\vv{x}\right)^{1/2},
\end{equation}
for $t_j\in[0,T_{max}]$. The integral above will be approximated by Gauss 
quadrature over LGL-points. 

For Pad\'e based discretization (NP) of the DtN maps, it is possible to carry out
evolution over larger number of time-steps. The Table~\ref{tab:ee-params} lists the 
numerical parameters for the NP20 and NP50 methods. The CQ based methods cannot 
be tested with the same parameters on account of the increasing memory 
requirements at each time-step. Therefore, we introduce a separate test-case 
for the sake of comparison with CQ methods by restricting $N_t=1000+1$ and
$T_{max}=2$ while other parameters are same as that of Table~\ref{tab:ee-params}. 

The numerical results for the evolution error on $\Omega_i$ corresponding to the 
chirped-Gaussian and the Hermite-Gaussian profiles are shown in Fig.~\ref{fig:wpcg2d-ee} 
and Fig.~\ref{fig:wphg2d-ee}, respectively. It can be seen that the diagonal Pad\'e 
approximants based method, namely, NP20 and NP50, perform equally well as compared 
to that of the convolution quadrature based method. Moreover, the TR methods perform
better than the BDF1 methods which is clear from the error peaks in figures~\ref{fig:wpcg2d-ee} 
and~\ref{fig:wphg2d-ee}. The accuracy deteriorates with increase in the value of $c$.
For faster moving profiles corresponding to IA and IIA type solutions, NP-20 becomes 
slightly less accurate compared to CQ and NP-50. There is no significant difference in 
the results for evolution error on $\Omega_i$ corresponding to the chirped-Gaussian 
profile and the Hermite-Gaussian profile.  

The numerical results for the evolution error on $\Omega_i$ particularly for NP methods
corresponding to the chirped-Gaussian and Hermite-Gaussian profiles are shown in 
Fig.~\ref{fig:wpcg2d-ee-np} and Fig.~\ref{fig:wphg2d-ee-np}, respectively. 
Note that the results for NP are superior than that of CQ on account of the
reduced size ($\Delta t$) of the time-step which was possible because increasing the
number of time-steps did not significantly increase the computation time as well
as memory requirements unlike
that of CQ\footnote{Attempts to lower the step-size for the same maximum
duration of time in case of CQ failed due to out-of-memory errors. We note that
these problems may severely limit the practical utility of the CQ method.}.

\begin{figure}[!h]
\begin{center}
\includegraphics[scale=0.95]{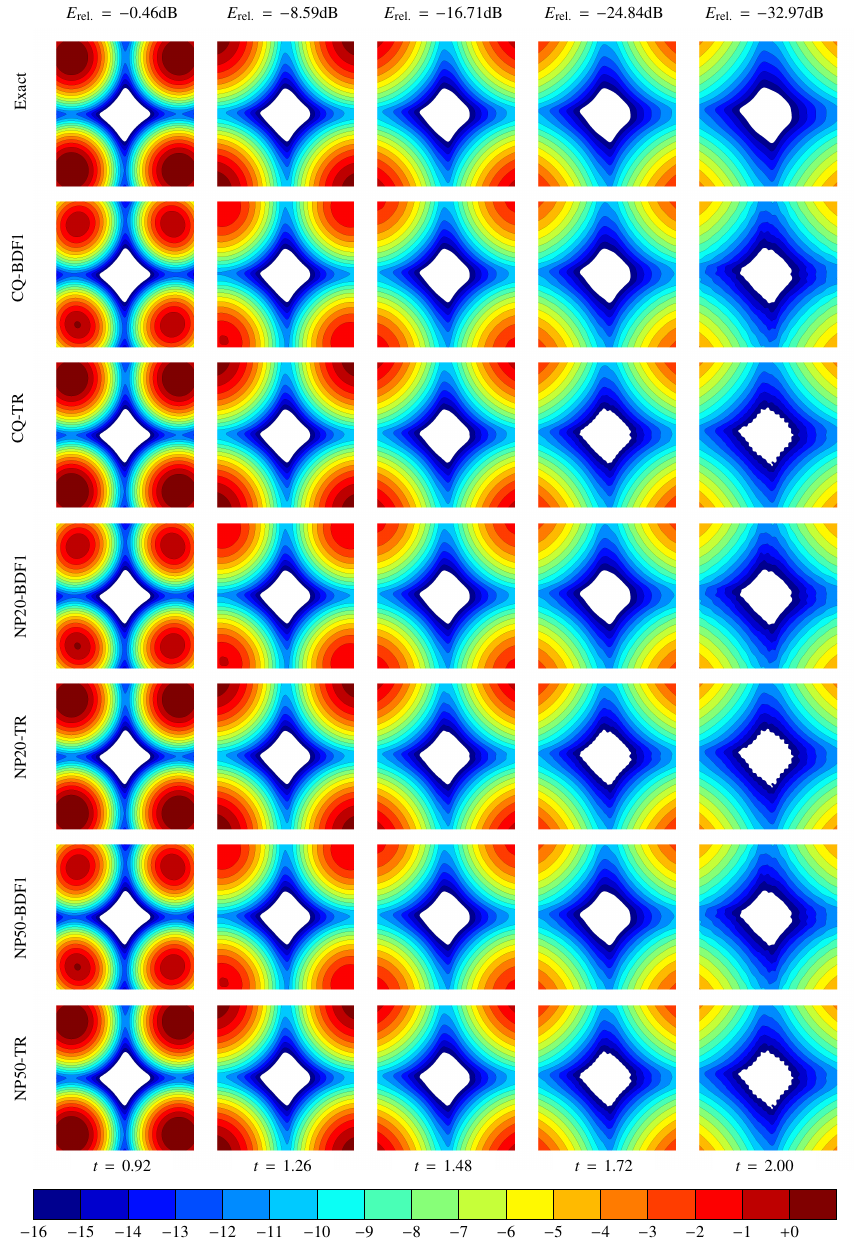}
\end{center}
\caption{\label{fig:wpcg2d_ee_typeIIB_c_12_contour}The figure shows contour plots 
of $4\log_{10}|u(\vv{x},t)|$ with various TBCs for the IIB-type
chirped-Gaussian profile (see Table~\ref{tab:cg2d}) for $c=12$ at different instants of time.
The numerical parameters and the labels are described in Sec.~\ref{sec:tests-ee}.}
\end{figure}

\begin{figure}[!h]
\begin{center}
\includegraphics[scale=0.95]{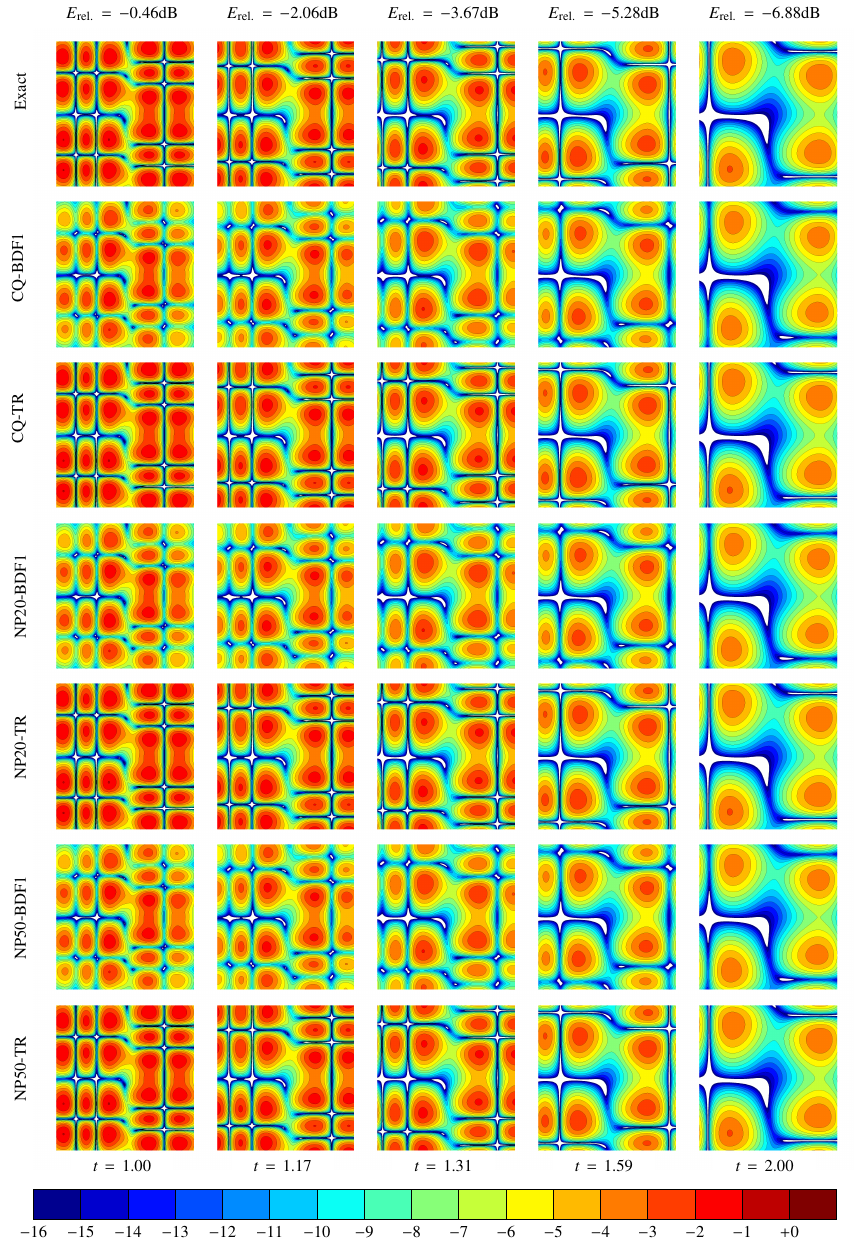}
\end{center}
\caption{\label{fig:wphg2d_ee_typeIIB_c_08_contour}The figure shows contour plots 
of $8\log_{10}|u(\vv{x},t)|$ with various TBCs for the IIB-type
Hermite-Gaussian profile (see Table~\ref{tab:hg2d}) for $c=8$ at different instants
of time. The numerical parameters and the labels are described in Sec.~\ref{sec:tests-ee}.}
\end{figure}

\def\arraystretch{2}
\setlength{\tabcolsep}{1mm}
\begin{table}[htb]
\centering
\caption{\label{tab:ee-params}Numerical parameters for studying the evolution 
error for the NP methods}
\begin{tabular}{m{80mm}m{50mm}}\hline
Computational domain ($\Omega_i$) & $(-10,10)\times(-10,10)$\\\hline
Maximum time ($T_{max}$)          & $5$\\\hline
No. of time-steps ($N_t$)         & $5000+1$\\\hline
Time-step ($\Delta t$)            & $10^{-3}=T_{max}/(N_t-1)$\\\hline
\>Number of LGL-points ($(N+1)\times (N+1)$) & $200\times 200$\\\hline
\end{tabular}
\end{table}
In order to visualize the reflections at the boundary, we resort to contour 
plots depicting the quantity $f_{\text{mag.}}\log_{10}|u(\vv{x},t)|$ where the
scalar quantity $f_{\text{mag.}}$ is adjusted to resolve the contour lines
properly. The specific snapshots correspond to the equispaced energy levels in
the logarithmic scale that is remaining in the computational domain. We restrict 
ourselves to presenting the contour plots only for the 
case of type IIB class of solutions which are directed towards the corners of 
the domain on account of being most complex among four different class of
solutions. The contour plots depicting
corresponding to the chirped-Gaussian ($f_{\text{mag.}}=4$) and the 
Hermite-Gaussian ($f_{\text{mag.}}=8$) profiles are shown in 
Fig.~\ref{fig:wpcg2d_ee_typeIIB_c_12_contour} and 
Fig.~\ref{fig:wphg2d_ee_typeIIB_c_08_contour}, respectively. The top row
corresponds to the exact solution for reference. From these plots, we conclude
that the results are similar for each of the types of solutions. There are 
minimal visible reflections in the contour plots proving the accuracy of the 
developed numerical schemes to implement the TBCs. Moreover, the TR based methods
better resolve the contour levels in comparison to the BDF1 based methods.

\def\arraystretch{2}
\setlength{\tabcolsep}{1mm}
\begin{table}[!h]
\centering
\caption{\label{tab:ca-params}Numerical parameters for studying the convergence 
error for the NP methods}
\begin{tabular}{m{80mm}m{50mm}}\hline
Computational domain ($\Omega_i$) & $(-10,10)\times(-10,10)$\\\hline
Maximum time ($T_{max}$)          & $5$\\\hline
Set of no. of time-steps ($\field{N}_t$)& $\{2^8,2^9,\ldots,2^{18}\}$\\\hline
Time-step               & $\{\Delta t=T_{max}/(N_t-1),\;N_t\in\field{N}_t\}$\\\hline
\>Number of LGL-points ($(N+1)\times (N+1)$) & $200\times 200$\\\hline
\end{tabular}
\end{table}
\subsection{Tests for convergence}\label{sec:tests-ca}
In this section, we analyze the convergence behaviour of the numerical schemes
for the IBVP in~\eqref{eq:2D-SE-CT} where the initial condition corresponds to the 
exact solutions described in Sec.~\ref{sec:exact-solution}. The numerical
solution is labelled according to that of the boundary maps described in
Sec.~\ref{sec:tests-mt}.

The error used to study the convergence behaviour is quantified by 
the maximum relative $\fs{L}^2(\Omega_i)$-norm: 
\begin{equation}\label{eq:max-error-ibvp}
e=\max\left\{e(t_j)|\;j=0,1,\ldots,N_t-1\right\},
\end{equation}
for $t_j\in[0,T_{max}]$.
For Pad\'e based discretization (NP) of the DtN maps, it is possible to analyze
the convergence behaviour over large number of time-steps. The 
Table~\ref{tab:ca-params} lists the numerical parameters for the NP20 and NP50 methods. 
The CQ based methods cannot be tested with the same parameters on account of the
increasing memory requirements at each time-step.
\begin{figure}[!htb]
\begin{center}
\includegraphics[width=\textwidth]{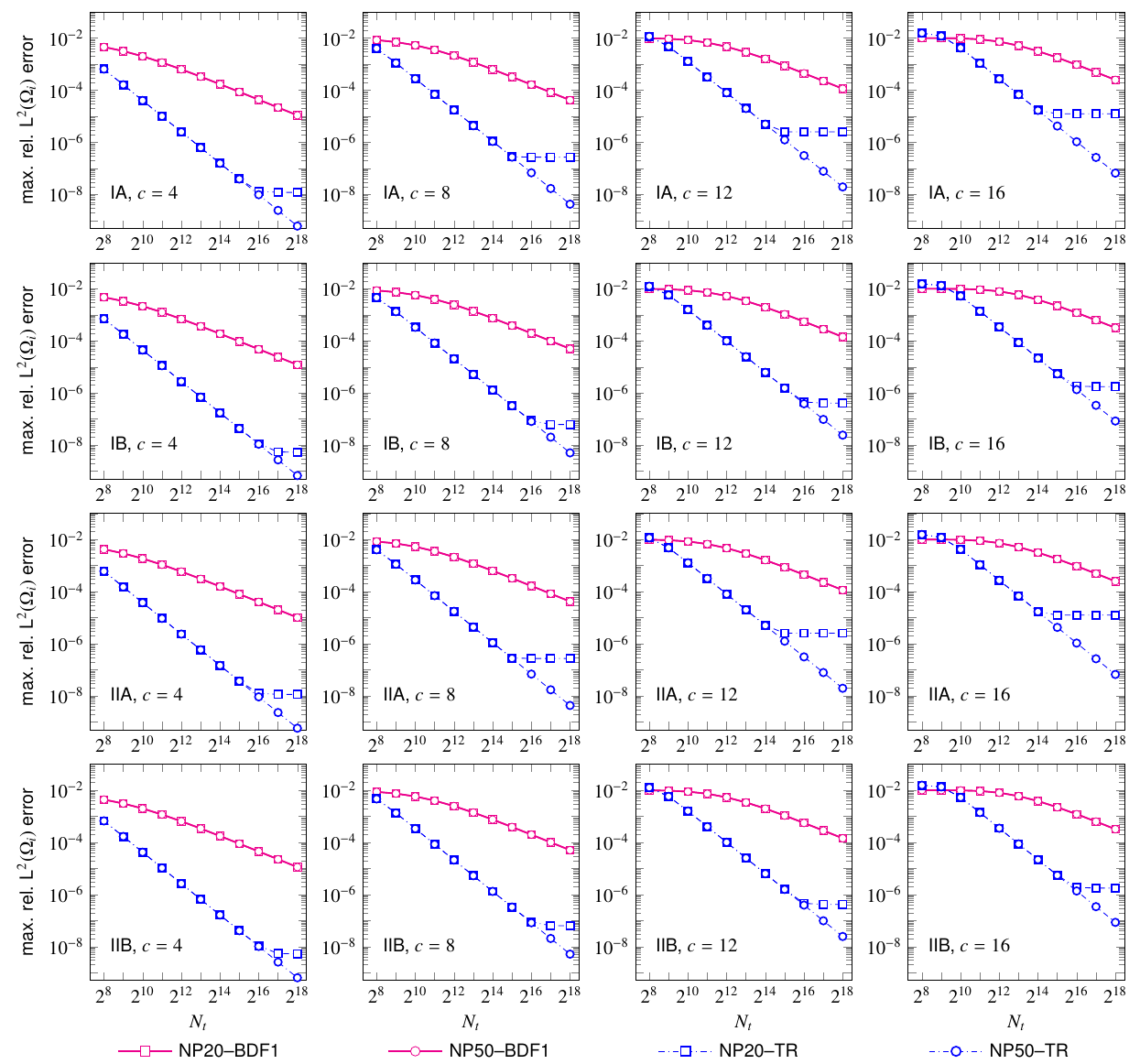}
\end{center}
\caption{\label{fig:wpcg2d-ca}The figure depicts the convergence behaviour of the
numerical solution of the IBVP~\eqref{eq:2D-SE-CT} with the TBCs developed using
novel Pad\'e approach for the chirped-Gaussian profile with different values of the speed `c' 
(see Table~\ref{tab:cg2d}). The numerical parameters and the labels are described in 
Sec.~\ref{sec:tests-ca} where the error is quantified by~\eqref{eq:max-error-ibvp}.}
\end{figure}
\begin{figure}[!h]
\begin{center}
\includegraphics[width=\textwidth]{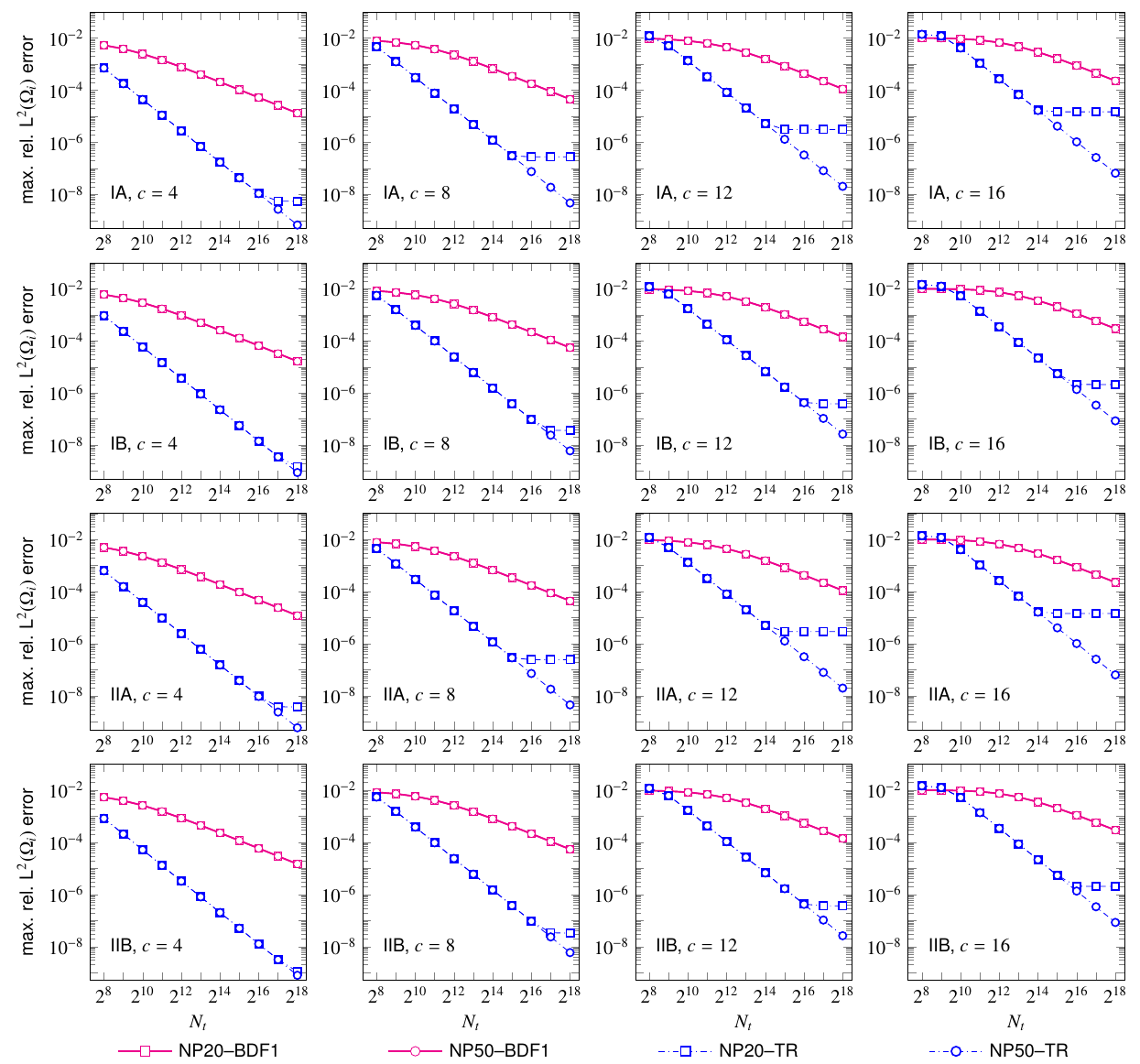}
\end{center}
\caption{\label{fig:wphg2d-ca}The figure depicts the convergence behaviour of the
numerical solution of the IBVP~\eqref{eq:2D-SE-CT} with the TBCs developed using
novel Pad\'e approach for the Hermite-Gaussian profile with different values of the speed `c' 
(see Table~\ref{tab:hg2d}). The numerical parameters and the labels are described in 
Sec.~\ref{sec:tests-ca} where the error is quantified by~\eqref{eq:max-error-ibvp}.}
\end{figure}

The numerical results for the convergence behaviour on $\Omega_i$ corresponding to the 
chirped-Gaussian and the Hermite-Gaussian profiles are shown in Fig.~\ref{fig:wpcg2d-ca} 
and Fig.~\ref{fig:wphg2d-ca}, respectively. It can be seen that the diagonal
Pad\'e approximant based method, namely, NP20 and NP50, shows stable behaviour for the
parameters specified in Table~\ref{tab:ca-params}. In the log-log scale, we can
clearly identify the error curve to be a straight line before it plateaus. The
orders can be recovered from the slope which we found consistent with the order
of the underlying one-step method. The TR methods perform better than the BDF1 
methods which is obvious from the slope of the error curves in 
figures~\ref{fig:wpcg2d-ca} and~\ref{fig:wphg2d-ca}. We note that the error 
background increases with increase in the value of $c$. This 
is on account of the increased oscillatory behaviour of the solution in $t$ as
$c$ increases and a method with higher order than $2$ would be required to
achieve the same accuracy as that of the slow moving profiles keeping the step-size
same. Note that NP20 methods starts
plateauing after a certain number of time-steps unlike NP50. This can be attributed 
to the approximate nature of the TBCs with a lower order of diagonal approximants 
chosen. Finally, let us observe that there is no significant difference in the 
results for convergence on $\Omega_i$ corresponding to the 
chirped-Gaussian profile and the Hermite-Gaussian profile.  


\section{Conclusion}\label{sec:conclusion}
In this paper, we presented a new effectively local approximation for the
numerical realization of the TBCs on a rectangular computational domain for 
the free Schr\"{o}dinger equation in $\field{R}^2$. Our starting point was a 
numerically tractable formulation of the transparent boundary operator, 
$\sqrt{\partial_t-i\triangle_{\Gamma}}$, where $\triangle_{\Gamma}$ is the 
Laplace-Beltrami operator~\cite{FP2011,V2019}, which is nonlocal in time 
as well as space. This formulation of the TBCs introduces 
an auxiliary function which happens to satisfy the 
one-dimensional Schr\"{o}dinger equation along the boundary segments of the 
computational domain. It addressed the nonlocality in space, however, the 
nonlocality in time remained a serious bottleneck on account of 
the growing computational complexity with number of time-steps which becomes 
prohibitive after a point. To overcome this challenge, we developed a
novel Pad\'e approximant based rational approximation for the $1/2$-order 
time-fractional derivative. This approach circumvents the need to store the
entire history of the field on the boundary segments and allows us to achieve 
a computational complexity that remains independent of the number of time-steps.
Let us emphasize that our approach is significantly different from that of 
Menza~\cite{Menza1997} which fails at the corners of the rectangular domain.
Further, we also present a convolution quadrature based scheme which 
happens to be computationally expensive but it is well-known as a golden 
standard for 1D problems. As far as time marching is concerned, we note that
one-step methods are best suited for our purpose and any higher-order method can
be readily used to implement our algorithm. These aspects will be further 
explored in a future publication.

For the spatial discretization, we use a Legendre-Galerkin spectral method with
a new boundary adapted basis which leads to a banded linear system for the 1D
problems on the boundary segments as well as for the interior problem. The
time-discrete version of the boundary condition are formulated as a Robin-type
condition in a manner reminiscent of the 1D problem~\cite{CiCP2008}. A compatible 
boundary-lifting procedure to homogenize the Robin-type boundary conditions 
is also presented. In order to demonstrate the effectiveness 
of the methods developed, we carried out numerical tests to study the behaviour 
of evolution of error as well as to verify the order of convergence empirically.
If the number of spatial grid points in every dimension is $\bigO{N}$, then the
computational complexity of our novel 
Pad\'e method for the DtN map scales as $\bigO{MN}$ where $M$ is the order of 
the diagonal Pad\'e approximant. Let us emphasize that this complexity is
independent of the number of time-steps. The overall complexity thus becomes dominated 
by the cost of the linear solver for interior problem which is given by 
$\Theta_2(N_{\text{dim}})=\bigO{N^3}$. A theoretical proof of the stability 
and convergence of the numerical scheme is beyond the scope of the paper. These 
issues are deferred to a future publication.

The ideas presented in this paper can be readily extended to three-dimensions. 
Feschenko and Popov have already shown that the ideas in 2D~\cite{FP2011} can 
be extended to 3D~\cite{FP2013}. The formulation of the operator 
$\sqrt{\partial_t-i\triangle_{\Gamma}}$ on the boundary of a cuboidal computational domain 
can be obtained along the same lines as that presented in~\cite{V2019}. For all
the cases discussed above, it remains to determine how competitive are PMLs.
Such a comparison should however be made for a more general type of Schr\"odinger
equation which is identical to the free Schr\"odinger on the
exterior domain $\Omega_e=\field{R}^d\setminus\ovl{\Omega}_i,\;d=2,3$.


\section*{Acknowledgements}
The first author thanks CSIR India (grant no. 09/086(1431)/2019-EMR-I) for 
providing the financial assistance.  

\bibliographystyle{elsarticle-num} 

\begin{thebibliography}{10}
\expandafter\ifx\csname url\endcsname\relax
  \def\url#1{\texttt{#1}}\fi
\expandafter\ifx\csname urlprefix\endcsname\relax\def\urlprefix{URL }\fi
\expandafter\ifx\csname href\endcsname\relax
  \def\href#1#2{#2} \def\path#1{#1}\fi

\bibitem{AK2003}
Y.~S. Kivshar, G.~P. Agrwal, Optical Solitons: From Fibers to Photonic
  Crystals, 1st Edition, Academic Press, San Diego, California, 2003.

\bibitem{LM1998}
D.~Lee, S.~T. McDaniel, Ocean Acoustic Propagation by Finite Difference
  Methods, Modern Applied Mathematics and Computer Science, Pergmon Press, New
  York, 1988.

\bibitem{CiCP2008}
X.~Antoine, A.~Arnold, C.~Besse, M.~Ehrhardt, A.~Sch\"{a}dle,
  \href{https://hal.archives-ouvertes.fr/hal-00347884}{A review of transparent
  and artificial boundary conditions techniques for linear and nonlinear
  {S}chr\"{o}dinger equations}, Comm. Comput. Phys. 4~(4) (2008) 729--796.
\newline\urlprefix\url{https://hal.archives-ouvertes.fr/hal-00347884}

\bibitem{FP2011}
R.~M. Feshchenko, A.~V. Popov, Exact transparent boundary condition for the
  parabolic equation in a rectangular computational domain, J. Opt. Soc. Am. A
  28~(3) (2011) 373--380.
\newblock \href {https://doi.org/10.1364/JOSAA.28.000373}
  {\path{doi:10.1364/JOSAA.28.000373}}.

\bibitem{V2019}
V.~Vaibhav, On the nonreflecting boundary operators for the general two
  dimensional {S}chr\"odinger equation, J.~Math.~Phys. 60~(1) (2019) 011509.
\newblock \href {https://doi.org/10.1063/1.5030875}
  {\path{doi:10.1063/1.5030875}}.

\bibitem{Menza1996}
L.~D. Menza, Absorbing boundary conditions on a hypersurface for the
  {S}chr\"{o}dinger equation in a half-space, Appl. Math. Lett. 9~(4) (1996)
  55--59.
\newblock \href {https://doi.org/10.1016/0893-9659(96)00051-1}
  {\path{doi:10.1016/0893-9659(96)00051-1}}.

\bibitem{Menza1997}
L.~D. Menza, Transparent and absorbing boundary conditions for the
  {S}chr\"{o}dinger equation in a bounded domain, Numer. Funct. Anal. Optim.
  18~(7) (1997) 759--775.
\newblock \href {https://doi.org/10.1080/01630569708816790}
  {\path{doi:10.1080/01630569708816790}}.

\bibitem{S2002}
A.~Sch\"adle, Non-reflecting boundary conditions for the two-dimensional
  {S}chr\"odinger equation, Wave Motion 35~(2) (2002) 181--188.
\newblock \href {https://doi.org/10.1016/S0165-2125(01)00098-1}
  {\path{doi:10.1016/S0165-2125(01)00098-1}}.

\bibitem{HH2004}
H.~Han, Z.~Huang, \href{https://projecteuclid.org/euclid.cms/1250880210}{Exact
  artificial boundary conditions for {S}chr\"odinger equation in
  $\mathbb{R}^2$}, Comm. Math. Sci. 2~(1) (2004) 79--94.
\newline\urlprefix\url{https://projecteuclid.org/euclid.cms/1250880210}

\bibitem{JPA2018}
S.~Ji, Y.~Yang, G.~Pang, X.~Antoine, Accurate artificial boundary conditions
  for the semi-discretized linear {S}chr{\"o}dinger and heat equations on
  rectangular domains, Comput. Phys. Commun. 222 (2018) 84--93.
\newblock \href {https://doi.org/10.1016/j.cpc.2017.09.019}
  {\path{doi:10.1016/j.cpc.2017.09.019}}.

\bibitem{AB2001}
{\color{myrem}
X.~Antoine, C.~Besse, Construction, structure and asymptotic approximations of
  a microdifferential transparent boundary condition for the linear
  {S}chr\"{o}dinger equation, J. Math. Pures Appl. 80~(7) (2001) 701--738.
\newblock \href {https://doi.org/10.1016/S0021-7824(01)01213-2}
  {\path{doi:10.1016/S0021-7824(01)01213-2}}.
}
\bibitem{S2005}
{\color{myrem}
J.~Szeftel, Design of absorbing boundary conditions for {S}chr\"odinger
  equations in $\mathbb{R}^d$, SIAM Journal on Numerical Analysis 42~(4) (2005)
  1527--1551.
}
\bibitem{ABK2012}
X.~Antoine, C.~Besse, P.~Klein, Absorbing boundary conditions for the
  two-dimensional {S}chr\"odinger equation with an exterior potential. {P}art
  {I}: {C}onstruction and a priori estimates, Math. Models Methods Appl. Sci.
  22~(10) (2012) 1250026.
\newblock \href {https://doi.org/10.1142/S0218202512500261}
  {\path{doi:10.1142/S0218202512500261}}.

\bibitem{ABK2013}
X.~Antoine, C.~Besse, P.~Klein, Absorbing boundary conditions for the
  two-dimensional {S}chr\"odinger equation with an exterior potential. {P}art
  {II}: {D}iscretization and numerical results, Numer. Math. 125~(2) (2013)
  191--223.
\newblock \href {https://doi.org/10.1007/S00211-013-0542-8}
  {\path{doi:10.1007/S00211-013-0542-8}}.

\bibitem{V2014}
V.~Vaibhav, Microlocal approach towards construction of nonreflecting boundary
  conditions, J. Comput. Phys. 272 (2014) 588--607.
\newblock \href {https://doi.org/10.1016/j.jcp.2014.04.050}
  {\path{doi:10.1016/j.jcp.2014.04.050}}.

\bibitem{S1994}
J.~Shen, Efficient spectral-{G}alerkin method {I}. {D}irect solvers of second-
  and fourth-order equations using {L}egendre polynomials, SIAM J. Sci. Comput.
  15~(6) (1994) 1489--1505.
\newblock \href {https://doi.org/10.1137/0915089} {\path{doi:10.1137/0915089}}.

\bibitem{Lubich1986}
C.~Lubich, Discretized fractional calculus, SIAM J. Math. Anal. 17~(3) (1986)
  704--719.
\newblock \href {https://doi.org/10.1137/0517050} {\path{doi:10.1137/0517050}}.

\bibitem{BP1991}
V.~Baskakov, A.~Popov, Implementation of transparent boundaries for numerical
  solution of the {S}chr\"odinger equation, Wave Motion 14~(2) (1991) 123--128.
\newblock \href {https://doi.org/10.1016/0165-2125(91)90053-Q}
  {\path{doi:10.1016/0165-2125(91)90053-Q}}.

\bibitem{Mayfield1989}
B.~Mayfield, Nonlocal boundary conditions for the {S}chr\"{o}dinger equation,
  Ph.D. thesis, University of Rhodes Island, Providence, RI (1989).

\bibitem{Zheng2007}
C.~Zheng, A perfectly matched layer approach to the nonlinear {S}chr\"odinger
  wave equations, J. Comput. Phys. 227~(1) (2007) 537--556.
\newblock \href {https://doi.org/10.1016/j.jcp.2007.08.004}
  {\path{doi:10.1016/j.jcp.2007.08.004}}.

\bibitem{AGT2020}
X.~Antoine, C.~Geuzaine, Q.~Tang, Perfectly matched layer for computing the
  dynamics of nonlinear {S}chr\"{o}dinger equations by pseudospectral methods.
  {A}pplication to rotating {B}ose-{E}instein condensates, Communications in
  Nonlinear Science and Numerical Simulation 90 (2020) 105406.
\newblock \href {https://doi.org/10.1016/j.cnsns.2020.105406}
  {\path{doi:10.1016/j.cnsns.2020.105406}}.

\bibitem{MR1993}
K.~Miller, B.~Ross, An Introduction to the Fractional Calculus and Fractional
  Differential Equations, 1st Edition, John Wiley and Sons Inc., 1993.

\bibitem{M2013}
{\color{myrem}
J.~Mennemann, A.~J\"ungel, H.~Kosina, Transient {S}chr{\"o}dinger--{P}oisson
  simulations of a high-frequency resonant tunneling diode oscillator, Journal
  of Computational Physics 239 (2013) 187--205.
\newblock \href {https://doi.org/10.1016/j.jcp.2012.12.009}
  {\path{doi:10.1016/j.jcp.2012.12.009}}.
}
\bibitem{V2011}
V.~Vaibhav, Artificial boundary conditions for certain evolution {PDE}s with
  cubic nonlinearity for noncompactly supported initial data, J. Comput. Phys.
  230~(8) (2011) 3205--3229.
\newblock \href {https://doi.org/10.1016/j.jcp.2011.01.024}
  {\path{doi:10.1016/j.jcp.2011.01.024}}.

\bibitem{L1985}
{\color{myrem}
E.~Lindmann, Free-space boundary conditions for time dependant wave equation,
  J. Compt. Phys. 18~(1) (1975) 16--78.
\newblock \href {https://doi.org/10.1016/0021-9991(75)90102-3}
  {\path{doi:10.1016/0021-9991(75)90102-3}}.
}
\bibitem{SV2023}
S.~Yadav, V.~Vaibhav, Nonreflecting boundary condition for the free
  {S}chr{\"o}dinger equation in 2d, in: 2023 Photonics and Electromagnetics
  Research Symposium (PIERS), 2023, pp. 328--337.
\newblock \href {https://doi.org/10.1109/PIERS59004.2023.10221299}
  {\path{doi:10.1109/PIERS59004.2023.10221299}}.

\bibitem{Higham2002}
N.~J. Higham, Accuracy and Stability of Numerical Algorithms, 2nd Edition,
  Society for Industrial and Applied Mathematics, Philadelphia, 2002.
\newblock \href {https://doi.org/10.1137/1.9780898718027}
  {\path{doi:10.1137/1.9780898718027}}.

\bibitem{Shen2011}
J.~Shen, T.~Tang, L.~L. Wang, Spectral Methods: Algorithms, Analysis and
  Applications, 1st Edition, Springer Berlin Heidelberg, 2011.

\bibitem{FP2013}
R.~M. Feshchenko, A.~V. Popov, Exact transparent boundary condition for the
  three-dimensional {S}chr\"odinger equation in a rectangular cuboid
  computational domain, Phys. Rev. E 88 (2013) 053308.
\newblock \href {https://doi.org/10.1103/PhysRevE.88.053308}
  {\path{doi:10.1103/PhysRevE.88.053308}}.

\end{thebibliography}
\providecommand{\noopsort}[1]{}\providecommand{\singleletter}[1]{#1}%

\clearpage

\appendix 

\section{Fractional Operators}\label{app:frac-op}
\begin{defn}[Riemann-Liouville fractional integrals]\label{def:RL-FI}
Let $\alpha>0$ and $f(t)$ be a piecewise continuous function on 
$(0,\infty)$, locally integrable on any finite subinterval of 
$[0,\infty)$. For $t>0$, the Riemann-Liouville fractional integral 
of order $\alpha$, denoted by $\partial_t^{-\alpha}$, is defined by
\begin{equation}\label{eq:fracintegral}
\partial_t^{-\alpha} f(t)=\frac{1}{\Gamma(\alpha)}
\int^{t}_0(t-s)^{\alpha-1}f(s)ds.
\end{equation}
where $\Gamma$ denotes the Euler's Gamma function.
\end{defn}

\begin{theorem}[Law of exponents]
Let $f(t)\in\fs{C}^{\infty}([0,\infty))$ and let $\mu,\nu>0$, then the composition 
of two fractional integrals is given by
\begin{equation*}
\partial_t^{-\nu}[\partial_t^{-\mu}f(t)]=\partial_t^{-(\mu+\nu)}f(t)
=\partial_t^{-\mu}[\partial_t^{-\nu}f(t)],\quad\forall t\geq0.
\end{equation*}
\end{theorem}

\begin{defn}[Riemann-Liouville fractional derivatives]
Let $f(t)\in\fs{C}^{\infty}([0,\infty))$ and let $\mu>0$. If $p$ is the smallest 
integer greater than $\mu$, then the fractional derivative of $f$ of order 
$\mu$ is defined by
\begin{equation}\label{def:fdc}
\partial_t^{\mu}f(t)=\partial^p_t[\partial_t^{-\nu}f(t)],\quad t>0,
\end{equation}
where $\nu=p-\mu>0$. For this fractional derivative to be defined at $t=0$, 
the function $f(t)$ must have first $p-1$ derivatives equal to zero at $t=0$.
\end{defn}

\section{Boundary Adapted Basis Representation: Useful Results and Identities}
\label{app:ip}

In order to leverage the discrete Legendre transforms over 
Legendre-Gauss-Lobatto nodes, we choose to favour 
storage of space dependent variables in terms of their Legendre coefficients. 
As a result of this choice, all operations involving the boundary adapted basis 
must be implemented via the Legendre transform.

Let us recall the definitions: $L_n(y)$ denote the 
Legendre polynomial of degree $n$ and the index set $\{0,1,\ldots,N-2\}$ is 
denoted by $\field{J}$. The boundary adapted basis polynomials are given by
\begin{equation}
\phi_{p}(y)=L_p(y)+b_pL_{p+2}(y),\quad 
b_p=-\frac{\kappa+\frac{1}{2}p(p+1)}{\kappa+\frac{1}{2}(p+2)(p+3)},\quad 
p\in\field{J},
\end{equation}
where $\kappa\in\field{C}$. Let us briefly discuss the parameters that determine
this quantity: We recall from 
Sec.~\ref{sec:linear-sys-1D} that $\beta=J^{-2}$ where $J=\partial_y x$ is the Jacobian
corresponding to the transformation from the reference interval $\field{I}=[-1,1]$ to the
actual domain $[x_-,x_+]$. We suppress the subscript which characterizes the
axis along which the domain lies for the sake of brevity. The time-stepping 
schemes discussed in Sec.~\ref{sec:numerical-implementation} dictate that 
either $\rho=1/\Delta t$ (for BDF1) or $\rho=2/\Delta t$ (for TR).
The form of $\kappa$ follows from the Robin-type formulation of the discretized
version of the DtN-maps as discussed in Sec.~\ref{sec:dtn-maps}.
For CQ methods $\kappa=\alpha=\sqrt{\rho/\beta}\exp(-i\pi/4)$ while 
for the NP methods, we have $\kappa=\alpha\varpi$ where $\varpi$ is defined in~\eqref{eq:varpi-defn}.

Any arbitrary function in 
$\fs{L}^2(\field{I})$ can be expanded in this basis as follows:
\begin{equation}
f(y)=\sum_{p=0}^{N-2}\hat{f}_{p}\phi_{p}(y)
+\eta_{0}L_0(y)+\eta_{1}L_1(y)
=\sum_{p'=0}^N\tilde{f}_{p'}L_{p'}(y).
\end{equation}
By virtue of linear independence of Legendre polynomials, we can equate their 
coefficients on both sides of the equality so that
\begin{equation}
\left\{\begin{aligned}
&\tilde{f}_{0}=\hat{f}_{0}+\eta_{0}, \tilde{f}_{1}=\hat{f}_{1}+\eta_{1},\\
&\tilde{f}_{p}=\hat{f}_{p}+b_{p-2}\hat{f}_{p-2},\quad p=2,3,\ldots N-2,\\
&\tilde{f}_{N-1}=b_{N-3}\hat{f}_{N-3},\tilde{f}_{N}=b_{N-2}\hat{f}_{N-2}.
\end{aligned}\right.
\end{equation}
The linear system can be written as
\begin{equation}
\wtilde{\vv{f}}=B\what{\vv{f}}+\eta_0\vv{e}_0+\eta_1\vv{e}_1
\end{equation}
where
\begin{equation}
B=
\begin{pmatrix}
1   &    &      &      &      &\\
0   &  1 &      &      &      &\\
b_0 &  0 &  1   &      &      &\\
    & b_1&  0   &     1&      &\\
    &    &\ddots&\ddots&\ddots&\\
    &    &      &b_{N-4}&    0&1\\
    &    &      &       &b_{N-3}&0\\
    &    &      &       &     &b_{N-2}
\end{pmatrix}\in\field{C}^{(N+1)\times (N-1)},
\end{equation}
Let the normalization factors be denoted by 
$\gamma_k=\|L_k\|^2_2=2/(2k+1),\,k\in\field{N}_0,$ and introduce the diagonal matrix 
$\Gamma=\diag(\gamma_0,\gamma_1,\ldots,\gamma_{N})$, then inner products of the form 
$g_{p}=\left(f,\phi_p\right)_{\field{I}},\,p\in\field{J},$ can be understood as
\begin{equation}
g_{p}=\left(f,\phi_p\right)_{\field{I}}
=\left[\gamma_p\tilde{f}_{p}+b_p\gamma_{p+2}\tilde{f}_{p+2}\right],
\quad p\in\field{J}.
\end{equation}
The relationship above can be compactly stated as $\vv{g}=(g_{0},g_1,\ldots,g_{N-2})^{\tp}
=Q\Gamma\wtilde{\vv{f}}$ where $Q=B^{\tp}$ is 
referred to as the \emph{quadrature matrix}.

Let the index set $\{0,1,\ldots,N_1-2\}$ be denoted by $\field{J}_1$ 
and the index set $\{0,1,\ldots,N_2-2\}$ be denoted by $\field{J}_2$.
The boundary adapted basis polynomials are given by
\begin{equation}
\phi^{(j)}_{p_j}(y_j)=L_{p_j}(y_j)+b^{(j)}_{p_j}L_{p_j+2}(y_j),
\quad p_j\in\field{J}_j,\;j=1,2.
\end{equation}
where $\{b^{(j)}_{p_j}\}$ is a sequence.
Any arbitrary function in $\fs{L}^2(\Omega_i^{\text{ref}})$ can be expanded in this basis 
as follows:
\begin{equation}
\begin{split}
f(\vv{y})
&=\sum_{p_1\in\field{J}_1}\sum_{p_2\in\field{J}_2}\hat{f}_{p_1,p_2}\phi^{(1)}_{p_1}(y_1)\phi^{(2)}_{p_2}(y_2)\\
&\quad+\sum_{p_1\in\field{J}_1}\left[\zeta_{p_1,0}\phi^{(1)}_{p_1}(y_1)L_0(y_2)
+\zeta_{p_1,1}\phi^{(1)}_{p_1}(y_1)L_1(y_2)\right]\\
&\quad+\sum_{p_2\in\field{J}_2}\left[\eta_{0,p_2}L_0(y_1)\phi^{(2)}_{p_2}(y_2)
+\eta_{1,p_2}L_1(y_1)\phi^{(2)}_{p_2}(y_2)\right]\\
&\quad+\left[\kappa_{0,0}L_0(y_1)L_0(y_2)
+\kappa_{0,1}L_0(y_1)L_1(y_2)+\kappa_{1,0}L_1(y_1)L_0(y_2)+\kappa_{1,1}L_1(y_1)L_1(y_2)\right]\\
&=\sum_{p_1=0}^{N_1}\sum_{p_2=0}^{N_2}\tilde{f}_{p_1,p_2}L_{p_1}(y_1)L_{p_2}(y_2).
\end{split}
\end{equation}
For $p_1=2,3,\ldots N_1-2$ and $p_2=2,3,\ldots N_2-2$, 
equating the coefficients on both sides of equality so that
\begin{equation}
\tilde{f}_{p_1,p_2}=\hat{f}_{p_1,p_2}
+b^{(1)}_{p_1-2}\hat{f}_{p_1-2,p_2}
+b^{(2)}_{p_2-2}\hat{f}_{p_1,p_2-2}
+b^{(1)}_{p_1-2}b^{(2)}_{p_2-2}\hat{f}_{p_1-2,p_2-2},
\end{equation}
The linear system can be written as
\begin{equation}
\begin{split}
\wtilde{F}
&=B_1\what{F}B_2^{\tp}
+B_1\vs{\zeta}_0\otimes\left(\vv{e}^{(2)}_0\right)^{\tp}
+B_1\vs{\zeta}_1\otimes\left(\vv{e}^{(2)}_1\right)^{\tp}
+\vv{e}^{(1)}_0\otimes\vs{\eta}_0B_2^{\tp}
+\vv{e}^{(1)}_1\otimes\vs{\eta}_1B_2^{\tp}\\
&\quad
+\kappa_{0,0}\vv{e}^{(1)}_0\otimes\left(\vv{e}^{(2)}_0\right)^{\tp}
+\kappa_{0,1}\vv{e}^{(1)}_0\otimes\left(\vv{e}^{(2)}_1\right)^{\tp}
+\kappa_{1,0}\vv{e}^{(1)}_1\otimes\left(\vv{e}^{(2)}_0\right)^{\tp}
+\kappa_{1,1}\vv{e}^{(1)}_1\otimes\left(\vv{e}^{(2)}_1\right)^{\tp}
\end{split}
\end{equation}
where
\begin{equation}
B_j=
\begin{pmatrix}
1   &    &      &      &      &\\
0   &  1 &      &      &      &\\
b^{(j)}_0 &  0 &  1   &      &      &\\
    & b^{(j)}_1&  0   &     1&      &\\
    &    &\ddots&\ddots&\ddots&\\
    &    &      &b^{(j)}_{N-4}&    0&1\\
    &    &      &       &b^{(j)}_{N-3}&0\\
    &    &      &       &     &b^{(j)}_{N-2}
\end{pmatrix}\in\field{C}^{(N+1)\times (N-1)},
\end{equation}
The inner products of the form 
$g_{p_1,p_2}=\left(f,\phi^{(1)}_{p_1}\phi^{(2)}_{p_2}\right)_{\Omega_i^{\text{ref}}}$ 
can be understood as
\begin{equation}
\begin{split}
g_{p_1,p_2}=\left(f,\phi^{(1)}_{p_1}\phi^{(2)}_{p_2}\right)_{\Omega_i^{\text{ref}}}
&=\gamma_{p_1}\gamma_{p_2}\tilde{f}_{p_1,p_2}
+b^{(1)}_{p_1}\gamma_{p_1+2}\gamma_{p_2}\tilde{f}_{p_1+2,p_2}\\
&\qquad+b^{(2)}_{p_2}\gamma_{p_1}\gamma_{p_2+2}\tilde{f}_{p_1,p_2+2}
+b^{(1)}_{p_1}b^{(2)}_{p_2}\gamma_{p_1+2}\gamma_{p_2+2}\tilde{f}_{p_1+2,p_2+2}.
\end{split}
\end{equation}
The relationship above can be compactly stated as $G=(g_{p_1,p_2})=Q_1\Gamma\wtilde{F}\Gamma Q_2^{\tp}$
where $Q_j=B_j^{\tp}$ are referred to as the \emph{quadrature matrices}.

\end{document}